\documentclass[12pt]{elsarticle_arx}
\usepackage{amssymb, amsmath, verbatim,floatpag}
\usepackage[english]{babel}

\newtheorem{theorem}{Theorem}

\newtheorem{lemma}{Lemma}[section]
\newtheorem{prop}{Proposition}

\newtheorem{theoremA}{Theorem}

\newtheorem{conjA}[theoremA]{Conjecture}

\newproof{prf}{Proof}

\newdefinition{prob}{Problem}

\newdefinition{rem}{Remark}

\renewcommand{\pmod}[1]{\,(\mathrm{mod}\,#1)}

\newcommand{\Aut}{\operatorname{Aut}}
\newcommand{\Out}{\operatorname{Out}}

\newcommand{\prk}{\operatorname{prk}}

\begin{document}

\begin{frontmatter}

\title{On finite groups isospectral to simple classical groups}

\author{A.V. Vasil$'$ev \fnref{l1}}
\ead{vasand@math.nsc.ru}
\address{Sobolev Institute of Mathematics and Novosibirsk State University, Novosibirsk, Russia}
\fntext[l1]{The author was supported by the RFBR (project 13-01-00505) and the
Target program of SB RAS for 2012-2014 (integration project No. 14).}

\begin{abstract}

The \emph{spectrum} $\omega(G)$ of a finite group $G$ is the set of element orders of $G$.
Finite groups $G$ and $H$ are \emph{isospectral} if their spectra coincide. Suppose that $L$
is a simple classical group of sufficiently large dimension (the lower bound
varies for different types of groups but is at most 62) defined over a finite field of characteristic $p$. It is proved
that a finite group $G$ isospectral to $L$ cannot have a nonabelian composition factor which is a
group of Lie type defined over a field of characteristic distinct from~$p$.
Together with a series of previous results this implies that every finite group $G$ isospectral to $L$ is `close'
to $L$. Namely, if $L$ is a linear or unitary group, then $L\leqslant G\leqslant\Aut{L}$, in
particular, there are only finitely many such groups $G$ for given $L$. If $L$ is a symplectic or orthogonal
group, then $G$ has a unique nonabelian composition factor $S$ and, for given $L$, there are at
most 3 variants for $S$ (including $S\simeq L$).

\end{abstract}


\end{frontmatter}

\section*{Introduction}

The \emph{spectrum} $\omega(G)$ of a finite group $G$ is the set of element orders of $G$. Finite groups having the same spectra are said to be \emph{isospectral}.
Recently, the following general result was obtained.

\begin{theoremA} Suppose that $L$ is a finite simple group, and $G$ is a finite group with
$\omega(G)=\omega(L)$ and $|G|=|L|$. Then $G$ is isomorphic to $L$.
\end{theoremA}

The statement was conjectured by W.~Shi \cite{89ShiS} in 1987, while its proof was completed in \cite{VasGrMaz2} (see the latter article for background information and complete list of references). It is worth mentioning that Theorem~A together with~\cite[Corollary~5.2]{08KimLuRC} implies that a finite simple group and a finite group with the same Burnside rings are isomorphic as well.

What happens if we omit the condition $|G|=|L|$ in Theorem~A? Then $G$ is not necessary
isomorphic to~$L$. For example, there are infinitely many groups with the same spectrum as the
spectrum of the alternating permutation group of degree six \cite{91BrShi}. On the other hand, it turns out that for a bulk of finite nonabelian simple groups $L$, a finite group isospectral to $L$ is isomorphic to $L$ or, at least, to a group $G$ with $L\leqslant G\leqslant \Aut{L}$. Investigations on this subject have 30-year history and resulted in more than a hundred papers of numerous authors. We do not intend to give a detailed review, rather prefer to draw an overall picture. In order to do that we formulate the following
conjecture attributed to V.D.~Mazurov:

\begin{conjA} For every finite nonabelian simple group $L$, apart from a finite number of sporadic, alternating and exceptional groups and apart from several series of classical groups of small dimensions, if a finite group $G$ is isospectral to~$L$ then $G$ is an almost simple group with socle isomorphic to~$L$.
\end{conjA}

The conjecture was proved for sporadic groups \cite{MazShi}, for alternating groups \cite{GorAlt}, and very recently for exceptional groups of Lie type \cite{14VasSt.t}. Here we deal with groups isospectral to finite simple classical groups. Observe that there are a lot of results on this topic concerning groups in particular characteristics or dimensions. The peculiarity of our approach is that we sacrifice groups of small dimensions in order to obtain a general result (cf. the theorems below with recent papers \cite{09AleKon.t, 09Kon, 09HeShi,10SheShiZin.t,12HeShi, 12GrLyt.t, 12ForIraAha, 13ForIraAha}). Throughout this paper we use single-letter names for simple classical groups, following \cite{Atlas}, i.\,e., for example, $L_n(q)$ means $PSL_n(q)$, as well as the standard abbreviation $L_n^\varepsilon(q)$, where $\varepsilon\in\{+,-\}$, $L_n^+(q)=L_n(q)$, and $L_n^-(q)=U_n(q)$.

\begin{theorem}\label{t:L&U}
Suppose that $L=L_n^\varepsilon(q)$ is a finite simple linear or unitary group and $n\geqslant45$.
Then a finite group isospectral to~$L$ is isomorphic to a group $G$ with $L\leqslant G\leqslant \Aut{L}$. In
particular, there are only finitely many pairwise non-isomorphic finite groups $G$ with $\omega(G)=\omega(L)$.
\end{theorem}

\begin{theorem}\label{t:S&O}
Suppose that $L$ is a finite simple symplectic or orthogonal group, where $n\geqslant28$ for
$L\in\{S_{2n}(q), O_{2n+1}(q)\}$, $n\geqslant31$ for $L=O_{2n}^+(q)$, and $n\geqslant30$ for
$L=O_{2n}^-(q)$. If $G$ is a finite group with $\omega(G)=\omega(L)$, then $G$ has a unique
nonabelian composition factor $S$, and one of the following holds:

\emph{(i)} $L\simeq S;$

\emph{(ii)} $L\in\{S_{2n}(q), O_{2n+1}(q)\}$ and $S\in\{O_{2n+1}(q), O_{2n}^-(q)\};$

\emph{(iii)} $n$ is even, $L=O_{2n}^+(q)$, and $S\in\{S_{2n-2}(q), O_{2n-1}(q)\}.$

In particular, there exist at most $3$ possibilities for~$S$ for given~$L$.
\end{theorem}

\begin{rem} It seems very likely that the conclusion of Theorem~\ref{t:L&U} is valid under the hypothesis of Theorem~\ref{t:S&O}.
\end{rem}

In fact, as shown in the last section, Theorems~\ref{t:L&U} and~\ref{t:S&O} are straightforward consequences of a
series of previous results~\cite{VasGrMaz1, VasGrSt, ZavLU, GrCovCl} and the following theorem
whose proof is the main goal of this paper.

\begin{theorem}\label{t:main} Let $L$ be a simple classical group over a finite field of
characteristic $p$, and $G$ be a finite group with $\omega(G)=\omega(L)$. Suppose that
$n\geqslant45$ for $L=L_n^\varepsilon(q)$, $n\geqslant28$ for $L\in\{S_{2n}(q), O_{2n+1}(q)\}$,
$n\geqslant31$ for $L=O_{2n}^+(q)$, and $n\geqslant30$ for $L=O_{2n}^-(q)$. Then $G$ has a unique
nonabelian composition factor~$S$, and $S$ is not isomorphic to a group of Lie type over a field of
characteristic distinct from~$p$.
\end{theorem}

We strongly believe that the conclusion of Theorem~\ref{t:main} remains true for all simple classical groups except the well-known examples of isomorphic groups in different
characteristics such as $L_2(4)\simeq L_2(5)$; three groups $L_3(3)$, $U_3(3)$, and $S_4(3)$, which are the only simple groups of Lie type isospectral to some solvable groups \cite[Corollary~1]{10Zav.t}; and a few exotic cases such as $\omega(U_3(5))=\omega(2^{18}:L_3(4))$ \cite{98Maz.t}. Moreover, for many classical groups of small dimensions and specific characteristics the conclusion of the theorem has already been proved. In the article we concentrate on groups of sufficiently large dimensions in order to achieve a generic proof covering classical groups of all types in all characteristics.

A final remark. The determination of properties of a group by means of its element orders is widely applied in computational group theory, especially in development of the so-called black-box algorithms, i.\,e. algorithms that do not exploit specific features of a group representation. We mention here just one of the numerous results on this subject, which is nearest to our main assertion. Namely, W. Kantor and {\'A}. Seress in \cite{09KanSer} proved that the characteristic of a finite simple group $G$ of Lie type can be determined if three greatest element orders of $G$ are known (it is additionally assumed that the characteristic of $G$ is odd). One may observe that Theorem~\ref{t:main} says the same thing but only for classical groups of large dimensions and involving the whole spectrum of~$G$. However, we do not presuppose that $G$ is a simple group, and this is an important distinction between our approach and that of~\cite{09KanSer}. On the other hand, in contrast to~\cite{09KanSer} we do not propose here any practical implementation of our results.

\section{Preliminaries: arithmetic of Zsigmondy primes}\label{s:arithmetic}

For nonzero integers $n_1,\ldots,n_k$, let $(n_1,\ldots,n_k)$ and $[n_1,\ldots,n_k]$ denote their
greatest common divisor and least common multiple, respectively. Given a nonzero integer $n$, we put
$\varphi(n)$ for the Euler totient function of $n$, $\pi(n)$ for the set of prime divisors of $n$,
and if $G$ is a finite group then, as usual, $\pi(G)$ stands for $\pi(|G|)$. If $\pi$ is a set of
primes, then $n_\pi$ denotes the $\pi$-part of $n$, that is, the largest divisor $k$ of $n$ with
$\pi(k)\subseteq \pi$; and $n_{\pi'}$ denotes the $\pi'$-part of $n$, that is, the ratio
$|n|/n_{\pi}$. If $n$ is a nonzero integer and $r$ is an odd prime with $(r,n)=1$, then $e(r,n)$
denotes the multiplicative order of $n$ modulo $r$. Given an odd integer $n$, we put $e(2,n)=1$ if
$n\equiv 1\pmod 4$, and $e(2,n)=2$ otherwise.

Fix an integer $a$ with $|a|>1$. A prime $r$ is said to be a \emph{primitive prime divisor} of
$a^i-1$ if $e(r,a)=i$. We write $r_i(a)$ to denote some primitive prime divisor of $a^i-1$, if such
a prime exists, and $R_i(a)$ to denote the set of all such divisors. Zsigmondy \cite{Zs} proved
that primitive prime divisors exist for almost all pairs $(a,i)$.

\begin{lemma}[Zsigmondy]\label{l:zsigmondy}
Let $a$ be an integer and $|a|>1$. For every natural number $i$ the set $R_i(a)$ is nonempty,
except for the pairs $(a,i)\in\{(2,1),(2,6),(-2,2),(-2,3),(3,1),(-3,2)\}$.
\end{lemma}

For $i\neq2$ the product of all primitive divisors of $a^i-1$ taken with multiplicities is denoted
by $k_i(a)$. Put $k_2(a)=k_1(-a)$. The number $k_i(a)$ is said to be the \emph{greatest primitive
divisor} of $a^i-1$. It follows from the definition that $(k_i(a),k_j(a))=1$ if $i\neq j$. It is
easy to check that $k_1(a)=|a-1|$ if $a\not\equiv3\pmod4$, and $k_1(a)=|a-1|/2$ if
$a\equiv3\pmod4$, as well as $k_2(a)=|a+1|$ if $a\not\equiv1\pmod4$, and $k_2(a)=|a+1|/2$ if
$a\equiv 1\pmod4$. It follows from \cite{R} that for $i>2$,
\begin{equation}\label{eq:ki}
k_i(a)=\frac{|\Phi_i(a)|}{(r,\Phi_{i_{\{r\}'}}(a))},
\end{equation}
where $\Phi_i(x)$ is the $i$th cyclotomic polynomial and $r$ is the largest prime dividing $i$;
moreover, if $i_{\{r\}'}$ does not divide $r-1$ then $(r,\Phi_{i_{\{r\}'}}(a))=1$.

Note that for a divisor, the property of being primitive depends on the pair $(a,i)$ and is not
determined by the number $a^i-1$. For example, $k_6(2)=1$ but $k_3(4)=7$, and $k_2(2)=3$ but
$k_2(-2)=1$.

\begin{lemma}\label{l:kjpu0} Suppose that $a$ and $i$ are integers with $|a|>1$ and $i>0$, and $p$ is a prime.
Then $k_{ip}(a)$ divides $k_i(a^p)$. Furthermore, if $p$ divides $i$ then $k_{ip}(a)=k_i(a^p)$.
\end{lemma}

\begin{prf}
Let $r$ be odd. Then $r\in R_{ip}(a)$ means that the multiplicative order of $a$ modulo $r$ equals
$ip$. Hence the order of $a^p$ modulo $r$ equals $i$, that is $r\in R_i(a^p)$. Vice versa, if the
order of $a^p$ modulo $r$ is equal to $i$, and $p$ divides $i$, then the order of $a$ modulo $r$ is
equal to $ip$. Thus, for odd $k_{ip}(a)$ and $k_i(a^p)$ the assertion holds by the definition of a
greatest primitive divisor.

Assume that $2\in R_{ip}(a)$. Then $a$ is odd, $p=2$, and $i=1$. Therefore, $k_{ip}(a)=k_2(a)$
divides $|a+1|$, and so divides $a^2-1=k_1(a^2)=k_i(a^p)$. If $2\in R_{i}(a^p)$ and $p$ divides
$i$, then $a$ is odd, $p=2$, and $i=2$, so $k_{i}(a^p)=k_2(a^2)=(a^2+1)/2=k_4(a)=k_{ip}(a)$.
\end{prf}

\begin{lemma}\label{l:kipma} Let $a$ and $i$ be integers with $|a|>1$ and $i>0$.
If $i$ is odd then $k_i(-a)=k_{2i}(a)$, and if $i$ is a multiple of $4$ then $k_i(-a)=k_i(a)$.
\end{lemma}

\begin{prf} Let $i$ be odd. By the definition of $k_2(a)$, we may assume that $i\geqslant3$, so both $k_i(-a)$ and
$k_{2i}(a)$ are odd. The order of $a$ modulo an odd prime $r$ is equal to $2i$ if and only if the
order of $-a$ modulo $r$ is equal to $i$. Hence $R_i(-a)=R_{2i}(a)$. It follows that
$k_i(-a)=k_{2i}(a)$ because $k_{2i}(a)$ divides $a^i+1$.

The latter assertion follows from Lemma~\ref{l:kjpu0}. Indeed, put $i=2\cdot2j$ and observe that $k_i(-a)=k_{2\cdot2j}(-a)=k_{2j}(a^2)=k_{2\cdot2j}(a)=k_i(a)$.
\end{prf}

\begin{lemma}\label{l:ki}
Suppose that $a$, $i$, and $\gamma$ are integers with $|a|>1$, $i>0$, and $\gamma>1$, $r$ is an odd
prime such that $(r,a)=1$ and $r$ divides $k_i(a)-1$. The following hold:

\emph{(i)} if $i=2^\gamma$, then $e(r,a)$ divides $2^{\gamma-1};$

\emph{(ii)} if $i=3\cdot2^\gamma$, then $e(r,a)$ divides $2^{\gamma};$

\emph{(iii)} if $i=5\cdot2^{\gamma+1}$, then $e(r,a)$ divides $2^{\gamma+1};$

\emph{(iv)} if $i=7\cdot2^{\gamma}$, then $e(r,a)$ divides $3\cdot2^{\gamma};$

\emph{(v)} if $i=9\cdot2^{\gamma}$, then $e(r,a)$ divides $3\cdot2^{\gamma-1};$

\emph{(vi)} if $i=11\cdot2^{\gamma}$, then $e(r,a)$ divides $5\cdot2^{\gamma}.$

In particular,  $e(r,a)\leqslant i/2$. Moreover, $e(r,a)\leqslant i/3$ in~\emph{(ii)}
and~\emph{(v)}.
\end{lemma}

\begin{prf} Observe that, by~(\ref{eq:ki}), for $\gamma\geqslant2$ and
$i\in\{2^\gamma,3\cdot2^\gamma,5\cdot2^{\gamma+1},7\cdot2^\gamma, 9\cdot2^\gamma,11\cdot2^\gamma\}$
we have $k_i(a)=k_i(-a)=\Phi_i(a)$. It is easy to verify the following equalities:

if $i=2^\gamma$, then $\Phi_{i}(a)-1=a^{2^{\gamma-1}}$ for even $a$, and
$\Phi_{i}(a)-1=(a^{2^{\gamma-1}}-1)/2$ for odd $a$;

if $i=3\cdot2^\gamma$, then $\Phi_{i}(a)-1=a^{2^{\gamma-1}}(a^{2^{\gamma-1}}+1)$;

if $i=5\cdot2^{\gamma+1}$, then
$\Phi_{i}(a)-1=a^{2^{\gamma-1}}(a^{2^\gamma}+1)(a^{2^{\gamma-1}}-1)$;

if $i=7\cdot2^\gamma$, then
$\Phi_{i}(a)-1=a^{2^{\gamma-1}}(a^{3\cdot2^{\gamma}}-1)/(a^{2^{\gamma-1}}+1);$

if $i=9\cdot2^\gamma$, then $\Phi_{i}(a)-1=a^{3\cdot2^{\gamma-1}}(a^{3\cdot2^{\gamma-1}}-1);$

if $i=11\cdot2^\gamma$, then
$\Phi_{i}(a)-1=a^{2^{\gamma-1}}(a^{5\cdot2^{\gamma}}-1)/(a^{2^{\gamma-1}}+1)$.

These equalities yield the lemma.
\end{prf}

\begin{lemma}\label{l:estimknbyvarphi}
Let $a$ and $i$ be integers, and $\varepsilon\in\{+,-\}$. If $a\geqslant 2$, $i\geqslant 3$, and
$(a,i)\not\in\{(2,3),(2,6)\}$, then $k_i(\varepsilon a)>a^{\varphi(i)/2}$.
\end{lemma}

\begin{prf}
We prove that $k_i(a)>a^{\varphi(i)/2}$ first. Let $r$ be the greatest prime divisor of $i$ and
$i=r^\alpha{k}$ where $(r,k)=1$. It follows that $k_i(a)=\Phi_i(a)/(r,\Phi_k(a))$, and if $r-1$ is
not a multiple of $k$, then $(r,\Phi_k(a))=1$. If $(r,\Phi_k(a))=1$, then the desired inequality
holds by \cite[Lemma 7]{R}, so we assume that $r$ divides $\Phi_k(a)$ and, in particular, $k$
divides $r-1$.

As observed in \cite{R}, the inequality $$\Phi_i(a)\geqslant
\left(\frac{b^r+1}{b+1}\right)^{\varphi(k)},$$ where $b=a^{r^{\alpha-1}}$, holds true. Since
$\Phi_k(a)\leqslant (a+1)^{\varphi(k)}$, we have
\begin{equation}\label{eq:estimkn}
k_i(a)\geqslant  \frac{(b^r+1)^{\varphi(k)}}{r(b+1)^{\varphi(k)}}\geqslant
\frac{(b^r+1)^{\varphi(k)}}{(a+1)^{\varphi(k)}(b+1)^{\varphi(k)}}\geqslant
\left(\frac{b^r+1}{(b+1)^2}\right)^{\varphi(k)}.
\end{equation}

Let $r\geqslant 5$ and $(r,b)\neq (5,2)$. Then $b^{(r+1)/2}>(b+1)^2$ and, therefore,
$k_i(a)>b^{\varphi(k)(r-1)/2}=a^{\varphi(i)/2}$ by~(\ref{eq:estimkn}).

Let $r=5$ and $b=2$. Then $a=2$, $\alpha=1$. Furthermore, $k$ divides $r-1=4$. If $k\in\{1,2\}$,
then $(5,\Phi_k(2))=1$, so we may assume that $k=4$ and $i=20$. In this case the assertion follows
because $k_{20}(2)=k_{20}(-2)=41>16=2^{\varphi(20)/2}$.

Let $r=3$. Then $k\in\{1,2\}$ and $\varphi(k)=1$. If $b\geqslant4$ then $b^2>3(b+1)$. Therefore,
by~(\ref{eq:estimkn}),
$$k_i(a)\geqslant
\frac{(b^3+1)^{\varphi(k)}}{3(b+1)^{\varphi(k)}}=\frac{b^3+1}{3(b+1)}>b=b^{\varphi(k)}=a^{\varphi(i)/2}.$$
Thus, $b\in\{2,3\}$, so $a\in\{2,3\}$ and $\alpha=1$. Since $i=r^\alpha{k}\in\{3,6\}$, the case
$a=2$ is impossible by the hypothesis. If $a=3$ then $(3,\Phi_k(3))=1$.

Let, finally, $r=2$. Then $k_i(a)=(a^{\varphi(i)}+1)/(2,a-1)\geqslant a^{\varphi(i)}/2\geqslant
a^{\varphi(i)/2}$, as required. Thus, the inequality $k_i(a)>a^{\varphi(i)/2}$ is proved.

Now we apply Lemma~\ref{l:kipma}. If $i\equiv0\pmod4$ then $k_i(-a)=k_i(a)$. If $i\equiv2\pmod4$
then $\varphi(i)=\varphi(i/2)$, so $k_i(-a)=k_{i/2}(a)>a^{\varphi(i/2)/2}=a^{\varphi(i)/2}$.
Finally, if $i$ is odd then $\varphi(i)=\varphi(2i)$, hence
$k_i(-a)=k_{2i}(a)>a^{\varphi(2i)/2}=a^{\varphi(i)/2}$. The lemma is proved.

\end{prf}

Define the following function on positive integers:

\begin{equation}\label{eq:eta}
\eta(k)=\left\{\begin{array}{l}
k,\text{ if }k \text{ is odd},\\
k/2,\text{ if }k\text{ is even.}\end{array}\right.
\end{equation}

\begin{lemma}\label{l:kjubyu7}
Let $u$ be a prime power, $\varepsilon\in\{+,-\}$, $p$ be an odd prime, and $j$ be a natural
number with $\eta(j)\geqslant 11$. Then $k_j(\varepsilon u)>u^7$ and $k_{jp}(\varepsilon
u)>u^{5p}$.
\end{lemma}

\begin{prf} By Lemma~\ref{l:estimknbyvarphi}, the assertion holds for sufficiently large~$j$. Indeed,
if $\varphi(j)\geqslant 15$, then $$k_j(\varepsilon u)>u^{\varphi(j)/2}>u^7,$$
$$k_{jp}(\varepsilon u)>u^{\varphi(jp)/2}\geqslant u^{\varphi(j)(p-1)/2}\geqslant u^{15(p-1)/2}\geqslant u^{5p},$$
where the latest inequality follows because $p\geqslant 3$. Therefore, we may assume that
$\varphi(j)\leqslant 14$ and, in particular, $j\leqslant \varphi(j)^2\leqslant 196$.

Let $M$ stand for the set of $j$, satisfying the hypothesis of the lemma and the inequality
$\varphi(j)\leqslant 14$. By brute-force attack, we obtain that $M=\{42$, $36$, $30$, $28$, $26$,
$24$, $22$, $21$, $15$, $13$, $11\}$. If $j=21, 15, 13, 11$, then $k_{2j}(\varepsilon
u)=k_{j}(-\varepsilon u)$ by~(\ref{eq:ki}), so it is sufficient to prove the assertion for all
$j\in M'= \{36$, $28$, $24$, $21$, $15$, $13$, $11\}$.

We use the following inequalities on $\Phi_j(\varepsilon{u})$ for $j\in M'$:
$$u^{12}>\Phi_{36}(\varepsilon u)=u^{12}-u^6+1>u^{11},$$
$$u^{12}>\Phi_{28}(\varepsilon u)=\frac{u^{14}+1}{u^2+1}>u^{11},$$
$$u^8>\Phi_{24}(\varepsilon u)=u^8-u^4+1>u^7,$$
$$u^{13}>\Phi_{21}(\varepsilon u)=u^{12}-\varepsilon u^{11}+\varepsilon u^9-u^8+u^6-u^4+\varepsilon u^3-\varepsilon u+1>u^{11},$$
$$u^{9}>\Phi_{15}(\varepsilon u)= u^8-\varepsilon u^7+\varepsilon u^5-u^4+ \varepsilon u^3- \varepsilon u+1>u^7,$$
$$u^{13}>\Phi_{13}(\varepsilon u)=\frac{u^{13}-\varepsilon1}{u-\varepsilon1}>u^{11},$$
$$u^{11}>\Phi_{11}(\varepsilon u)=\frac{u^{11}-\varepsilon1}{u-\varepsilon1}>u^{9}.$$

By (\ref{eq:ki}), for $j=36$, $28$, $24$, $15$  the equality $k_j(\varepsilon
u)=\Phi_j(\varepsilon u)$ holds, so $k_j(\varepsilon u)>u^7$ by the above inequalities. For the other
numbers from $M'$ we have:
$$k_{21}(\varepsilon u)=\frac{\Phi_{21}(\varepsilon u)}{(7,\Phi_3(\varepsilon u))}>\frac{u^{11}}{u^2+u+1}>u^8,$$
$$k_{13}(\varepsilon u)=\frac{\Phi_{13}(\varepsilon u)}{(13, u-\varepsilon 1)}>\frac{u^{11}}{u+1}>u^{9},$$
$$k_{11}(\varepsilon u)=\frac{\Phi_{11}(\varepsilon u)}{(11, u-\varepsilon 1)}>\frac{u^{9}}{u+1}>u^7.$$
Thus, the first required inequality is proved.

If $p$ divides $j$ then $k_{jp}(\varepsilon u)=k_j(\varepsilon u^p)$ by Lemma~\ref{l:kjpu0}, so $k_{jp}(\varepsilon u)>u^{7p}$. Therefore, we assume that $p$ does not divide $j$. Then
$\Phi_{jp}(x)=\Phi_j(x^p)/\Phi_j(x)$. Suppose that $p\geqslant13$. Then $p$ is the greatest prime
divisor of~$jp$, hence, using~(\ref{eq:ki}), we obtain
$$k_{jp}(\varepsilon u)=\frac{\Phi_{j}(\varepsilon u^p)}{\Phi_{j}(\varepsilon u)
(p,\Phi_{j}(\varepsilon u))}\geqslant \frac{\Phi_{j}(\varepsilon u^p)}{\Phi^2_{j} (\varepsilon
u)}>\frac{u^{7p}}{u^{26}}\geqslant u^{5p}.$$ Suppose that $5\leqslant p\leqslant 11$. In this case
the greatest prime divisor of $jp$ is at most $13$. It follows that
$$k_{jp}(\varepsilon u)=\frac{\Phi_{j}(\varepsilon u^p)}{\Phi_{j}(\varepsilon u)(p,\Phi_{j}(\varepsilon u))}
\geqslant \frac{\Phi_{j}(\varepsilon u^p)}{13\Phi_{j}(\varepsilon
u)}>\frac{u^{9p}}{13u^{13}}>\frac{u^{9p}}{u^{17}}> u^{5p},$$ unless $j$ equals $15$ or $24$. For
$j\in\{15,24\}$, the greatest prime divisor of $jp$ is equal to $p$, and $p-1$ is not a multiple of
$j$, hence
$$k_{jp}(\varepsilon u)=\frac{\Phi_{j}(\varepsilon u^p)}{\Phi_{j}(\varepsilon u)}>\frac{u^{7p}}{u^{9}}>u^{5p}.$$
Let $p=3$. Then $j=28, 13, 11$, and the assertion follows from the inequalities:
$$k_{54}(\varepsilon u)=\frac{\Phi_{28}(\varepsilon u^3)}{\Phi_{28}(\varepsilon u)}>\frac{u^{33}}{u^{12}}>u^{15},$$
$$k_{39}(\varepsilon u)=\frac{\Phi_{13}(\varepsilon u^3)}{13\Phi_{13}(\varepsilon u)}>\frac{u^{33}}{13u^{13}}>u^{15},$$
$$k_{33}(\varepsilon u)=\frac{\Phi_{11}(\varepsilon u^3)}{\Phi_{11}(\varepsilon u)}>\frac{u^{27}}{u^{11}}>u^{15}.$$

The lemma is proved.
\end{prf}

\begin{lemma}\label{l:rpart} Let $q$ and $m$ be integers greater than $1$, and $\varepsilon\in\{+,-\}$.

\emph{(i)} If an odd prime $r$ divides $\varepsilon{q}-1$, then
$((\varepsilon{q})^m-1)_{\{r\}}=m_{\{r\}}(\varepsilon{q}-1)_{\{r\}}$.

\emph{(ii)} If an odd prime $r$ divides $(\varepsilon{q})^m-1$, then $r$ divides
$(\varepsilon{q})^{m_{\{r\}'}}-1$.

\emph{(iii)} If $\varepsilon{q}-1$ is a multiple of $4$, then
$((\varepsilon{q})^m-1)_{\{2\}}=m_{\{2\}}(\varepsilon{q}-1)_{\{2\}}$.
\end{lemma}
\begin{prf}
See, for example, \cite[Chapter IX, Lemma 8.1]{HB}.
\end{prf}

Let $[x]$ denote the integer part of a real number~$x$.

\begin{lemma}\label{l:eta}
Let $a,b$ be positive integers, $b>a$, and $A=\{i\in \mathbb{N}\mid b-a< \eta(i)\leqslant b\}$. If
$b$ is even then $|A|=[3a/2]$. If $b$ is odd then $|A|=[(3a+1)/2]$.
\end{lemma}

\begin{prf}
Since the equation $\eta(x)=j$ has two integer solutions for odd $j$ and one solution for even $j$,
the cardinality of $A$ is equal to the sum of $a$ and the quantity of odd numbers in the interval
$(b-a,b]$. If $b$ is even then this quantity is equal to $[a/2]$, and if $b$ is odd then it equals
$[(a+1)/2]$.
\end{prf}

\begin{lemma}\label{l:intervalwithprime} If $n$ is a natural number and $n\geqslant30$, then the interval $(5n/6,n)$ contains
a prime. If, in addition, $n\neq 35, 36, 37, 53$, then the interval $(8n/9,n)$ contains a prime.
\end{lemma}

\begin{prf} The first part is proved in~\cite{Nag}. It is shown in~\cite{RW} that for every
natural number $m\geqslant 119$, the interval $[m, 1.073m]$ contains a prime. Since $9/8>1.073$,
the second part holds for all sufficiently large $n$ (precisely, for $n\geqslant 129$). For
smaller $n$ the assertion can be verified directly.
\end{prf}

\section{Preliminaries: the prime graph and the spectrum of a finite group}\label{s:primegraph}

The \emph{prime graph} $GK(G)$ of a finite group $G$ is the nonoriented graph with the vertex set
$\pi(G)$ and two distinct vertices $r$ and $s$ are adjacent if and only if $rs\in\omega(G)$. The
notion of prime graph was introduced by G.K.~Gruenberg and O.~Kegel (for this reason it is also
called the Gruenberg--Kegel graph). They established that a finite group with disconnected prime
graph is either Frobenius or $2$-Frobenius group, or has a unique nonabelian composition factor
with disconnected prime graph. J.S.~Williams \cite{Wil} published this result and started the
classification of finite simple groups with disconnected prime graph. The classification was
completed by A.S.~Kondrat$'$ev~\cite{Kon}. The full list of nonabelian simple groups with
disconnected prime graph can be found, e.g., in~\cite[Tables~1a-1c]{MazS4}.

Recall that an independent set of vertices or a \emph{coclique} of a graph $\Gamma$ is any subset
of pairwise nonadjacent vertices of $\Gamma$. We write $t(\Gamma)$ to denote the independence
number of $\Gamma$, that is the greatest size of coclique in~$\Gamma$. For a group $G$, put
$t(G)=t(GK(G))$. By analogy, for each prime $r$, define the $r$-independence number $t(r,G)$ to be
the greatest size of cocliques containing the vertex $r$ in $GK(G)$. For convenience, we refer to a
coclique containing $r$ as an $\{r\}$-coclique. In~\cite{VasStr} there was proved the following
assertion which is, in some sense, a generalization of the Gruenberg--Kegel theorem (below we give
the statement of this result from~\cite[Theorem~1]{VasGor}).

\begin{lemma}\label{l:structure}
Let $G$ be a finite group with $t(G)\geqslant3$ and $t(2,G)\geqslant2$. Then the following hold:

\emph{(i)} There exists a nonabelian simple group $S$ such that
$S\leqslant\overline{G}=G/K\leqslant\Aut S$, where $K$ is the maximal normal soluble subgroup of
$G$.

\emph{(ii)} For every coclique $\rho$ of $GK(G)$ of size at least~$3$, at most one prime of $\rho$
divides the product $|K|\cdot|\overline{G}/S|$. In particular, $t(S)\geqslant t(G)-1$.

\emph{(iii)} One of the following holds:

\hspace{0.5cm}\emph{(a)} every prime $r\in\pi(G)$ nonadjacent to $2$ in $GK(G)$ does not divide
the product $|K|\cdot|\overline{G}/S|;$ in particular, $t(2,S)\geqslant t(2,G);$

\hspace{0.5cm}\emph{(b)}  there exists a prime $r\in\pi(K)$ nonadjacent to $2$ in $GK(G)$; in
which case $t(G)=3$, $t(2,G)=2$, and $S\simeq A_7$ or $L_2(q)$ for some odd $q$.
\end{lemma}

Let $L$ and $G$ be isospectral finite groups. It follows by the definition of prime graph that
$GK(L)=GK(G)$. Therefore, if $L$ satisfies the hypothesis of Lemma~\ref{l:structure}, then so
does~$G$.  In~\cite{VasVd,VasVdM}, for every finite nonabelian simple group, an adjacency criterion
of its prime graph is developed and all cocliques and $\{2\}$-cocliques of greatest size in this
graph are found, as well as $\{p\}$-cocliques of greatest size for groups of Lie type in
characteristic $p$. This information and Williams--Kondrat$'$ev's classification, imply the
following assertion which is the first step toward a proof of main results of the present paper.

\begin{prop}\label{p:UniqueS}
Let $L$ be a finite simple group of Lie type different from the groups $L_3(3)$, $U_3(3)$, and
$S_4(3)$. If $G$ is a finite group isospectral to~$L$, then the following hold true for $G$:

\emph{(i)} There exists a nonabelian simple group $S$ such that
$S\leqslant\overline{G}=G/K\leqslant\Aut S$, where $K$ is the maximal normal soluble subgroup of
$G$.

\emph{(ii)} For every coclique $\rho$ of $GK(G)$ of size at least~$3$, at most one prime of $\rho$
divides the product $|K|\cdot|\overline{G}/S|$. In particular, $t(S)\geqslant t(L)-1$.

\emph{(iii)} Every prime $r\in\pi(G)$ nonadjacent to $2$ in $GK(G)$ does not divide the product
$|K|\cdot|\overline{G}/S|$. In particular, $t(2,S)\geqslant t(2,L).$
\end{prop}

\begin{prf}
If $L$ is a finite nonabelian simple group with connected prime graph and is different from an
alternating group, then it satisfies the hypothesis of Lemma~\ref{l:structure} by~\cite{VasVd}. If
$L$ has a disconnected prime graph then the existence and uniqueness of a nonabelian composition
factor $S$ follow from the Gruenberg--Kegel theorem and~\cite{Ale}.
By~\cite[Propositions~2,3]{VasStr}, the inequality $t(2,G)\geqslant2$ and the insolubility of $G$
imply that (ii) and (iii) of Lemma~\ref{l:structure} hold. It remains to observe that the
exceptional case (b) of assertion (iii) of Lemma~\ref{l:structure} does not hold by
\cite[Theorem~2]{VasGor}. The proposition is proved.
\end{prf}

For a classical group $L$, we put $\prk(L)$ to denote its dimension if $L$ is a linear or unitary
group, and its Lie rank if $L$ is a symplectic or orthogonal group. Observe that $n=\prk(L)$ in
Theorems~1--3 in Introduction.

\begin{prop}\label{p:notexcept} Suppose that $L$ is a finite simple classical group,
$\prk(L)\geqslant27$ if $L$ is linear or unitary, and $\prk(L)\geqslant19$ if $L$ is symplectic or
orthogonal. Suppose that $G$ is a finite group isospectral to $L$, and $S$ is a unique nonabelian
composition factor of $G$. If $S$ is a group of Lie type, then $S$ is a classical group,
$\prk(S)\geqslant25$ if $S$ is linear or unitary, and $\prk(S)\geqslant16$ if $S$ is symplectic or
orthogonal.
\end{prop}

\begin{prf} Applying \cite[Tables~2,3]{VasVdM}, it is easy to obtain that $t(L)\geqslant14$ provided
the hypothesis of the proposition. By Proposition~\ref{p:UniqueS}(ii), we have $t(S)\geqslant13$.
On the other hand, it follows from~\cite[Table~4]{VasVdM} that $t(H)\leqslant t(E_8(u))=12$ for
every simple exceptional group $H$ of Lie type. Thus, $S$ is a classical group. The required
inequalities on $\prk(S)$ hold by~\cite[Tables~2,3]{VasVdM}. The proposition is proved.
\end{prf}

Given a classical group $L$ over a field of order $q$, put
\begin{equation}
\delta(L)=\begin{cases} \pi(\varepsilon{q}-1), \text{ if } L=L_n^\varepsilon(q),\\
                        \pi((2,q-1)), \text{ if $L$ is symplectic or orthogonal}. \end{cases}
\end{equation}

\begin{lemma}\label{l:indepofcrit}
Let $L$ be a simple classical group over a field of order $q$ and characteristic~$p$. Suppose
that $r$ and $s$ are distinct primes, $r,s\not\in\delta(L)$, $r\in R_i(q)$, and $s\in R_j(q)$.

\emph{(i)} If $rs\in\omega(L)$, then $r's'\in\omega(L)$ for every distinct odd primes $r'\in
R_i(q)$ and $s'\in R_j(q)$.

\emph{(ii)} If $pr\in\omega(L)$, then $pr'\in\omega(L)$ for every odd prime $r'\in R_i(q)$.
\end{lemma}

\begin{prf} It follows from~\cite{VasVd}.
\end{prf}

Thus, for two distinct primes $r,s\in\pi(L)\setminus\delta(L)$, where $r\neq p$, the answer to the
question whether they are adjacent in $GK(L)$ depends only on $e(r,q)$, if $s=p$, and $e(r,q)$,
$e(s,q)$, if $s\neq p$.

In~\cite{VasVd} several functions of natural argument are used to formulate an adjacency criterion,
one of them is~$\eta$ defined in~(\ref{eq:eta}), and two others we define here.

\begin{equation}\label{eq:nu}
\nu(k)=\left\{\begin{array}{l}
k,\text{ if }k\equiv0\pmod4,\\
k/2,\text{ if }k\equiv2\pmod4,\\
2k,\text{ if }k\text{ is odd.}\end{array}\right.
\end{equation}

For $\varepsilon\in\{+,-\}$, put

\begin{equation}\label{eq:nue}
\nu_\varepsilon(k)=\left\{\begin{array}{l}
k,\text{ if }\varepsilon=+,\\
\nu(k)\text{ if }\varepsilon=-.\end{array}\right.
\end{equation}

It is an easy observation that $\nu_\varepsilon$ is a bijection and $\nu_\varepsilon^2$ is an
identity.

For linear and unitary groups, we exploit also a reformulation of an adjacency criterion (see
\cite[Lemmas 2.1--2.3]{VasGrSt}), if it is more convenient for our goals than an initial
formulation from \cite{VasVd} which used the function~$\nu_\varepsilon$. This reformulation is
based on the equality $k_{\nu_\varepsilon(i)}(q)=k_i(\varepsilon{q})$, which follows from
Lemma~\ref{l:kipma} and the definition of~$\nu_\varepsilon$.

Now we introduce a new function in order to unify further arguments. Namely, given a simple
classical group $L$ over a field of order $q$ and a prime~$r$ coprime to $q$, we put

\begin{equation}\label{eq:varphi}
\varphi(r,L)=\left\{\begin{array}{l}
e(r,\varepsilon q),\text{ if }L=L_n^{\varepsilon}(q),\\
\eta(e(r,q)),\text{ if }L\text{ is symplectic or orthogonal.}\end{array}\right.
\end{equation}

It follows that

\begin{equation}\label{eq:varphireverse}
e(r,q)=\left\{\begin{array}{ll}
2\varphi(r,L), & \text{if either }e(r,q)\text{ is even and }L\text{ is symplectic}\\
& \text{or orthogonal,}\\
& \text{or }e(r,q)\equiv2\pmod4\text{ and }L\text{ is unitary;}\\
\varphi(r,L)/2, & \text{if }e(r,q)\equiv1\pmod2\text{ and }L\text{ is unitary;}\\
\varphi(r,L) & \text{otherwise.}\end{array}\right.
\end{equation}

Observe that $e(r,-q)=\varphi(r,L)$ in the case of $e(r,q)=\varphi(r,L)/2$.

\begin{lemma}\label{l:fermat}
Let $L$ be a simple classical group over a field of order $q$ and characteristic~$p$. If $r$ is
an odd prime from $\pi(L)\setminus\{p\}$ then $\varphi(r,L)$ divides $r-1$, and if $L$ is a
symplectic or orthogonal group then $2\varphi(r,L)$ divides $r-1$.
\end{lemma}

\begin{prf}
If $L$ is not unitary, then $\varphi(r,L)$ divides $e(r,q)$, which divides $r-1$ by Fermat's little
theorem. If $L$ is unitary, then $\varphi(r,L)$ does not divide $e(r,q)$ only when $e(r,q)$ is odd.
But then $\varphi(r,L)=e(r,-q)=2e(r,q)$ and $\varphi(r,L)$ divides $r-1$ because $r$ is odd. Let
$L$ be a symplectic or orthogonal group. If $\varphi(r,L)$ is even then $2\varphi(r,L)=e(r,q)$, and
if not, then $\varphi(r,L)$ divides $(r-1)/2$.
\end{prf}

\begin{lemma}\label{l:auxvarphi}
Let $L$ be a simple classical group over a field of order $q$ and characteristic~$p$, and let
$\prk(L)=n\geqslant 4$.

\emph{(i)} If $r\in\pi(L)\setminus\{p\}$, then $\varphi(r,L)\leqslant n$.

\emph{(ii)} If $r$ and $s$ are distinct primes from $\pi(L)\setminus\{p\}$ with
$\varphi(r,L)\leqslant n/2$ and $\varphi(s,L)\leqslant n/2$, then $r$ and $s$ are adjacent
in~$GK(L)$.

\emph{(iii)} If $r$ and $s$ are distinct primes from $\pi(L)\setminus\{p\}$ with
$n/2<\varphi(r,L)\leqslant n$ and $n/2<\varphi(s,L)\leqslant n$, then $r$ and $s$ are adjacent
in~$GK(L)$ if and only if $e(r,q)=e(s,q)$.

\emph{(iv)} If $r$ and $s$ are distinct primes from $\pi(L)\setminus\{p\}$ and $e(r,q)=e(s,q)$,
then $r$ and $s$ are adjacent in~$GK(L)$.
\end{lemma}

\begin{prf} It follows from~\cite{VasVd, VasVdM}.
\end{prf}

Let $L$ be a simple classical group over a field of order~$q$ and characteristic~$p$. For
$\sigma\subseteq\pi(L)\setminus\{p\}$, set $E(\sigma, L)=\{e(r,q)\mid r\in\sigma\}$. If
$\prk(L)=n\geqslant13$ then, by~\cite{VasVdM}, every coclique $\rho$ of greatest size in $GK(L)$
does not contain~$p$, so the set $E(\rho, L)$ is well-defined for~$\rho$. Define $J(L)$ as the
union of sets $E(\rho, L)$, and $E(L)$ as the intersection of these sets, where $\rho$ runs over
all cocliques of greatest size in~$GK(L)$. The next lemma is a particular case of the main theorem
of~\cite{VasVdM}.

\begin{lemma}\label{l:tL}
Let $L$ be a simple classical group over a field of order $q$ and characteristic $p$, and let
$\prk(L)=n\geqslant13$. Let $\rho$ be a coclique of greatest size in~$GK(L)$. If $J(L)=E(L)$ then
$E(\rho,L)=E(L)$. If $J(L)\neq E(L)$ then $E(\rho,L)=E(L)\cup\{j\}$ for some $j\in J(L)\setminus
E(L)$. In particular, $|E(L)|\leqslant t(L)\leqslant |E(L)|+1$. The sets $E(L)$, $J(L)\setminus E(L)$ and numbers
$t(L)$ are listed in Table~\emph{\ref{tab:tL}}.
\end{lemma}

\begin{prf}
See \cite[Tables 2,3]{VasVdM}.
\end{prf}

\begin{table}[!th]
\caption{Cocliques of greatest size}\label{tab:tL}
\begin{center}
{\begin{tabular}{|c|l|c|c|c|}
  \hline
  $L$ & Conditions & $t(L)$ & $E(L)$ & $J(L)\setminus E(L)$\\
  \hline
  $L_n^\varepsilon(q)$ & $n$ odd & $\frac{n+1}{2}$ & $\{i\mid \frac{n}{2}< \nu_\varepsilon\left(i\right)\le n\}$ & $\varnothing$\\
  & $n$ even & $\frac{n}{2}$ & $\{ {i}\mid \frac{n}{2}< \nu_\varepsilon\left(i\right)<n\}$ & $\{\frac{n}{2}, n\}$\\
  \hline
  $S_{2n}(q)$ or  & $n\equiv{0}\pmod 4$ & $\frac{3n+4}{4}$ &
  $\{i\mid \frac{n}{2}\leqslant\eta(i)\leqslant n\}$ &  $\varnothing$\\
  $O_{2n+1}(q)$ & $n\equiv{1}\pmod 4$ & $\frac{3n+5}{4}$ &
  $\{i\mid \frac{n}{2}<\eta(i)\leqslant n\}$ &  $\varnothing$\\
  &  $n\equiv{2}\pmod 4$ & $\frac{3n+2}{4}$ & $\{i\mid\frac{n}{2}<\eta(i)\leqslant n\}$ & $\{\frac{n}{2},n\}$\\
  & $n\equiv{3}\pmod 4$ & $\frac{3n+3}{4}$ & $\{i\mid \frac{n+1}{2}<\eta(i)\leqslant n\}$
  &  $\{ \frac{n-1}{2},n-1,$\\
 &  &  & & $n+1\}$\\
 \hline
 $O_{2n}^+(q)$ & $n\equiv{0}\pmod 4$ & $\frac{3n}{4}$ & $\{i\mid
\frac{n}{2}\leqslant\eta(i)\leqslant n,$ & $\varnothing$\\
 & & & $i\neq2n\}$ & \\
 & $n\equiv{1}\pmod 4$ & $\frac{3n+1}{4}$ & $\{i\mid
 \frac{n}{2}<\eta(i)\leqslant n,$ & $\{n-1, n+1\}$\\
 &  &  & $i\neq2n,n+1\}$ &\\
 & $n\equiv{2}\pmod 4$ & $\frac{3n-2}{4}$ & $\{i\mid
\frac{n}{2}<\eta(i)\leqslant n,$ & $\{\frac{n}{2},n\}$\\
 &  &  & $i\neq2n\}$ &\\
 & $n\equiv{3}\pmod 4$ & $\frac{3n+3}{4}$ & $\{i\mid
\frac{n-1}{2}\leqslant\eta(i)\leqslant n,$ & $\varnothing$\\
 & & & $i\neq2n,n-1\}$ &\\
 \hline
 $O_{2n}^-(q)$& $n\equiv{0}\pmod 4$ & $\frac{3n+4}{4}$ & $\{i\mid
\frac{n}{2}\leqslant\eta(i)\leqslant n\}$ & $\varnothing$\\
 & $n\equiv{1}\pmod 4$ & $\frac{3n+1}{4}$ & $\{i\mid
\frac{n}{2}<\eta(i)\leqslant n,$
 & $\{\frac{n+1}{2}, n-1\}$\\
 & & & $i\neq n,\frac{n+1}{2}\}$ & \\
 & $n\equiv{2}\pmod 4$ & $\frac{3n+2}{4}$ & $\{i\mid
\frac{n}{2}<\eta(i)\leqslant n\}$ &
 $\{\frac{n}{2}, n-2, n\}$\\
 & $n\equiv{3}\pmod 4$ & $\frac{3n+3}{4}$ & $\{i\mid
\frac{n-1}{2}\leqslant\eta(i)\leqslant n,$ & $\varnothing$\\
 & &  &$i\neq n,\frac{n-1}{2}\}$ &\\
 \hline
\end{tabular}}
\end{center}
\end{table}

Define $J(p,L)$ as the union of sets $E(\rho\setminus\{p\}, L)$, and $E(p,L)$ as the intersection
of these sets, where $\rho$ runs over all $\{p\}$-cocliques of greatest size in $GK(L)$.

\begin{lemma}\label{l:tpL}
Let $L$ be a simple classical group over a field of order $q$ and characteristic $p$, and let
$\prk(L)=n$. Suppose that $n\geqslant4$ and $(n,\varepsilon{q})\not\in\{(4,-2),(6,2),(7,2)\}$ for
$L=L_n^\varepsilon(q);$ $n\geqslant3$ and $(n,q)\neq(3,2)$ for $L\in\{S_{2n}(q), O_{2n+1}(q)\};$
$n\geqslant4$ and $(n,q)\neq(4,2)$ for $L=O_{2n}^\pm(q)$. If $\rho$ is a $\{p\}$-coclique of
greatest size in~$GK(L)$, then $E(\rho\setminus\{p\},L)=J(p,L)=E(p,L)$, in particular,
$t(p,L)=|J(p,L)|+1$. The sets $J(p,L)$ and numbers $t(p,L)$ are listed in
Table~\emph{\ref{tab:tpL}}. In particular, if $n\geqslant9$ then $t(p,L)<t(L)$. Furthermore,
$\varphi(r,L)>n/2$ for every prime $r$ nonadjacent to~$p$ in~$GK(L)$.
\end{lemma}

\begin{prf} It follows from~\cite[Proposition~6.3, Table~4]{VasVd}.
\end{prf}

\begin{table}%[!th]
\caption{Cocliques containing the characteristic}\label{tab:tpL}
\begin{center}
{\begin{tabular}{|c|l|c|c|}
  \hline
  $L$ & Conditions & $t(p,L)$ & $J(p,L)$\\
  \hline
  $L_n^\varepsilon(q)$ & & 3 & $\{\nu_\varepsilon(n-1),\nu_\varepsilon(n)\}$\\
  \hline
  $S_{2n}(q)$ or  & $n$ is even & 2 & $\{2n\}$\\
  $O_{2n+1}(q)$ & $n$ is odd & 3 & $\{n,2n\}$\\
  \hline
  $O_{2n}^+(q)$ & $n$ is even & 3 & $\{n-1,2n-2\}$\\
  & $n$ is odd & 3 & $\{n, 2n-2\}$\\
  \hline
  $O_{2n}^-(q)$& $n$ is even & 4 & $\{n-1,2n-2,2n\}$\\
  & $n$ is odd & 3 & $\{2n-2,2n\}$\\
  \hline
\end{tabular}}
\end{center}
\end{table}

A prime $r\in\pi(L)$ is called {\em large} (with respect to $L$), if $r$ lies in some coclique of
greatest size in the prime graph $GK(L)$, and {\em small} (with respect to $L$) otherwise.

\begin{lemma}\label{l:estimlarge} Let $L$ be a simple classical group over a field of order $q$ and characteristic $p$, and let
$\prk(L)=n\geqslant13$.

\emph{(i)} If $\varphi(r,L)\geqslant n/2$, then $r$ is large with respect to~$L$.

\emph{(ii)} If  $r$ is large with respect to~$L$, then $\varphi(r,L)\geqslant n/2-1$.

\emph{(iii)} If  $r$ is large with respect to~$L$, then

\begin{equation}\label{eq:estimlarge2}
\varphi(r,L)\geqslant\left\{\begin{array}{ll}
t(L), &\text{if }L\text{ is linear or unitary;}\\
(2t(L)-4)/3, &\text{if }L\text{ is symplectic or orthogonal.}\end{array}\right.
\end{equation}

\emph{(iv)} If $\rho$ is a coclique in~$GK(L)$ and $n/2<\varphi(r,L)$ for every $r\in\rho$, then
$GK(L)$ has a coclique $\sigma$ of size $t(L)$ with $\rho\subseteq\sigma$.
\end{lemma}

\begin{prf} Apply Table~\ref{tab:tL}.
\end{prf}

\begin{lemma}\label{l:rcoclique} Let $L$ be a simple classical group over a field of order $q$ and characteristic $p$, and let
$\prk(L)=n\geqslant13$. Suppose that $r\in\pi(L)$ and $\rho$ is an $\{r\}$-coclique of greatest size
in~$GK(L)$. Then every $s\in\rho'=\rho\setminus\{r\}$ is large with respect to~$L$. Further, if $r$
is small with respect to $L$, then $\varphi(s,L)>n/2$ for every $s\in\rho'$, and $E(\rho',L)$ is
uniquely determined by~$r$. If, in addition, $i=e(r,q)>2$, then $E(\rho',L)$ is uniquely determined
by~$i$.
\end{lemma}

\begin{prf} If $r$ is large with respect to~$L$, then the assertion is obvious.
Suppose that $r$ is small and $s\in\rho'$. If $r=p$ then $\varphi(s,L)>n/2$ and $E(\rho',L)=J(p,L)$
is uniquely determined by Lemma~\ref{l:tpL}. Let $r\neq p$. Lemma~\ref{l:tpL} implies that $s\neq
p$. Since $\varphi(r,L)<n/2$ by Lemma~\ref{l:tL}, the inequality $\varphi(s,L)>n/2$ follows from
Lemma~\ref{l:auxvarphi}(ii). Let $\sigma$ be another $\{r\}$-coclique of greatest size and
$\sigma'=\sigma\setminus\{r\}$. Assume that $E(\rho',L)\neq E(\sigma',L)$. Then there is
$w\in\sigma'$ with $e(w,q)\not\in E(\rho',L)$. It follows from Lemma~\ref{l:auxvarphi}(iii) that
$\{w\}\cup\rho$ is an $\{r\}$-coclique in $GK(L)$, which contradicts to the maximality of~$\rho$. If
$i=e(r,q)>2$ then $R_i(q)$ and $\delta(L)$ are disjoint. As shown above, $\rho'$ consists of primes
large with respect to $L$, and so $\rho'$ and $\delta(L)$ are disjoint as well.
Lemma~\ref{l:indepofcrit} yields that $E(\rho',L)$ depends only on~$i$.
\end{prf}

Let $L$ be a simple classical group over a field of order $q$ and characteristic $p$, and let
$\prk(L)=n\geqslant13$. Let $r$ be small with respect to $L$ and $\rho$ be an $\{r\}$-coclique of
greatest size. As shown in Lemma~\ref{l:rcoclique}, the set $E(\rho\setminus\{r\},L)$ is contained
in $J(L)$ and does not depend on a choice of $\rho$, so we denote it by~$J(r,L)$.

\begin{lemma}\label{l:trLeq2&3&4}
Let $L$ be a simple classical group over a field of order $q$ and characteristic $p$, and let
$\prk(L)=n\geqslant13$. Suppose that $r\in\pi(L)$, $2\neq r\neq p$, and $t(r,L)\leqslant4$. Then
$e(r,q)$, $t(r,L)$, and $J(r,L)$ are listed in Table~\emph{\ref{tab:trL}}.
\end{lemma}

\begin{prf}
The application of an adjacency criterion from \cite{VasVd,VasVdM} reduces the proof to easy
arithmetical calculations.
\end{prf}

\begin{table}%[!ht]
\caption{Cocliques of size at most 4}\label{tab:trL}
\thisfloatpagestyle{empty}
\begin{center}{\small
{\begin{tabular}{|c|c|l|c|c|}
  \hline
  $L$ & $e(r,q)$ & Conditions & $t(r,L)$ & $J(r,L)$\\
  \hline
  $L_n^\varepsilon(q)$ & $\nu_\varepsilon(1)$ & $|\varepsilon{q}-1|_r=n_r$ & 3 & $\{\nu_\varepsilon(n-1), \nu_\varepsilon(n)\}$\\
  & & $|\varepsilon{q}-1|_r>n_r$  & 2 & $\{\nu_\varepsilon(n)\}$\\
  & & $|\varepsilon{q}-1|_r<n_r$  & 2 & $\{\nu_\varepsilon(n-1)\}$\\
  & $\nu_\varepsilon(2)$ & $n\equiv0\pmod2$ & 2 & $\{\nu_\varepsilon(n-1)\}$\\
  &  & $n\equiv1\pmod2$ & 2 & $\{\nu_\varepsilon(n)\}$\\
  & $\nu_\varepsilon(3)$ & $n\equiv0\pmod3$ & 3 & $\{\nu_\varepsilon(n-2),\nu_\varepsilon(n-1)\}$\\
  &  & $n\equiv1\pmod3$ & 3 & $\{\nu_\varepsilon(n-2),\nu_\varepsilon(n)\}$\\
  &  & $n\equiv2\pmod3$ & 3 & $\{\nu_\varepsilon(n-1),\nu_\varepsilon(n)\}$\\
  & $\nu_\varepsilon(4)$ & $n\equiv0\pmod4$ & 4 & $\{\nu_\varepsilon(n-3),\nu_\varepsilon(n-2),\nu_\varepsilon(n-1)\}$\\
  & &  $n\equiv1\pmod4$ & 4 & $\{\nu_\varepsilon(n-3),\nu_\varepsilon(n-2),\nu_\varepsilon(n)\}$\\
  & &  $n\equiv2\pmod4$ & 4 & $\{\nu_\varepsilon(n-3),\nu_\varepsilon(n-1),\nu_\varepsilon(n)\}$\\
  & &  $n\equiv3\pmod4$ & 4 & $\{\nu_\varepsilon(n-2),\nu_\varepsilon(n-1),\nu_\varepsilon(n)\}$\\
  \hline
  $S_{2n}(q)$ or & $1$ & & 2 & $\{2n\}$\\
  $O_{2n+1}(q)$ & $2$ & $n\equiv0\pmod2$ & 2 & $\{2n\}$\\
  &  & $n\equiv1\pmod2$ & 2 & $\{n\}$\\
  & $4$ & $n\equiv0\pmod4$ & 4 & $\{n-1,2n-2,2n\}$\\
  & &  $n\equiv1\pmod4$ & 4 & $\{n,2n-2,2n\}$\\
  & &  $n\equiv2\pmod4$ & 3 & $\{n-1,2n-2\}$\\
  & &  $n\equiv3\pmod4$ & 3 & $\{n,2n\}$\\
  & $3$ & $n\equiv4\pmod6$ & 4 & $\{2n-4,2n-2,2n\}$\\
  & $6$ & $n\equiv4\pmod6$ & 4 & $\{2n-4,n-1,2n\}$\\
  \hline
  $O_{2n}^+(q)$ & $1$ & & 2 & $\{2n-2\}$\\
  & $2$ & $n\equiv0\pmod2$ & 2 & $\{n-1\}$\\
  &  & $n\equiv1\pmod2$ & 2 & $\{n\}$\\
  & $4$ & $n\equiv0\pmod4$ & 3 & $\{n-1,2n-2\}$\\
  & &  $n\equiv1\pmod4$ & 4 & $\{n-2,2n-2,n\}$\\
  & &  $n\equiv2\pmod4$ & 3 & $\{n-1,2n-2\}$\\
  & &  $n\equiv3\pmod4$ & 3 & $\{n-2,n\}$\\
  & $3$ & $n\equiv4\pmod6$ & 4 & $\{2n-6,2n-4,2n-2\}$\\
  & $6$ & $n\equiv4\pmod6$ & 4 & $\{2n-4,n-3,n-1\}$\\
  & & $n\equiv5\pmod6$ & 4 & $\{2n-2,n-2,n\}$\\
  \hline
  $O_{2n}^-(q)$ & $1$ & & 2 & $\{2n\}$\\
  & $2$ & $n\equiv0\pmod2$ & 2 & $\{2n\}$\\
  &  & $n\equiv1\pmod2$ & 2 & $\{2n-2\}$\\
  & $4$ & $n\equiv0\pmod4$ & 4 & $\{n-1,2n-2,2n\}$\\
  & &  $n\equiv1\pmod4$ & 4 & $\{2n-4,2n-2,2n\}$\\
  & &  $n\equiv2\pmod4$ & 4 & $\{n-1,2n-4,2n-2\}$\\
  & &  $n\equiv3\pmod4$ & 3 & $\{2n-4,2n\}$\\
  & $3$ & $n\equiv5\pmod6$ & 4 & $\{2n-4,2n-2,2n\}$\\
\hline
\end{tabular}}}
\end{center}
\end{table}

\begin{lemma}\label{l:graph}
Suppose that $L$ is a finite simple classical group over a field of characteristic~$p$. Then for
every $r\in\pi(L)$ there is $s\in\pi(L)$ such that $p\neq s\neq r$ and $rs\not\in\omega(L)$.
\end{lemma}

\begin{prf} It follows from~\cite{VasVd} (e.g., see \cite[Lemma~4]{VasGr} and \cite[Lemma~12]{GrCovCl}).
\end{prf}

We note without proof a more general result following from~\cite{VasVd}: for every finite
nonabelian simple group $L$, which is not alternating, and for every $r\in\pi(L)$ there exists
$s\in\pi(L)$ not equal to $r$ and satisfying $rs\not\in\omega(L)$.

\begin{lemma}\label{l:smallerq}
Let $r$ be a prime divisor of the order of a simple classical group $L$ over a field of order
$q$, and let $e(r,q)$ divide $l\cdot2^k$, where $l\in\{1,3,5\}$, $k$ is a nonnegative integer. If
$L=L_n^\varepsilon(q)$, then either $\varphi(r,L)\leqslant 2l$ or $\varphi(r,L)=e(r,q)$. If $L$ is
a symplectic or orthogonal group, then either $\varphi(r,L)\leqslant l$ or $\varphi(r,L)=e(r,q)/2$.
\end{lemma}

\begin{prf} It follows from~(\ref{eq:varphi}) by direct verification.
\end{prf}

\begin{lemma}\label{l:hallkiq} Let $L$ be a simple classical group over a field of order $q$ and characteristic $p$, and let
$\prk(L)=n\geqslant4$. If $r\in\pi(L)\setminus\{p\}$, $i=e(r,q)$, and $n/2<\varphi(r,L)\leqslant n$,
then $L$ includes a cyclic Hall subgroup of order~$k_i(q)$.
\end{lemma}

\begin{prf}
It follows from formulae for orders of simple classical groups and information on cyclic structure
of their maximal tori (see, e.g., \cite{ButGr}).
\end{prf}

Our main source on spectra of classical group is a series of papers~\cite{ButGr,ButLU,ButSO}
containing an explicit arithmetical criterion for a natural number to lie in the spectrum of a
classical group. In particular, the following lemma is a direct corollary of these results and
properties of~$e(r,q)$.

\begin{lemma}\label{l:criterion}
Let $L$ be a simple classical group over a field of order $q$ and characteristic $p$, and let
$\prk(L)=n$. Let $k$ and $l$ be integers, $k\geqslant0$, $l>0$, and $\delta=\delta(L)$. For
$j=1,\ldots,l$, suppose that pairwise distinct primes $r_j$ lie in
$\pi(L)\setminus(\delta\cup\{p\})$ and put $i_j=e(r_j,q)$. The product $p^k r_1r_2\cdots r_l$ lies
in $\omega(L)$ if and only if the $\delta'$-part of $p^k a$ lies in $\omega(L)$, where
$$a=\begin{cases}
[(\varepsilon{q})^{\nu_\varepsilon(i_1)}-1,(\varepsilon{q})^{\nu_\varepsilon(i_2)}-1,\ldots,(\varepsilon{q})^{\nu_\varepsilon(i_l)}-1],
\text{ if }L=L_n^\varepsilon(q), \\
[q^{\eta(i_1)}+(-1)^{i_1},q^{\eta(i_2)}+(-1)^{i_2},\ldots,q^{\eta(i_l)}+(-1)^{i_l}]\text{
otherwise}.
\end{cases}$$
In particular, if $i_1,i_2,\ldots,i_l$ are greater than $2$ and pairwise distinct, then $p^k
r_1r_2\cdots r_l\in\omega(L)$ if and only if $p^k k_{i_1}(q)k_{i_2}(q)\cdots
k_{i_l}(q)\in\omega(L)$.
\end{lemma}

\begin{prf} See~\cite{ButGr,ButLU,ButSO}.
\end{prf}

\begin{rem} The last assertion of the lemma in the particular case $k=0$, $l=1$ means the following: If
$r$ with $e(r,q)=i>2$ divides $|L|$, then $k_i(q)\in\omega(L)$. If in addition $n\geqslant4$, then
obviously $k_i(q)\in\omega(L)$ for $i=1,2$. In particular, it gives another proof of
Lemma~\ref{l:auxvarphi}(iv).
\end{rem}

\begin{lemma}\label{l:bigk} Let $L$ be a simple classical group over a field of order~$q$ and characteristic~$p$, and let $\prk(L)=n$.

\emph{(i)} If $L=L_n^\varepsilon(q)$ and $n\geqslant 23$, then $\omega(L)$ contains a number $k$
with $k\geqslant q^{4t(L)/3}$ and all prime divisors of $k$ are large with respect to~$L$.

\emph{(ii)} If $L\in\{S_{2n}(q), O_{2n+1}(q)\}$ and $n\geqslant29$, or $L=O_{2n}^\varepsilon(q)$
and $n\geqslant30$, then $\omega(L)$ contains a number $k$ with $k\geqslant q^{10t(L)/9}$ and all
prime divisors of $k$ are large with respect to~$L$.

\emph{(iii)} The numbers from $\omega(L)$ do not exceed $q^{2t(L)}$.

\emph{(iv)} If $p^\gamma>2n-1$, then the exponent of a Sylow $p$-subgroup of $L$ does not
exceed~$p^\gamma$.
\end{lemma}

\begin{prf}
(i) Suppose that $L=L^\varepsilon_n(q)$, where $n\geqslant 29$. Then $t(L)=[(n+1)/2]\geqslant 15$
and $n+1\geqslant 2t(L)$. By Lemma~\ref{l:intervalwithprime}, there is a prime $j$ such that
$5(n+1)/6<j<n+1$. The inequalities $j\leqslant n$ and $j>n/2$ imply that $k_j(\varepsilon q)$ lies
in $\omega(L)$ and all its prime divisors are large with respect to $L$. Furthermore,
applying~(\ref{eq:ki}) it is easy to get the inequality $k_j(\varepsilon q)>q^{j-3}$ (see, for
example, \cite[Lemma~3.1]{VasGrMaz1}). It follows
$$j-3>\frac{5(n+1)}{6}-3\geqslant \frac{5t(L)}{3}-3\geqslant \frac{4t(L)}{3}.$$ If $23\leqslant
n\leqslant 28$, then we prove the assertion by putting $j=23$.

(ii) Let $L$ be a symplectic or orthogonal group. To prove the lemma it is sufficient to find a
prime $j$ such that $j-2\geqslant 10t(L)/9$, and either $j<n$ or $L\in\{S_{2n}(q), O_{2n+1}(q)\}$
and $j\leqslant n$. Indeed, if these conditions hold, then both numbers $k_j(q)$ and $k_{j}(-q)$
lie in $\omega(G)$, all their prime divisors are large, and at least one of $k_j(q)$ and
$k_{j}(-q)$ is greater than $q^{j-2}$ (applying~(\ref{eq:ki}) again). The required assertion will
be also proved, if we find $j$, which is a power of $2$, satisfying $n/2<j<n$ and
$k_{2j}(q)=(q^j+1)/(2,q-1)\geqslant q^{10t(L)/9}$.

Suppose that $n\geqslant 54$. Then $t(L)\geqslant (3n-2)/4\geqslant 40$. We find desired $j$
applying Lemma~\ref{l:intervalwithprime}.

Let $n$ be even. Then $t(L)\leqslant (3n+4)/4$ and, therefore, $n+1\geqslant (4t(L)-1)/3$. There
exists a prime $j$ with $8(n+1)/9<j<n+1$ and, in particular, $j>8(4t(L)-1)/27$. Since $n$ is even,
we have $j<n$.

Let $n$ be odd and $L\in\{S_{2n}(q), O_{2n+1}(q)\}$. Then $t(L)\leqslant (3n+5)/4$, and so
$n+2\geqslant (4t(L)+1)/3$. There is a prime $j$ with $8(n+2)/9<j<n+2$ and, in particular,
$j>8(4t(L)+1)/27$. Since $n$ is odd, the inequality $j\leqslant n$ holds.

Let, finally, $n$ be odd and $L=O_{2n}^\varepsilon(q)$. Then $t(L)\leqslant (3n+3)/4$, and so
$n\geqslant (4t(L)-3)/3$. There is a prime $j$ with $8n/9<j<n$ and, in particular,
$j>8(4t(L)-3)/27$.

In all cases $j>8(4t(L)-3)/27$, hence
$$j-2> \frac{8(4t(L)-3)}{27}-2=\frac{32t(L)-78}{27}\geqslant \frac{10t(L)}{9}.$$

For $n\leqslant53$ we point out $j$ explicitly.

If $48\leqslant n\leqslant 53$ then $t(L)\leqslant 41$, and if $t(L)=41$ then $n=53$ and
$L\in\{S_{2n}(q), O_{2n+1}(q)\}$. Put $j=47$, if $t(L)\leqslant 40$, and $j=53$, if $t(L)=41$. If
$44\leqslant n\leqslant 47$, then $t(L)\leqslant 36$ and $j=43$ can be taken. If $n=42, 43$, then
$t(L)\leqslant 33$ and $j=41$. If $38\leqslant n\leqslant 41$, then $t(L)\leqslant 32$ and $j$
equals 41 or 37 according to a type of group. If $n=32$, then $t(L)\leqslant 25$ and $j=31$. If
either $29\leqslant n\leqslant 31$ and $L\in\{S_{2n}(q), O_{2n+1}(q)\}$, or $n=30,31$ and
$L=O_{2n}^\varepsilon(q)$, then $t(L)\leqslant 24$ and $j=29$.

It remains to treat the case $33\leqslant n\leqslant 37$. It follows that $t(L)\leqslant 29$. If
$t(L)=29$, then $n=37$ and $L\in\{S_{2n}(q), O_{2n+1}(q)\}$, so we put $j=37$. Let $t(L)\leqslant
28$. We show that $j=32$ is suitable in this case. If $q$ is even then
$k_{2j}(q)=q^{32}+1>q^{280/9}$. If $q$ is odd then $q^{8/9}>2$, so
$k_{2j}(q)>q^{32}/2>q^{32}/q^{8/9}=q^{280/9}$.

(iii) It follows from \cite[Lemma~1.3]{VasGrMaz2} that numbers from $\omega(L)$ do not exceed
$q^{m+1}/(q-1)$, where $m$ is the Lie rank of $L$. Now the required assertion can be easily
obtained by using the formulae for $t(L)$ from Table~\ref{tab:tL}.

(iv) By \cite[Proposition~0.5]{Test}, the exponent of a Sylow $p$-subgroup of $L$ is equal to the
minimal power of $p$ greater than the maximal height $h(L)$ of a root in the root system of $L$.
Since $h(L)=n-1$ for linear and unitary groups, $h(L)=2n-3$ for orthogonal groups of even
dimension, and $h(L)=2n-1$ for symplectic groups and orthogonal groups of odd dimension, the
required assertion follows.
\end{prf}

\section{Preliminaries: actions and automorphisms}\label{s:subgroups}

In this section we collect some facts concerning spectra of covers  and automorphic extensions of
classical groups. Our main tools are well-known results on Frobenius actions, Hall--Higman type
theorems, as well as a description of parabolic subgroups and centralizers of field automorphisms
of classical groups.

We start with two known results on covers of finite groups.

\begin{lemma}[{\cite[Lemma~10]{ZavMazAS}}]\label{l:semidirect}
Suppose that $K$ is a normal elementary abelian $p$-subgroup of a finite group $G$, $H\simeq G/K$, and
$G_1=K\rtimes H$ is the natural semidirect product under the action of $H$ on $K$ via
conjugation. Then $\omega(G_1)\subseteq\omega(G)$.
\end{lemma}

\begin{lemma}[{\cite[Lemma~1.5]{VasGrSt}}]\label{l:action2} Let $G$ be a finite group, $K$ be a normal subgroup of $G$, and
$r\in\pi(K)$. Suppose that the factor group $G/K$ has a section isomorphic to a non-cyclic abelian $p$-group for some odd prime $p$ distinct from~$r$. Then
$rp\in\omega(G)$.
\end{lemma}

Now we put a result on a faithful Frobenius action and its corollary for a cover of Frobenius group.

\begin{lemma}\label{l:action1} If a Frobenius group $FC$ with kernel $F$ and cyclic complement $C=\langle c\rangle$
of order $n$ acts faithfully on a vector space $V$ of positive characteristic $p$ coprime to the order of $F$, then
the minimal polynomial of $c$ on $V$ is equal to $x^n-1$. In particular, the natural semidirect product $V\rtimes C$
contains an element of order $pn$ and $\dim C_V(c)>0$.
\end{lemma}

\begin{prf} See, e.g., \cite[Lemma~2]{ZavMazAS}.
\end{prf}

\begin{lemma}[{\cite[Lemma 1]{Maz97}}] \label{l:action} Let $G$ be a finite group,
$K$ be a normal subgroup of~$G$, and let $G/K$ be a Frobenius group with kernel $F$ and cyclic complement~$C.$
If $(|F|, |K|)=1$ and $F$ does not lie in $KC_G(K)/K,$ then
$r|C|\in\omega(G)$ for some prime divisor $r$ of~$|K|.$
\end{lemma}

Next step is to discuss several results based on theorems of Hall--Higman type. Our main source
here is Di Martino and Zalesskii's theorem \cite{DMZal} on minimal polynomials of elements of
prime-power order lying in proper parabolic subgroups of classical groups. We begin with a direct
corollary of this theorem.

\begin{lemma}\label{l:ZalCor} Let $L$ be a simple classical group over a field of order~$q$ and characteristic~$p$,
$r\in\pi(L)$, $r^s\in\omega(P)$, where $P$ is a proper parabolic subgroup of $L$, and
$(r,6p(q+1))=1$. If $L$ acts faithfully on a vector space $V$ over the field of characteristic $t$ distinct from~$p$,
then $tr^{s}\in\omega(V\rtimes L)$.
\end{lemma}

\begin{prf} Let an element $g\in P$ have an order $r^s$. Since $(r,6p(q+1))=1$, it follows from the main result of
\cite{DMZal} that the minimal polynomial of $g$ on $V$ has degree
$r^s$, so the lemma holds.
\end{prf}

The proof of Di Martino and Zalesskii's theorem is based on the application of the Hall--Higman
theorem and its cross-characteristic analogues to the restriction of faithful representation of a
classical group on its parabolic subgroups. Next two lemmas contain very similar arguments, and we
include them because their formulations are convenient for our further purposes.

\begin{lemma}\label{l:hallhigman} Let $s$ and $p$ be distinct primes, a group $H$ be a semidirect product of a normal
$p$-subgroup $T$ and a cyclic subgroup $C=\langle g\rangle$ of order~$s$, and let $[T,g]\neq1$.
Suppose that $H$ acts faithfully on a vector space $V$ of positive characteristic
$t$ not equal to~$p$. If the minimal polynomial of $g$ on $V$ does not equal $x^s-1$, then

\emph{(i)} $C_T(g)\neq1;$

\emph{(ii)} $T$ is nonabelian\emph{;}

\emph{(iii)} $p=2$ and $s=2^{2^\delta}+1$ is a Fermat prime.
\end{lemma}

\begin{prf}
If $C_T(g)=1$, then $TC$ is a Frobenius group and the minimal polynomial of $g$ on $V$ equals
$x^s-1$ by Lemma~\ref{l:action1}. If $T$ is abelian, then $T=[T,g]\times C_T(g)$ and $[T,g]C$ is a
Frobenius group acting on $V$ faithfully. Therefore, (ii) also holds. The last assertion follows
from the Hall--Higman theorem \cite[Theorem~2.1.1]{HH} in the case $t=s$, and can be easily derived
from \cite[Satz~17.13]{Hup} (see, e.g., \cite{HartM} and \cite{KhM}) for $t\neq s$.
\end{prf}

\begin{lemma}\label{l:good} Let $L$ be a simple classical group over a field of order~$q$ and characteristic $p$, and $\prk(L)>4$.
Suppose that a prime $s$ divides the order of a proper parabolic subgroup of $L$, and $(s,p(q^2-1))=1$. Then $L$ includes a subgroup
$H$ such that $H$ is a semidirect product of a normal $p$-subgroup $T$ and a cyclic subgroup
$C=\langle g\rangle$ of order $s$ with $[T,g]\neq1$, and at least one of three assertions from the conclusion of Lemma~\emph{\ref{l:hallhigman}} does not hold for~$H$.
\end{lemma}

\begin{prf}
First of all, observe that $s>3$ because $(s,p(q^2-1))=1$. Let $g$ be an element of order~$s$
in~$L$. There exists a proper parabolic subgroup $P$ of $L$ admitting the Levi decomposition $A:B$,
where $A$ is the unipotent radical, $B$ is the Levi factor, and $g\in B$. By~\cite[13.2]{GLy}, we
have $g\not\in C_L(A)$, so $[A,g]\neq1$. Suppose that the conclusion of the lemma does not hold.
Then $q=2^\beta$, $s=2^{2^\delta}+1$ is a Fermat prime, $A$ is nonabelian, and $C_A(g)\neq1$. First
we assume that $s$ is coprime to the order of a stabilizer of one-dimensional totally isotropic
(totally singular) subspace in the natural representation of~$L$. If it occurs, the one of the
following holds: $L=U_n(q)$, $e(s,-q)$ is even and $e(s,-q)\geqslant n-1$;
$L=S_{2n}(q)=O_{2n+1}(q)$ or $L=O_{2n}^+(q)$ and $e(s,q)=n$ is odd; $L=O_{2n}^-(q)$ and
$e(s,q)=n-1$ is odd. However, in all these cases, it follows from~\cite[Table~4]{VasVd} that $s$
and $p$ are not adjacent in $GK(L)$, and so $C_A(g)=1$, a contradiction. Thus, we may assume that
$P$ is a stabilizer of one-dimensional totally isotropic subspace. Since $p=2$, the unipotent
radical $A$ of $P$ is abelian, unless $L=U_n(q)$ (see, for example, \cite[Lemma~3.1]{DMZal}).
Therefore, $L=U_n(q)$. Now $s=2^{2^\delta}+1$ is a Fermat prime greater than 3, so
$e(s,2)=2^{\delta+1}\geqslant4$. Putting $e(s,q)=l$, we have $2^{\delta+1}$ divides $\beta{l}$ and
does not divide $\beta i$ for $i<l$. If $l_{\{2\}}\leqslant2$, then $2^{\delta+1}$ divides
$2\beta$, which is impossible because $l>2$. Hence $e(s,-q)=e(s,q)=l\equiv0\pmod4$. Furthermore, as
proved, $e(s,-q)<n-1$. By \cite[Lemma~5]{GrCovCl}, $L$ has a Frobenius subgroup $TC$ such that its
kernel $T$ is a $p$-subgroup, and a complement $C$ has the order~$s$. This completes the proof.
\end{prf}

The subgroup $H$ from the conclusion of Lemma~\ref{l:good} is said to be \emph{good in} $L$
\emph{with respect to a prime}~$s$.

The next lemma gives an easy criterion whether a prime divisor of the order of a classical group
divides the order of some its proper parabolic subgroup.

\begin{lemma}\label{l:adjanisotrop}
For a simple classical group  $L$ over a field of order $q$ and characteristic $p$ with
$\prk(L)=n\geqslant4$, put
$$j=\left\{\begin{array}{l}
n, \text{ if }L\simeq L_{n}(q);\\
2n-2, \text{ if either }L\simeq O^+_{2n}(q)\text{ or }L\simeq U_{n}(q)\text{ and }n\text{ is even,}\\
2n, \text{ otherwise.}\end{array}\right.$$ Then $(k_j(q),|P|)=1$ for every proper parabolic
subgroup $P$ of~$L$. If $i\neq j$ and a primitive prime divisor $r_i(q)$ lies in $\pi(L)$, then
there is a proper parabolic subgroup $P$ of $L$ such that $k_i(q)$ lies in~$\omega(P)$. In
particular, if two distinct primes $r,s\in\pi(L)$ do not divide the order of any proper parabolic
subgroup of $L$, then $r$ and $s$ are adjacent in $GK(L)$.
\end{lemma}

\begin{prf}
The order and structure of parabolic subgroups of finite classical groups are well-known (see, for
example, \cite[Propositions 4.1.17--4.1.20]{KL}). So it is easy to verify that $(k_j(q),|P|)=1$ for every
proper parabolic subgroup $P$ of~$L$. Let $r\in\pi(L)$ and $e(r,q)=i\neq j$. Since $n\geqslant4$,
it is clear that if $i\leqslant2$ then there is a proper parabolic subgroup $P$ with
$k_i(q)\in\omega(P)$. So we may assume that $i>2$. Applying \cite[Propositions 4.1.17--4.1.20]{KL}
again, we obtain a proper parabolic subgroup $P$ with the Levi factor having a section isomorphic
to a simple classical group $M$ over a field of order $q$ (with one exception when $M$ is over
a field of order $q^2$) such that $r$ divides $M$. If $M$ is over a field of order $q$, then
$k_i(q)\in\omega(M)\subseteq\omega(P)$ by Lemma~\ref{l:criterion}. In the exceptional case
$L=U_n(q)$, $n$ even, $M\simeq L_{n/2}(q^2)$. Again, using Lemma~\ref{l:criterion}, we have
$k_{n/2}(q^2)\in\omega(M)$, and $k_n(q)\in\omega(M)\subseteq\omega(P)$ by Lemma~\ref{l:kjpu0}.
Finally, if two distinct primes $r,s\in\pi(L)$ do not divide the order of any proper parabolic
subgroup of $L$, then $r,s\in R_j(q)$, so they are adjacent to each other by
Lemma~\ref{l:auxvarphi}(iv).
\end{prf}

In the end of the section we handle the spectra of extensions of classical groups by field
automorphisms.

\begin{lemma}\label{l:autoZav}
Let a symbol $X$ be chosen from the set
$\{SL^\varepsilon_n,Sp_{2n},\Omega_{2n+1},\Omega^\varepsilon_{2n}\}$. Suppose that $q$ is a prime power,
and $\tau$ is a field automorphism of odd order $t$ of the group $X(q)$. Then
\begin{equation}\label{eq:autoZav}
\omega(X(q)\rtimes\langle\tau\rangle)=\bigcup_{k\mid t}k\omega(X(q^{1/k})).
\end{equation}
\end{lemma}

\begin{prf}
If $X=SL^\varepsilon_n$, then the assertion is proved in~(i) of Corollary~14 in \cite{ZavU3}. The
proof is based on the general result on connected linear algebraic groups
\cite[Proposition~13]{ZavU3}. So it can be extended to the other cases exactly by the same way as
in the proof of Corollary~14 from \cite{ZavU3}.
\end{prf}

Given a finite group $G$ and a prime $r$, let $\exp_r(G)$ stand for the $r$-exponent of $G$, i.e.,
the exponent of its Sylow $r$-subgroup.

\begin{lemma}\label{l:AutoOrderPreserve}
Let $L$ be a simple classical group over a field of order~$q$, and $r\in\pi(L)$. Suppose
that $(r,q|\operatorname{Inndiag}L/L|)=1$, and if $L=O^\varepsilon_{8}(q)$ then $r\neq3$.  If
$L\leqslant G\leqslant\Aut{L}$, then $\exp_r(L)=\exp_r(G)$.
\end{lemma}

\begin{prf}
It is nothing to prove if $|G/L|_r=1$. Let $|G/L|_r=r^\kappa>1$. Since
$(r,q|\operatorname{Inndiag}L/L|)=1$ and $r\neq3$ for $L=O^\varepsilon_{8}(q)$, the group $G$
includes a subgroup $H$ isomorphic to the extension of $L$ by the field automorphism $\tau$ of
order~$r^\kappa$, and $\exp_r(G)=\exp_r(H)$. Choose the symbol $X$ from the statement of
Lemma~\ref{l:autoZav} so that $L$ is the nonabelian composition factor of~$X(q)$. Then
$\exp_r(X(q))=\exp_r(L)$. Furthermore, $\tau$ can be lifted to the field automorphism of $X(q)$,
which we denote by the same letter~$\tau$. Thus, $\exp_r(X(q))=\exp_r(L)\leqslant
\exp_r(H)\leqslant \exp_r(X(q)\rtimes\langle\tau\rangle)$, and it is sufficient to prove
that the $r$-exponents of $X(q)$ and $X(q)\rtimes\langle\tau\rangle$ are equal.

Put $q_0=q^{1/r}$. Since $i=e(r,q)$ divides $r-1$, we have $(r,i)=1$. It follows from
Lemma~\ref{l:rpart}(ii) that $i=e(r,q_0)$ and $(q^i-1)_{\{r\}}=r((q_0)^i-1)_{\{r\}}$. Applying
Lemma~\ref{l:rpart}(i), it is easy to see that $\exp_{r}(X(q))=(q^{ir^l}-1)_{\{r\}}$ for some
nonnegative integer~$l$. Therefore, $\exp_r(X(q))=(q^{ir^l}-1)_{\{r\}}=r^l(q^i-1)_{\{r\}}=r^{l+1}
((q_0)^i-1)_{\{r\}}=r((q_0)^{ir^l}-1)_{\{r\}}=r\exp_r(X(q_0))$. The equality
$\exp_r(X(q))=r\exp_r(X(q_0))$ yields the validity of the following chain of equalities:
\begin{equation}\label{eq:autoexpr}
\exp_r(X(q))=r\exp_r(X(q^{1/r}))=\ldots=r^\kappa\exp_r(X(q^{1/r^\kappa})).
\end{equation}
By Lemma~\ref{l:autoZav}, we have
\begin{equation}\label{eq:autoZav1}
\omega(X(q)\rtimes\langle\tau\rangle)=\bigcup_{0\leqslant l\leqslant
\kappa}r^l\omega(X(q^{1/{r^l}})).
\end{equation}
The lemma follows from (\ref{eq:autoexpr}) and~(\ref{eq:autoZav1}).
\end{prf}

\section{Proof: restrictions on $K$ and $\overline{G}/S$}\label{s:K&GS}

The following four sections contain the proof of Theorem~\ref{t:main}. Throughout $L$ is a simple
classical group over a field of order $q$ and characteristic $p$, and $\prk(L)=n$. We will prove
the theorem by contradiction. So we assume that there exists a finite group $G$ isospectral to $L$
with a unique nonabelian composition factor $S$ isomorphic to a simple group of Lie type over a
field of order $u$ and characteristic $v$ distinct from $p$. Further, since the assumptions on
$\prk(L)$ in the hypothesis of Proposition~\ref{p:notexcept} are weaker than ones from the
hypothesis of Theorem~\ref{t:main}, we obtain that $S$ is a classical group and put $\prk(S)=m$.
Here and below $K$ is the soluble radical of $G$, $\overline{G}=G/K$, $S$ is treated as a subgroup
of $\overline{G}$, so $S\trianglelefteq \overline{G}\leqslant\Aut S$ and
$\overline{G}/S\leqslant\Out S$.

The purpose of this section is to prove the following three propositions under the assumption that
the dimension of the natural representation of $L$, denoted by $\dim L$, is at least $40$. Observe
that the hypothesis of Theorem~\ref{t:main} yields $\dim L\geqslant40$, and the later inequality
implies the validity of the assumptions on $\prk(L)$ from Proposition~\ref{p:notexcept}.

\begin{prop}\label{p:structureK}
Suppose that $\dim L\geqslant40$. Then the soluble radical $K$ of $G$ is nilpotent. If
$r\in\pi(K)\setminus\{v\}$, then $t(r,L)=2$, and $(s,|K|\cdot|\overline{G}/S|\cdot|P|)=1$ for
every $s\in\pi(L)$ nonadjacent to $r$ in $GK(L)$ and every proper parabolic subgroup $P$ of~$S$.
\end{prop}

\begin{prop}\label{p:structureGS}
Suppose that $\dim L\geqslant40$. If a prime $r$ not equal to $p$ divides the order of
$\overline{G}/S$, then either $\varphi(r,L)\leqslant n/3$, or $L=L^\varepsilon_n(q)$,
$n\in[2^{\gamma+3},9\cdot2^\gamma)$, and $e(r,\varepsilon{q})=\varphi(r,L)=3\cdot2^\gamma$ for some
integer $\gamma\geqslant3$. In particular, $r$ is small with respect to~$L$.
\end{prop}

\begin{prop}\label{p:largeonlyinS}
Suppose that $\dim L\geqslant40$. If a prime $r$ is large with respect to~$L$, then
$(r,pv|K|\cdot|\overline{G}/S|)=1$ and $k_{e(r,q)}(q)\in\omega(S)$. In particular, $t(S)\geqslant
t(L)$.
\end{prop}

The following six lemmas show that Propositions~\ref{p:structureK} and~\ref{p:structureGS} hold.

\begin{lemma}\label{l:auto1} If $r\in\pi(L)\setminus\{p\}$ and
$\varphi(r,L)>n/2$, then $r$ does not divide $|\overline{G}/S|$.
\end{lemma}

\begin{prf} Assume to the contrary that $r$ divides $|\overline{G}/S|$.

Let $L$ be a linear or unitary group first. Since $\dim L=\prk(L)=n\geqslant40$ there are integers
$\alpha,\beta\geqslant2$ and $\gamma\geqslant3$ with
$n/2<2^\alpha,3\cdot2^\beta,5\cdot2^\gamma\leqslant n$. Set $I=\{2^\alpha, 3\cdot2^\beta,
5\cdot2^\gamma\}$. Then the numbers $k_{i}(q)$, where $i\in I$, lie in $\omega(L)$ and, for every
prime divisor $s$ of any of these numbers, $\varphi(s,L)>n/2$. By hypothesis $\varphi(r,L)>n/2$,
hence Lemma~\ref{l:auxvarphi}(iii) implies that $r$ is adjacent to $s$ if and only if
$e(r,q)=e(s,q)$. Therefore, there exist at least two numbers $a$ and $b$ from $\{k_i(q)\mid i\in
I\}$ such that $rs\not\in\omega(L)$ for every prime $s$ dividing $ab$. If $s\in\pi(a)$,
$w\in\pi(b)$, then $\{r,s,w\}$ is a coclique in $GK(L)$, so both $a$ and $b$ lie in $\omega(S)$ by
Proposition~\ref{p:UniqueS}(ii). If at least one prime divisor $w$ of one of these numbers
satisfies $\varphi(w,S)\leqslant m/2$, then, for every prime divisor $s$ of the other, we have
$\varphi(s,S)>m/2$ by Lemma~\ref{l:auxvarphi}(ii). Thus, the set $\{k_i(q)\mid i\in I\}$ contains a
number $k$ such that $rs\not\in\omega(L)$ and $\varphi(s,S)>m/2$ for every prime $s$ dividing~$k$.
Lemma~\ref{l:auxvarphi}(iii) implies that $e(s,u)$ is the same for every prime divisor $s$ of~$k$.
Therefore, $k$ divides $k_{e(s,u)}(u)$. Since $\varphi(s,S)>m/2$, Lemma~\ref{l:hallkiq} yields that
$S$ has a cyclic Hall subgroup of order $k_{e(s,u)}(u)$. The $s$-exponent of $L$ equals $k_{\{s\}}$
for every $s$ dividing $k$, so $(k,k_{e(s,u)}(u)/k)=1$ and $S$ has a cyclic Hall subgroup $H$ of
order~$k$. Let $s\in\pi(k)$. The normalizer $N_{\overline{G}}(P)$ in $\overline{G}$ of a Sylow
$s$-subgroup $P$ of $S$ contains an element $x$ of order $r$ by Frattini argument. Therefore,
$H\langle x\rangle$ is a Frobenius group, and so $r$ divides $|P|-1$. Since the last inference is
valid for every prime $s\in\pi(k)$ and $H$ is cyclic of order $k$, the prime $r$ divides
$k-1=k_{i}(q)-1$  for some $i\in I$. It follows from Lemma~\ref{l:ki} that $e(r,q)$ divides
$c=2^\delta$ with $c\leqslant n/2$. Lemma~\ref{l:smallerq} implies that $\varphi(r,L)$ either
divides $c$ or equals $2$. Both possibilities contradict the hypothesis $\varphi(r,L)>n/2$.

Now let $L$ be a symplectic or orthogonal group. Since $\dim L=2\prk(L)=2n\geqslant40$, there are
integers $\alpha,\beta\geqslant2$ and $\gamma\geqslant3$ with
$n<2^\alpha,3\cdot2^\beta,5\cdot2^\gamma\leqslant 2n$. Put $I=\{2^\alpha, 3\cdot2^\beta,
5\cdot2^\gamma\}$. Exclude the case when $L=O_{2n}^+(q)$ and $2n\in I$ for a while. Then numbers
$k_{i}(q)$, where $i\in I$, lie in~$\omega(L)$ and $\varphi(s,L)=i/2>n/2$ for every prime divisor
$s$ of any of these numbers. Repeating the arguments of the preceding paragraph word for word, we
derive that $r$ divides $k_{i}(q)-1$ for some $i\in I$. Therefore, $e(r,q)$ divides $c=2^\delta$
with $c\leqslant n$ by Lemma~\ref{l:ki}. Lemma~\ref{l:smallerq} yields that $\varphi(r,L)=e(r,q)/2$
divides $c/2$, contrary to $\varphi(r,L)>n/2$. Let, finally, $L=O_{2n}^+(q)$ and $2n\in I$. We
consider the set $I'=(I\cup\{n\})\setminus\{2n\}$ instead of~$I$. The numbers $k_{i}(q)$, where
$i\in I'$, lie in $\omega(L)$. By adjacency criterion (see \cite{VasVd} or Table~\ref{tab:tL}),
$rs\not\in\omega(L)$ for every primes $s$ with $e(s,q)=n$. Therefore, there exist at least two
numbers $a$ and $b$ from $\{k_i(q)\mid i\in I'\}$ such that $rs\not\in\omega(L)$ for every prime
$s$ dividing $ab$. If $s\in\pi(a)$ and $w\in\pi(b)$, then $\{r,s,w\}$ is the coclique in $GK(L)$
(again see Table~\ref{tab:tL}), so both $a$ and $b$ lie in $\omega(S)$ by
Proposition~\ref{p:UniqueS}(ii). Repeating the preceding arguments once more, we conclude that $r$
divides $k_{i}(q)-1$ for one of $i\in I'$. Again $e(r,q)$ divides $c=2^\delta$ with $c\leqslant n$,
which is impossible because $\varphi(r,L)>n/2$.
\end{prf}

\begin{lemma}\label{l:vadj} If $r\in\pi(K)\cup\pi(\overline{G}/S)$ and $r\neq v$, then $vr\in\omega(G)$. There
exists $s\in\pi(S)\setminus(\pi(K)\cup\pi(\overline{G}/S)\cup\{v\})$ such that
$vs\not\in\omega(G)$.
\end{lemma}

\begin{prf} If $v=2$, then the first assertion of the lemma holds by
Proposition~\ref{p:UniqueS}(iii). Assume that $v$ is odd. It follows from
Proposition~\ref{p:UniqueS}(ii) that the dimension of $S$ is large enough, so $S\not\simeq L_2(v)$
and a Sylow $v$-subgroup of $S$ includes a noncyclic abelian subgroup. Therefore, by
Lemma~\ref{l:action2}, we have $rv\in\omega(G)$ for every $r\in\pi(K)\setminus\{v\}$. Furthermore,
it is well-known that the centralizer in $S$ of any outer automorphism of $S$ contains an element
of order~$v$. Thus, the first assertion of the lemma is completely proved. On the other hand,
Lemma~\ref{l:graph} yields that there is $s\in\pi(L)\setminus\{v\}$ with $sv\not\in\omega(L)$. Now
the latter assertion of the lemma follows from the former one.
\end{prf}

\begin{lemma}\label{l:parabolicaction}
If $r,s\in\pi(G)$, $r$ divides $|K|$, $s$ does not lie in $\pi(K)$ and divides the order of some
proper parabolic subgroup of~$S$, then $rs\in\omega(G)$.
\end{lemma}

\begin{prf}
We may assume that $s\neq v$ by Lemma~\ref{l:vadj}, and $s\neq2$ by
Proposition~\ref{p:UniqueS}(iii). Put $C=C_G(K)$ for the centralizer of $K$ in $G$. If
$C\not\subseteq K$, then the preimage of $S$ in $G$ lies in $CK$, and so the factor group
$C/(C\cap K)\simeq CK/K$ has a subgroup isomorphic to $S$, hence $rs\in\omega(G)$. Therefore,
$C\leqslant K$. If $S$ has an elementary abelian subgroup of order $s^2$, then $rs\in\omega(G)$ by
Lemma~\ref{l:action2}, so we may assume that there is no such subgroup in $S$, in particular, it
follows $(s,u^2-1)=1$.

Consider a normal $r$-series of $K$:
$$1=R_0\leqslant K_1<R_1\leqslant K_2\leqslant\ldots\leqslant R_{t-1}<K_t\leqslant R_t=K,$$ where
$K_i/R_{i-1}=O_{r'}(K/R_{i-1})$ and $R_i/K_i=O_r(K/K_i)$.

First we suppose that $R=R_t/K_t\neq1$. Put $V=R/\Phi(R)$ for the factor group of $R$ by its
Frattini subgroup $\Phi(R)$. It follows from Lemma~\ref{l:semidirect} that $\omega(V\rtimes
S)\subseteq\omega(\overline{G})$, where $V\rtimes S$ is the natural semidirect product
under the action of $S$ on $V$ by conjugation. By Lemma~\ref{l:good}, the group $S$ has a subgroup
$H$ which is good in $S$ with respect to~$s$. Since $C\leqslant K$, the action of $H$ on $V$ by
conjugation is faithful. Therefore, by Lemma~\ref{l:hallhigman}, there exists an element $g$ of
order $s$ in $H$ such that the minimal polynomial of $g$ under this action is equal to $x^s-1$.
Thus, $rs\in\omega(G)$ in this case.

Let $K=R_t=K_t$ and put $\widetilde{K}=K/R_{t-1}$. Assume that $v$ does not divide
$|\widetilde{K}|$. Once more we take a subgroup $H$ that is good in $S$ with respect to~$s$. By the
Schur--Zassenhaus theorem, the factor group $\widetilde{G}=G/R_{t-1}$ has a subgroup
$\widetilde{H}$ isomorphic to~$H$. Put $V=R_{t-1}/\Phi(R_{t-1})$. If $C_{\widetilde{G}}(V)$ does
not lie in $\widetilde{K}$, then $C_{\widetilde{G}}(V)\widetilde{K}/\widetilde{K}$ has a subgroup
isomorphic to~$S$, hence $s$ divides $C_{\widetilde{G}}(V)$ and $rs\in\omega(G)$. Furthermore,
$\widetilde{H}\cap\widetilde{K}=1$. So $\widetilde{H}$ acts on $V$ faithfully and we derive
$rs\in\omega(G)$ by Lemma~\ref{l:hallhigman}.

Finally, suppose that $v$ divides $|\widetilde{K}|$. Let $\widetilde{T}$ be a Sylow $v$-subgroup of
$\widetilde{K}$ and $\widetilde{G}=G/R_{t-1}$. Applying, if necessary, the Frattini argument, we
may assume that $\widetilde{T}\trianglelefteq\widetilde{G}$. Let $U$ be a minimal normal in
$\widetilde{G}$ subgroup of the center $Z(\widetilde{T})$ of $\widetilde{T}$ and put
$V=R_{t-1}/\Phi(R_{t-1})$. If an element $g$ of order $s$ from $\widetilde{G}$ centralizes $U$,
then $C_{\widetilde{G}}(U)\widetilde{K}/\widetilde{K}$ includes a subgroup isomorphic to $S$, so
$v$ is adjacent to every prime divisor of the order of $S$, which contradicts Lemma~\ref{l:vadj}.
Thus, $A=[U,g]\neq1$ and $A\langle g\rangle$ is a Frobenius group with a cyclic complement of
order~$s$. Since $C_{\widetilde{G}}(V)$ is normal in $\widetilde{G}$, either $U\leqslant
C_{\widetilde{G}}(V)$ or $C_{\widetilde{G}}(V)\cap U=1$. The former is impossible because
$K_{t-1}/R_{t-2}=O_{r'}(K/R_{t-2})$. Therefore, $A$ does not lie in $VC_{\widetilde{K}}(V)$.
Lemma~\ref{l:action1} yields $rs\in\omega(H)$.
\end{prf}

\begin{lemma}\label{l:rsnonadj}
If $r\in\pi(K)\setminus\{v\}$, then $t(r,L)=2$, and $(s,|K|\cdot|\overline{G}/S|\cdot|P|)=1$ for
every $s\in\pi(L)$ nonadjacent to $r$ in $GK(L)$ and every proper parabolic subgroup $P$ of~$S$.
\end{lemma}

\begin{prf}
Observe first that $t(r,L)\geqslant2$ due to Lemma~\ref{l:graph}. Suppose that $t(r,L)>2$ and $\rho$ is any $\{r\}$-coclique of size at least~$3$. It
follows from Proposition~\ref{p:UniqueS}(ii) that $(s,|K|\cdot|\overline{G}/S|)=1$ for every
$s\in\rho\setminus\{r\}$. By Lemma~\ref{l:parabolicaction}, all such $s$ do not
divide the order of any proper parabolic subgroup of $S$. It contradicts to
Lemma~\ref{l:adjanisotrop}.

Thus, $t(r,L)=2$. If $r=2$, then the lemma holds due to Proposition~\ref{p:UniqueS}(iii) and
Lemma~\ref{l:parabolicaction}. So we assume that $r$ is odd. Let $\{r,s\}$ be a coclique
in~$GK(L)$. Using Table~\ref{tab:trL}, we have $n/2<\varphi(s,L)\leqslant n$. Therefore, $s$ does
not divide $|\overline{G}/S|$ by Lemma~\ref{l:auto1}. Furthermore, applying Table~\ref{tab:tL}, we
obtain a coclique $\rho$ in $GK(L)$ of size $3$ with $s\in\rho$. Suppose that $s$ divides $|K|$.
Then $t$ does not divide $|\overline{G}/S|\cdot|K|$ for every $t\in\rho\setminus\{s\}$ by
Proposition~\ref{p:UniqueS}(ii). Now Lemma~\ref{l:adjanisotrop} yields that there exists a prime
$t$ from $\rho\setminus\{s\}$ dividing the order of some proper parabolic subgroup of~$S$. If
$t=v$, then $ts\in\omega(G)$ by Lemma~\ref{l:vadj}, and if $t\neq v$, then $ts\in\omega(G)$ by
Lemma~\ref{l:parabolicaction}, a contradiction with the choice of~$t$. Therefore, $s$ does not
divide $|\overline{G}/S|\cdot|K|$. Yet another application of Lemma~\ref{l:parabolicaction}
completes the proof.
\end{prf}

\begin{lemma}\label{p:KNilp} The soluble radical $K$ is nilpotent.
\end{lemma}

\begin{prf}
Assume the contrary. Then the Fitting subgroup $F=F(K)$ is a proper subgroup of~$K$. Put
$\widetilde{G}=G/F$, $\widetilde{K}=K/F$, and set $\widetilde{H}$ for the preimage of $S$
in~$\widetilde{G}$.  Let $\widetilde{T}$ be a minimal normal subgroup of~$\widetilde{G}$ lying
in~$\widetilde{K}$, and put $T$ for its preimage in~$G$. The solubility of $\widetilde{K}$ implies
that $\widetilde{T}$ is an elementary abelian $t$-group for a prime~$t$. Let
$r\in\pi(F)\setminus\{t\}$, $R$ be the Sylow $r$-subgroup of~$F$, $C_r=C_G(R)$ be the centralizer
of $R$ in~$G$, and $\widetilde{C}_r$ be the image of this centralizer in~$\widetilde{G}$. Since
$\widetilde{C}_r$ is a normal subgroup of~$\widetilde{G}$, the minimality of $\widetilde{T}$ yields
that either $\widetilde{C}_r\cap\widetilde{T}=1$ or $\widetilde{T}\leqslant\widetilde{C}_r$. If
$\widetilde{T}\leqslant\widetilde{C}_r$ for every $r\in\pi(F)\setminus\{t\}$, then $T$ is a normal
nilpotent subgroup of~$K$, a contradiction. Thus, there is a prime $r\in\pi(F)\setminus\{t\}$ with
$\widetilde{C}_r\cap\widetilde{T}=1$. If $\widetilde{K}$ does not include
$\widetilde{C}=C_{\widetilde{G}}(\widetilde{T})$, then $\widetilde{C}\widetilde{K}$ includes
$\widetilde{H}$ and $t$ is adjacent to every prime from $\pi(S)$, which is impossible due to
Lemma~\ref{l:rsnonadj}. So $\widetilde{C}\leqslant\widetilde{K}$. Lemma~\ref{l:rsnonadj} implies
that there exists $s\in\pi(S)\setminus\pi(K)$ with $rs\not\in\omega(G)$. Consider a cyclic subgroup
$\langle x\rangle$ of order $s$ in $\widetilde{H}$. Observe that $t\neq s$ because $t\in\pi(K)$.
Therefore, $\widetilde{T}=[\widetilde{T},x]\times C_{\widetilde{T}}(x)$ and
$C_{\widetilde{T}}(x)\neq\widetilde{T}$. Hence $A=[\widetilde{T},x]:\langle x\rangle$ is a
Frobenius group with a cyclic complement of order~$s$. Since $\widetilde{C}_r\cap\widetilde{T}=1$,
the action of $A$ on $R/\Phi(R)$ is faithful. By Lemma~\ref{l:action1}, we obtain $rs\in\omega(G)$
and so derive a contradiction. The lemma and Proposition~\ref{p:structureK} are proved.
\end{prf}

\begin{lemma}\label{l:auto2} If a prime $r$ not equal to $p$ divides the order of~$\overline{G}/S$,
then either $\varphi(r,L)\leqslant n/3$, or $L=L^\varepsilon_n(q)$,
$n\in[2^{\gamma+3},9\cdot2^\gamma)$, and $e(r,\varepsilon{q})=\varphi(r,L)=3\cdot2^\gamma$ for some
integer $\gamma\geqslant3$. In particular, $r$ is small with respect to~$L$.
\end{lemma}

\begin{prf}
We start with two simple observations. First, for all primes $s\in\pi(L)$ treated in the further
proof, we have $s\neq p$ and $e(s,q)>2$, so, by Lemma~\ref{l:criterion}, two such primes $s$ and
$w$ with distinct $e(s,q)$ and $e(w,q)$ are adjacent if and only if $k_{e(s,q)}(q)\cdot
k_{e(w,q)}(q)\in\omega(L)$. Due to this fact, $k_i(q)$ and $k_j(q)$ are said to be adjacent
(nonadjacent) if their prime divisors are adjacent (nonadjacent). Second, if $\varphi(s,L)>n/2$ for
$s\in\pi(L)$, then $k_{e(s,q)}(q)$ is coprime to $|K|\cdot|\overline{G}/S|$ and
$k_{e(s,q)}(q)\in\omega(S)$, hence, if $\varphi(w,L)>n/2$ and $e(w,q)\neq e(s,q)$ for $w\in\pi(L)$,
then $k_{e(s,q)}(q)\cdot k_{e(w,q)}(q)\in\omega(L)$ if and only if $k_{e(s,q)}(q)\cdot
k_{e(w,q)}(q)\in\omega(S)$. Applying Lemmas~\ref{l:auxvarphi}(iii), \ref{l:hallkiq} and arguing as
in the proof of Lemma~\ref{l:auto1}, we obtain the following: if $k_i(q)$ and $k_j(q)$ are distinct
and nonadjacent in $GK(L)$, then $S$ includes a cyclic Hall subgroup of order $k_i(q)$ or $k_j(q)$.

We will prove the lemma from the contrary. Given $r$ dividing $|\overline{G}/S|$, put
$\varphi(r,L)=l$, and observe that we may assume $n/3<l\leqslant n/2$ due to Lemma~\ref{l:auto1}.
The proof is similar to the proof of Lemma~\ref{l:auto1} but more labor-consuming. Its idea is as
follows. Suppose that we find out an integer $i$ such that $S$ includes a cyclic Hall subgroup of
order $k_i(q)$ and $r$ is nonadjacent to any prime divisor of $k_i(q)$ in~$GK(L)$. Then
$\overline{G}$ contains a Frobenius subgroup with kernel of order $k_i(q)$ and complement of order
$r$ (see proof of Lemma~\ref{l:auto1}), so $r$ divides $k_i(q)-1$. If, according to
Lemma~\ref{l:ki}, we take $i$ with $\varphi(s,L)\leqslant n/3$ for every prime divisor $s$ of
$k_i(q)-1$, then the desired contradiction is obtained.

Since the case of symplectic and orthogonal groups is slightly easier than the case of linear and
unitary groups, let us assume that $L$ is a symplectic or orthogonal group at first.

If $s\in\pi(L)$ is chosen so that $e(s,q)$ is a multiple of $4$ and $n/2<\varphi(s,L)\leqslant n$,
then $r$ and $s$ are nonadjacent provided $\varphi(r,L)+\varphi(s,L)>n$. Indeed, the adjacency
criterion for symplectic and orthogonal groups \cite[Propositions~2.4 and~2.5]{VasVdM} implies that
the adjacency of $r$ and $s$ with $\varphi(r,L)+\varphi(s,L)>n$ is possible only if $e(s,q)/e(r,q)$
is an odd integer. Therefore, $e(r,q)$ is also a multiple of $4$, and $e(r,q)=2\varphi(r,L)<n$
cannot be equal to $e(s,q)=2\varphi(s,L)>n$. So
$\varphi(s,L)/\varphi(r,L)=e(s,q)/e(r,q)\geqslant3$, which leads to the impossible chain of the
inequalities: $n\geqslant\varphi(s,L)\geqslant3\varphi(r,L)>3n/3=n$.

Since $L$ is symplectic or orthogonal, it follows that $\dim L=2n\geqslant40$. So there is an
integer $\gamma\geqslant3$ with $2n\in J=[5\cdot2^\gamma,5\cdot2^{\gamma+1})$. We partition $J$
into six intervals $J=J_1\cup\ldots\cup J_6$, putting $J_1=[5\cdot2^\gamma,11\cdot2^{\gamma-1})$,
$J_2=[11\cdot2^{\gamma-1},3\cdot2^{\gamma+1})$, $J_3=[3\cdot2^{\gamma+1},7\cdot2^{\gamma})$,
$J_4=[7\cdot2^{\gamma},2^{\gamma+3})$, $J_5=[2^{\gamma+3},9\cdot2^\gamma)$, and
$J_6=[9\cdot2^\gamma,5\cdot2^{\gamma+1})$.

Suppose that $2n\in J_1\cup J_2$. Let $a=9\cdot2^{\gamma-1}$, and let $s$ be an arbitrary prime
divisor of~$k_a(q)$. The number $a$ is divided by $4$ and satisfies the inequalities $4n/3<a<2n$.
It follows that $2n/3<\varphi(s,L)=a/2<n$, so $r$ and $s$ is nonadjacent due to the inequality
$a/2+l>2n/3+n/3=n$ and the adjacency criterion (see the observation above). The group $L$ has a
cyclic subgroup of order $k_a(q)$, and so does $S$. If such subgroup of $S$ is a Hall subgroup,
then $r$ divides $k_a(q)-1$. In this case Lemma~\ref{l:ki} yields $\varphi(r,L)\leqslant
e(r,q)\leqslant a/6<n/3$, a contradiction. Thus, we assume that $S$ has not a cyclic Hall subgroup
of order~$k_a(q)$. Let $b=2^{\gamma+2}$. The number $b$ is also divided by $4$ and satisfies
$2n/3<\varphi(s,L)=b/2<n$ for every prime divisor $s$ of $k_b(q)$. Since $k_a(q)$ are $k_b(q)$
distinct and nonadjacent, $S$ includes a cyclic Hall subgroup of order $k_b(q)$. It follows from
$b/2+l>n$ that $r$ and $k_b(q)$ are nonadjacent. Hence $r$ divides $k_b(q)-1$. By Lemma~\ref{l:ki},
$e(r,q)$ divides $2^{\gamma+1}$. If $e(r,q)\neq2^{\gamma+1}$, then
$l=e(r,q)/2\leqslant2^{\gamma-1}=3\cdot2^{\gamma-1}/3<n/3$, contrary to our assumption. Therefore,
$e(r,q)=2^{\gamma+1}$.

If $2n\in J_1$ then put $c=7\cdot2^{\gamma-1}$. Given $s\in R_c(q)$, we have $\varphi(s,L)=c/2$.
Due to the inequalities $n<c<2n$ and $c/2+l=7\cdot2^{\gamma-2}+2^\gamma=11\cdot2^{\gamma-2}>n$ the
group $S$ has a cyclic Hall subgroup of order $k_c(q)$ and $r$ is nonadjacent to every prime
divisor of~$k_c(q)$. Lemma~\ref{l:ki} implies that $e(r,q)=2^{\gamma+1}$ divides
$3\cdot2^{\gamma-1}$, a contradiction.

If $2n\in J_2$, set $c=11\cdot2^{\gamma-1}$. Then $n<c\leqslant 2n$. If $c=2n$, let $L\neq
O^+_{2n}(q)$ at first. Then we may suppose that $S$ has a cyclic Hall subgroup of order $k_c(q)$.
It follows from $c/2+l=11\cdot2^{\gamma-2}+2^\gamma=15\cdot2^{\gamma-2}>n$ that $r$ is not adjacent
to $k_c(q)$, hence $r$ divides $k_c(q)-1$. By Lemma~\ref{l:ki}, $e(r,q)=2^{\gamma+1}$ divides
$5\cdot2^{\gamma-1}$, a contradiction. If $c=2n$ and $L=O_{2n}^+(q)$, then $L$ does not contain an
element of order $k_c(q)$, which forces us to complicate a little our arguments. Consider the
interval $I=(2l+1,n)=(2^{\gamma+1}+1,11\cdot2^{\gamma-2})$ and show that it always contains a
prime, which we denote by~$w$. If $\gamma=3$ then $w=19\in I$. If $\gamma\geqslant4$, then $n>30$
and $I$ contains a prime $w$ due to the inequality $2l+1<5n/6$ and Lemma~\ref{l:intervalwithprime}.
Since $(q-1,q+1)\leqslant2$ and $w$ is odd, there is $\varepsilon\in\{+,-\}$ such that
$(w,\varepsilon{q}-1)=1$. By our assumption, a cyclic subgroup of order $k_a(q)$ is not a Hall
subgroup of~$S$. It follows that a subgroup $H$ of order $k_w(\varepsilon{q})$ must be a Hall
subgroup of~$S$. Hence $r$ divides
$k_w(\varepsilon{q})-1=\varepsilon{q}((\varepsilon{q})^{w-1}-1)/(\varepsilon{q}-1)$, so
$e(r,\varepsilon{q})=e(r,q)=2l$ divides $w-1$, which is impossible because
$e(r,q)=2^{\gamma+1}<w-1<11\cdot2^{\gamma-2}<2e(r,q)$. Thus, the lemma holds in the case $2n\in J_1\cup
J_2$.

Suppose that $2n\in J_3\cup J_4\cup J_5$. Put $a=9\cdot2^{\gamma-1}$ if $2n=3\cdot2^{\gamma+1}$,
and put $a=3\cdot2^{\gamma+1}$ otherwise. As in the previous case, $a$ is a multiple of $4$,
$4n/3<a<2n$, and $r$ cannot divide $k_a(q)-1$ by Lemma~\ref{l:ki}. There exists a cyclic subgroup
of order $k_a(q)$ in~$S$. If it is a Hall subgroup, then we derive a contradiction immediately.

If $2n\in J_3$, set $b=5\cdot2^\gamma$. Then $2n/3<b/2<n$, so $S$ includes a cyclic Hall subgroup
of order~$k_b(q)$. It follows from $b/2+l>n$ that $r$ is not adjacent to $k_b(q)$, so $r$ divides
$k_b(q)-1$. By Lemma~\ref{l:ki}, $e(r,q)$ divides $2^{\gamma+1}$, and the condition
$e(r,q)/2=l>n/3$ yields $e(r,q)=2^{\gamma+1}$. Let $c=11\cdot2^{\gamma-1}$. Following the same way,
we obtain that $r$ divides $k_c(q)-1$. Applying Lemma~\ref{l:ki}, we conclude that
$e(r,q)=2^{\gamma+1}$ must divide $5\cdot2^{\gamma-1}$, which is impossible.

If $2n\in J_4$, put $b=7\cdot2^\gamma$ and $c=11\cdot2^{\gamma-1}$. Then $S$ has cyclic Hall
subgroups of orders~$k_b(q)$ (except the case: $2n=b$ and $L=O_{2n}^+(q)$, treated separately)
and~$k_c(q)$. It follows from $b/2+l>n$ that $r$ divides $k_b(q)-1$, so $e(r,q)=3\cdot2^{\gamma}$.
Then $c/2+l=11\cdot2^{\gamma-2}+3\cdot2^{\gamma-1}>2^{\gamma+2}>n$, hence $e(r,q)$ divides
$5\cdot2^{\gamma-1}$, a contradiction. Let $b=2n$ and $L=O_{b}^+(q)$. The existence of a cyclic
Hall subgroup of order $k_c(q)$ in $S$ yields $e(r,q)=2l=5\cdot2^{\gamma-1}$. There is a prime $w$
in $I=(2l+1,n)=(5\cdot2^{\gamma-1}+1,7\cdot2^\gamma)$. Indeed, one can put $w=23$ for $\gamma=3$,
and apply Lemma~\ref{l:intervalwithprime} for $\gamma\geqslant4$. Now we choose
$\varepsilon\in\{+,-\}$ so that $(w,\varepsilon{q}-1)=1$. The group $S$ has a cyclic Hall subgroup
of order $k_w(\varepsilon{q})$, so $r$ divides $k_w(\varepsilon{q})-1$. Therefore, $e(r,q)$ divides
$w-1$, which is impossible because $e(r,q)<w-1<2e(r,q)$.

If $2n\in J_5$, put $b=2^{\gamma+3}$ and $c=7\cdot2^{\gamma}$. If $L\neq O_{b}^+(q)$, then there
are cyclic Hall subgroups of $S$ having orders $k_b(q)$ and~$k_c(q)$. Then $r$ divides $k_b(q)-1$,
so $e(r,q)=2^{\gamma+2}$, and $r$ and $k_c(q)$ are nonadjacent due to $c/2+l>n$. However,
$2^{\gamma+2}$ does not divide $3\cdot2^\gamma$, a contradiction. Let $L=O_{b}^+(q)$. The existence
of a cyclic Hall subgroup of order $k_c(q)$ provides $e(r,q)=2l=3\cdot2^\gamma$. Since
$n=b/2=2^{\gamma+2}>30$, there is a prime $w$ in $I=(2l+1,n)$ by Lemma~\ref{l:intervalwithprime}.
Arguing as in the previous case, we obtain that $e(r,q)$ divides $w-1$ and so derive a
contradiction.

Suppose, finally, that $2n\in J_6$. Put $a=3\cdot2^{\gamma+1}$ for $2n=9\cdot2^\gamma$, and
$a=9\cdot2^\gamma$ otherwise. By the choice of $a$, we must assume that a cyclic subgroup of order
$k_a(q)$ in $S$ is not a Hall subgroup in order to avoid the immediate contradiction (similarly to
the previous cases). Putting $b=2^{\gamma+3}$, we obtain that $S$ has a cyclic Hall subgroup of
order $k_b(q)$, and derive $e(r,q)=2^{\gamma+2}$. Handling a cyclic Hall subgroup of order
$k_c(q)$, where $c=7\cdot2^\gamma$, and the inequality $c/2+l=11\cdot2^\gamma>n$, we conclude that
$e(r,q)$ divides $3\cdot2^\gamma$, a contradiction. Thus, the lemma is proved for symplectic and
orthogonal groups.

Let $L=L_n^\varepsilon(q)$. It follows that $l=\varphi(r,L)=e(r,\varepsilon q)$. If $s\in\pi(L)$ is
chosen so that $e(s,q)$ is a multiple of $4$, then $e(s,q)=e(s,-q)=\varphi(s,L)$. Unfortunately,
unlike the case of symplectic and orthogonal groups, the inequalities $n/2<\varphi(s,L)\leqslant n$
and $l+\varphi(s,L)>n$ do not guarantee that $r$ and $s$ are nonadjacent, because $\varphi(s,L)/l$
can be an integer greater than $1$ (see \cite[Lemma~2.1]{VasGrSt}).

By the hypothesis $\dim L=n\geqslant40$, so there is an integer $\gamma\geqslant3$ with $n\in
J=[5\cdot2^\gamma,5\cdot2^{\gamma+1})$. We consider the same partition $J=J_1\cup\ldots\cup J_6$ as
in the case of symplectic and orthogonal groups.

Suppose that $n\in J_1\cup J_2$. Let $a=9\cdot2^{\gamma-1}$ and $b=2^{\gamma+2}$. The numbers $a$
and $b$ are divided by $4$, and the inequalities $2n/3<b<a<n$ hold. Therefore, due to the adjacency
criterion \cite[Lemma~2.1]{VasGrSt}, $r$ is adjacent to $k_a(q)$ ($k_b(q)$ respectively) if and
only if $l$ divides $a$ ($l$ divides $b$). Assume that $l$ does not divide~$a$. The group $S$ has a
cyclic subgroup $H$ of order $k_a(q)$. If $H$ is a Hall subgroup, then $r$ divides $k_a(q)-1$, so
$e(r,q)$ divides $3\cdot2^{\gamma-2}$ by Lemma~\ref{l:ki}. Then $l\leqslant3\cdot2^{\gamma-2}<n/3$,
a contradiction. If $H$ is not a Hall subgroup, then a cyclic subgroup of order $k_b(q)$ in $S$
must be Hall due to the nonadjacency of $k_a(q)$ and~$k_b(q)$. If $r$ is adjacent to $k_b(q)$ then
$l$ divides $b$. The inequalities $l>n/3$ and $b<n$ yield $l=b/2=2^{\gamma+1}$. Then $l$ is a
multiple of $4$, and so $e(r,q)=l=2^{\gamma+1}$. If $r$ is nonadjacent to $k_b(q)$, then $r$
divides $k_b(q)-1$, and Lemma~\ref{l:ki} implies that $e(r,q)$ divides $2^{\gamma+1}$. It follows
from the inequality $l>n/3$ that $e(r,q)=2^{\gamma+1}$. Assume, finally, that $l$ divides $a$, then
$l=9\cdot2^{\gamma-2}$ because $n/3<l\leqslant n/2$. Thus, the following alternative holds: either
$l=e(r,q)=2^{\gamma+1}$ or $l=9\cdot2^{\gamma-2}$. Furthermore, subgroups of orders $k_a(q)$ and
$k_b(q)$ in $S$ cannot be Hall simultaneously, otherwise $2^{\gamma+1}=l=9\cdot2^{\gamma-2}$.

If $n\in J_1$, put $c=7\cdot2^{\gamma-1}$. For both of alternative values, $l$ does not divide $c$
and $c+l>n$. Indeed,
$7\cdot2^{\gamma-1}+9\cdot2^{\gamma-2}>7\cdot2^{\gamma-1}+2^{\gamma+1}=11\cdot2^{\gamma-1}>n$.
Therefore, $r$ is nonadjacent to~$k_c(q)$. Furthermore, $S$ includes a cyclic Hall subgroup of
order $k_c(q)$, hence $r$ divides $k_c(q)-1$. It follows that $e(r,q)$ divides
$3\cdot2^{\gamma-1}$, which is impossible.

If $n\in J_2$, put $c=11\cdot2^{\gamma-1}$. Since $l$ does not divide $c$ and the inequalities
$11\cdot2^{\gamma-1}+9\cdot2^{\gamma-2}>11\cdot2^{\gamma-1}+2^{\gamma+1}=15\cdot2^{\gamma-1}>n$
hold, $r$ is nonadjacent to~$k_c(q)$. Therefore, $r$ divides $k_c(q)-1$. Due to Lemma~\ref{l:ki},
$e(r,q)$ must divide $5\cdot2^{\gamma-1}$, a contradiction. This completes the proof for $n\in
J_1\cup J_2$.

Let $n\in J_3\cup J_4\cup J_5$. Put $a=3\cdot2^{\gamma+1}$ and observe that $2n/3<a\leqslant n$. If
$r$ divides $k_a(q)-1$, then $e(r,q)$ divides $2^{\gamma+1}\leqslant n/3$ by Lemma~\ref{l:ki}.
Therefore, similarly to the previous case,  we can avoid the immediate contradiction just in two
cases: either a cyclic subgroup $H$ of order $k_a(q)$ is not a Hall subgroup of $S$ or
$l=e(r,q)=3\cdot2^\gamma$.

Let $n\in J_3$ and $b=5\cdot2^\gamma$. If $H$ is not a Hall subgroup, then a cyclic subgroup of
order $k_b(q)$ must be a Hall subgroup of~$S$. Therefore, either $l$ divides $b$ and so
$l=b/2=5\cdot2^{\gamma-1}$ or $l$ divides $k_b(q)-1$ and so $l=e(r,q)=2^{\gamma+1}$. Thus, one of
two subgroups of orders $k_a(q)$ and $k_b(q)$ is not a Hall subgroup of $S$, and $l$ possesses one
of the following values: $3\cdot2^\gamma,5\cdot2^{\gamma-1},2^{\gamma+1}$. If
$l\neq5\cdot2^{\gamma-1}$, put $c=11\cdot2^{\gamma-1}$. It follows from
$11\cdot2^{\gamma-1}+3\cdot2^\gamma>11\cdot2^{\gamma-1}+2^{\gamma+1}=14\cdot2^{\gamma-1}>n$ that
$r$ is nonadjacent to~$k_c(q)$, so $r$ divides $k_c(q)-1$. Then $l=e(r,q)$ divides
$5\cdot2^{\gamma-1}$, contrary to our assumptions. If $l=5\cdot2^{\gamma-1}$, put
$c=9\cdot2^{\gamma-1}$. Then $c+l>n$. Hence $r$ divides $k_c(q)-1$ and $l$ divides
$3\cdot2^{\gamma-2}$, which is impossible.

Let $n\in J_4$ and $b=11\cdot2^{\gamma-1}$. If $H$ is not a Hall subgroup, then a cyclic subgroup
of order $k_b(q)$ is a Hall subgroup of $S$ and, arguing as in the previous paragraph, we obtain
that either $l$ divides $b$ and $l=11\cdot2^{\gamma-2}$, or $r$ divides $k_b(q)-1$ and
$l=e(r,q)=5\cdot2^{\gamma-1}$. Thus, at least one of cyclic subgroups of orders $k_a(q)$ and
$k_b(q)$ is not Hall in $S$, and $l\in\{3\cdot2^\gamma,5\cdot2^{\gamma-1},11^{\gamma-2}\}$. If
$l=3\cdot2^\gamma$, put $c=5\cdot2^\gamma$, otherwise put $c=7\cdot2^\gamma$. The verification
analogous to the one in the previous paragraph leads to a contradiction for every possible value
of~$l$.

Let $n\in J_5$ and $b=2^{\gamma+3}$. If $H$ is not a Hall subgroup, then either $l$ divides $b$ or
$r$ divides $k_b(q)-1$. In both cases, $l=e(r,q)=2^{\gamma+2}$. Thus, one of two subgroups of
orders $k_a(q)$ and $k_b(q)$ is not a Hall subgroup of $S$, and
$l\in\{3\cdot2^\gamma,2^{\gamma+2}\}$. If $l=2^{\gamma+2}$, then $r$ is nonadjacent to $k_c(q)$,
where $c=7\cdot2^\gamma$. So $r$ divides $k_c(q)-1$, together with Lemma~\ref{l:ki} this leads to
a contradiction due to $l$ does not divide $3\cdot2^\gamma$. Therefore, for $n\in J_5$ the number
$l=e(r,q)$ must be equal to $3\cdot2^\gamma$. In this case $l<n/2-1$, so
Lemma~\ref{l:estimlarge}(ii) yields the conclusion of the lemma for $n\in J_5$.

Let $n\in J_6$, $a=9\cdot2^\gamma$, $b=2^{\gamma+3}$, and $c=7\cdot2^\gamma$. By the choice of $a$,
if a cyclic subgroup of order $k_a(q)$ is a Hall subgroup of $S$, then $l=9\cdot2^{\gamma-1}$.
Similarly, if a cyclic subgroup of order $k_b(q)$ is a Hall subgroup of $S$, then $l=2^{\gamma+2}$.
Since $l$ cannot be equal to $9\cdot2^{\gamma-1}$ and~$2^{\gamma+2}$ simultaneously, one of these
subgroups is not Hall. Therefore, $S$ has a cyclic Hall subgroup of order~$k_c(q)$. It follows from
$9\cdot2^{\gamma-1}+7\cdot2^\gamma>2^{\gamma+2}=11\cdot2^\gamma>n$ that $r$ is nonadjacent
to~$k_c(q)$. Then $l$ divides $3\cdot2^\gamma$ by Lemma~\ref{l:ki}, a contradiction. The lemma and
Proposition~\ref{p:structureGS} are proved.
\end{prf}

It remains to prove Proposition 5.

Since $\dim L\geqslant40$, it follows from Table~\ref{tab:tL} that $t(L)\geqslant15$. Suppose that
$r$ is large with respect to $L$. By Table~\ref{tab:tpL}, we have $t(p,L)\leqslant4$ and
$t(v,S)\leqslant4$. Lemma~\ref{l:vadj} yields that $t(v,L)=t(v,G)\leqslant t(v,S)\leqslant4$. Hence
$p\neq r\neq v$. Propositions~\ref{p:structureK} and~\ref{p:structureGS} provide
$(r,|K|\cdot|\overline{G}/S|)=1$. Applying Table~\ref{tab:trL} and \cite[Table~6]{VasVd}, we obtain
that $r\not\in\delta(L)$. Therefore, by Lemma~\ref{l:indepofcrit}, every prime $w\in R_{e(r,q)}(q)$
is large with respect to $L$, so $(k_{e(r,q)}(q),|K|\cdot|\overline{G}/S|)=1$. On the other hand,
$k_{e(r,q)}(q)\in\omega(L)=\omega(G)$. Hence $k_{e(r,q)}(q)\in\omega(S)$. Finally, if $\rho$ is a
coclique of size $t(L)$ in $GK(L)$, then every $r\in\rho$ is large with respect to $L$. It follows
that $\rho\subseteq\pi(S)$. If $\rho$ is not a coclique in $GK(S)$, then it is not a coclique in
$GK(G)=GK(L)$, which is impossible. Thus, $t(S)\geqslant t(L)$. Proposition 5 is proved.

\section{Proof: characteristic 2}\label{s:Char2}

Here we prove Theorem~\ref{t:main} provided the characteristic $p$ of $L$ equals~$2$. Since we
apply Proposition~\ref{p:largeonlyinS}, we preserve the condition on the dimension of $L$ from the
previous section.

\begin{prop}\label{p:solutionfor2}
Suppose that $\dim L\geqslant40$ and $p=2$. Then $S$ cannot be a group of Lie type over the field
of odd characteristic.
\end{prop}

\begin{prf}
Given a finite group $H$, consider cocliques $\rho$ of $GK(H)$ such that $4r\not\in\omega(H)$ for
every $r\in\rho$. Choose among them a coclique of greatest size and denote it by $\rho^*(4,H)$. Put
$t^*(4,H)=|\rho^*(4,H)|$.

\begin{lemma}\label{l:t4L} Every prime lying in $\rho^*(4,L)$ is large with respect to~$L$, and
$t^*(4,L)\geqslant3$.
\end{lemma}

\begin{prf} Apply \cite{ButLU,ButSO}.
\end{prf}

\begin{lemma}\label{l:t4S} $t^*(4,S)\leqslant2$.
\end{lemma}

\begin{prf} Again apply \cite{ButLU,ButSO}, keeping in mind that $u^2\equiv1\pmod8$ for odd~$u$.
\end{prf}

We are ready to complete the proof of Proposition~\ref{p:solutionfor2}. Consider a coclique
$\rho=\rho^*(4,L)$. By Lemma~\ref{l:t4L}, all primes from $\rho$ are large with respect to~$L$.
Therefore, $\rho\subseteq\pi(S)$ due to Proposition~\ref{p:largeonlyinS}. It follows from
Lemmas~\ref{l:t4L} and~\ref{l:t4S} that $t^*(4,L)>2\geqslant t^*(4,S)$. Thus, for at least one
$r\in\rho$, we have $4r\in\omega(S)\subseteq\omega(G)=\omega(L)$, a contradiction.
\end{prf}

\begin{rem}
If $L=L_n^\varepsilon(q)$ and $q$ is even, then the conclusion of Theorem~\ref{t:main} has already
been obtained under much more weaker hypothesis (see \cite{GrLin, VasGrLs} for linear and
\cite{GrUQ,GrShUn} for unitary groups). Moreover, if $L$ is a finite simple linear or unitary group
over a field of even order, $L\not\in\{U_4(2),U_5(2)\}$, and $G$ is a finite group with $\omega(G)=\omega(L)$,
then $L\leqslant G\leqslant\Aut L$, and all such groups $G$ are determined for every given~$L$.
\end{rem}

\section{Proof: independence and $p$-independence numbers of $S$}\label{s:tSeqtL}

Here we prove that $t(L)=t(S)$ provided $t(L)\geqslant 23$. Using this equality, we eliminate the
exceptional case of Proposition~\ref{p:structureGS}, i.e., we establish that every prime divisor
$r$ of the order of $L$ with $\varphi(r,L)>n/3$ does not divide $|K|\cdot|\overline{G}/S|$. In
conclusion, we show that under some additional assumptions the $p$-independence numbers of $L$ and
$S$ coincide as well.

\begin{prop}\label{p:tSlesstLp2}
Suppose that $t(L)\geqslant23$ and for some positive integer $a$ the number $k_a(u)$ has a prime
divisor large with respect to~$S$. Then $k_a(u)$ has a prime divisor large with respect to~$L$. In
particular, $t(L)=t(S)$ and every prime $r$ large with respect to~$L$ is large with respect to~$S$.
\end{prop}

\begin{prf} Suppose to the contrary that none of prime divisors of $k_a(u)$ is large with respect to~$L$.

\begin{lemma}\label{l:setJ}
There exists a set $J$ of positive integers of size $d=\max\{1,t(S)-t(L)\}$, satisfying the following:

\emph{(i)} for any $j\in J$, every $r\in R_j(u)$ is large with respect to $S$ and divides the
number $|\overline{G}/S|\cdot|K|$;

\emph{(ii)} if $t(S)>t(L)$ and every coclique $\rho$ of the greatest size in $GK(L)$ contains a
prime $s$ with $\varphi(s,S)\leqslant m/2$, then $\varphi(r,S)>m/2$ for any $j\in J$ and every
$r\in R_j(u)$.
\end{lemma}

\begin{prf}
First we prove that there is a set $J$ satisfying~(i). Put $t=t(S)$ and $t(L)=l$. There exists a
set $I$ of positive integers of size $t$ containing $a$ such that for any $i\in I$ every prime from
$R_i(u)$ is large with respect to~$S$. Assume that there are $l+1$ numbers $i\in I$ such that
$\widetilde{R}_i(u)=R_i(u)\setminus(\pi(K)\cup\pi(\overline{G}/S))\neq\varnothing$ and put $\rho$
for a set consisting of $l+1$ primes from different $\widetilde{R}_i(u)$. Then $\rho$ is a coclique
in $GK(L)$ of size $l+1$, a contradiction. Thus, if $t>l$, then there is a subset $J$ of $I$ with
$|J|=t-l\geqslant1$ and such that for any $j\in J$ every prime from $R_j(u)$ divides
$|K|\cdot|\overline{G}/S|$. Let $t=l$ (one may observe that $t(L)$ does not exceed $t(S)$ due to
Proposition~\ref{p:largeonlyinS}). If for every $i\in I$ we have
$\widetilde{R}_i(u)=R_i(u)\setminus(\pi(K)\cup\pi(\overline{G}/S))\neq\varnothing$, then a set
$\rho$ consisting of $l$ primes from different $\widetilde{R}_i(u)$ forms a coclique of greatest
size in~$GK(L)$. Therefore all primes from $\rho$ are large with respect to~$L$. On the other hand,
$\rho$ contains a prime from $R_a(u)$, which is impossible due to our assumption on primes
dividing~$k_a(u)$.

Now we show that a set $J$ can be chosen to satisfy (ii) as well. It follows from
Lemma~\ref{l:auxvarphi}(ii) that the set $I$ includes a subset $I'$ of size at least $t-1$ such
that for any $i\in I'$ and every prime $w\in R_i(u)$ the inequality $\varphi(w,S)>m/2$ holds. If
$|I'|=|I|=t$, then $J\subseteq I=I'$, as required. So we may assume that $|I'|=t-1$. Suppose that
there are $l$ numbers from $I'$ with
$\widetilde{R}_i(u)=R_i(u)\setminus(\pi(K)\cup\pi(\overline{G}/S))\neq\varnothing$. Then a set
$\rho$ consisting of $l$ primes from different $\widetilde{R}_i(u)$ forms a coclique in~$GK(L)$.
However, if the conditions of (ii) hold, then there are at most $l-1$ such numbers in $I'$. Thus,
there exists a subset $J$ of $I'$ such that $|J|=t-1-(l-1)=t-l\geqslant1$, and for any $j\in J$
every prime $r$ from $R_j(u)$ is large with respect to $S$, divides~$|\overline{G}/S|\cdot|K|$, and
satisfies $\varphi(r,S)>m/2$. The lemma is proved.
\end{prf}

\begin{lemma}\label{l:rsinLnotinS} Let a set $J$ be as in Lemma~\emph{\ref{l:setJ}}. For each $j\in J$ and every
$r\in R_j(u)$ there is a large with respect to $L$ prime $s$ with
$rs\in\omega(L)\setminus\omega(S)$.
\end{lemma}

\begin{prf}
Fix $j\in J$ and $r\in R_j(u)$. Let $\rho=\{s_1,\ldots,s_l\}$ be a coclique of greatest size
in~$GK(L)$. By Proposition~\ref{p:largeonlyinS}, every prime from $\rho$ does not divide
$|\overline{G}/S|\cdot|K|$, and is not equal to~$v$. It follows that the set
$I=E(\rho,S)=\{e(s,u)\mid s\in\rho\}$ is well defined (see Section~\ref{s:primegraph}), has the
size $l$, and $j=e(r,u)\not\in I$. Choose a coclique $\sigma$ of size $t$ in~$GK(S)$ containing $r$
and put $Y=\{e(w,u)\mid w\in\sigma\}$.

If $t=l$, then $\rho$ is also a coclique of greatest size in~$GK(S)$. Lemma~\ref{l:tL} yields
$I\cap Y=Y\setminus\{j\}$. Therefore, $\rho$ includes a subset $\rho'$ of size $l-1$ such that the
set $\{r\}\cup\rho'$ is a coclique in~$GK(S)$. Since the size of this set equals $l$ and $r$ is
small with respect to~$L$, this set cannot be a coclique in~$GK(L)$. Hence there is $s\in\rho$ with
$rs\in\omega(L)\setminus\omega(S)$. Thus we may assume that $t>l$.

Suppose that $\varphi(r,S)>m/2$. By Lemma~\ref{l:auxvarphi}(ii), the set $\rho$ contains a subset
$\rho'$ of size $l-1$ with $\varphi(s,S)>m/2$ for any $s\in\rho'$. Since $j\not\in I$,
Lemma~\ref{l:auxvarphi}(iii) implies that $\{r\}\cup\rho'$ is a coclique in~$GK(S)$, and is not
in~$GK(L)$. Hence the assertion of the lemma again holds.

Let, finally, $\varphi(r,S)\leqslant m/2$. According to Lemma~\ref{l:setJ}(ii), a set $\rho$ can be
chosen in a way that $\varphi(s,S)>m/2$ for any $s\in\rho$.  Then, by Lemma~\ref{l:estimlarge}(iv),
there is a coclique $\sigma'$ of greatest size in $GK(S)$ with $\rho\subseteq\sigma'$. Set
$X=\{e(w,u)\mid w\in\sigma'\}$. Applying Lemma~\ref{l:tL}, we obtain that $Y\cap X\supseteq
Y\setminus\{j\}$. Therefore, $\rho$ includes a subset $\rho'$ of size $l-1$ such that
$\{r\}\cup\rho'$ is a coclique in~$GK(S)$, and is not in~$GK(L)$. The lemma is proved.
\end{prf}

\begin{lemma}\label{l:estimphirsbyn}
Let a prime $r$ be large with respect to~$S$. If $S$ is linear or unitary, then $\varphi(r,S)\geqslant t(L)\geqslant n/2$ and $r\geqslant n/2+1$.
If $S$ is symplectic or orthogonal, then $\varphi(r,S)\geqslant (2t(L)-4)/3\geqslant (n-4)/3$ and $r\geqslant(2n-5)/3>n/2$.
\end{lemma}

\begin{prf} By Lemma~\ref{l:estimlarge}(iii), we have $\varphi(r,S)\geqslant t(S)\geqslant
t(L)\geqslant n/2$, if $S$ is linear or unitary, and $\varphi(r,S)\geqslant (2t(S)-4)/3\geqslant
(2t(L)-4)/3\geqslant(n-4)/3$, if $S$ is symplectic or orthogonal. Lemma~\ref{l:fermat} yields
$r\geqslant\varphi(r,S)+1\geqslant n/2+1$ in the first case, and
$r\geqslant2\varphi(r,S)+1\geqslant(2n-5)/3$ in the second one.
\end{prf}

\begin{lemma}\label{l:estimkju}
Let a set $J$ be as in Lemma~\emph{\ref{l:setJ}}. Then for each $j\in J$ and every prime $r$ from
$R_j(u)$ the number $(k_j(u))_{\{r\}}$ divides $p(q^2-1)\log_vu$. Moreover, the inequality
$$\frac{\prod_{j\in J}k_j(u)}{\log_vu}\leqslant p(q^2-1)$$ holds true.
\end{lemma}

\begin{prf} Let us fix $j\in J$ and a prime $r\in R_j(u)$. Put $(k_j(u))_{\{r\}}=r^\gamma$.
Taking into account that $r$ is large with respect to $S$ and applying Lemma~\ref{l:rpart}, one can
observe that $r>3$, $r\neq v$, and $r^\gamma$ is the $r$-exponent of~$S$.

Suppose that $r$ divides $|K|$. Since $r\neq v$, Proposition~\ref{p:structureK} provides a prime
$s\in\pi(S)$ nonadjacent to~$r$ in $GK(L)$ and coprime to the order of any parabolic subgroup of~$S$. It
follows from Lemma~\ref{l:adjanisotrop} that $r$ divides the order of some parabolic subgroup $P$
of~$S$ and, by the same lemma,  $r^\gamma\in\omega(P)$. If $R$ is the Sylow $r$-subgroup of $K$
(recall that $K$ is nilpotent by Proposition~\ref{p:structureK}), then $S$ acts on $R/\Phi(R)$ by
conjugation.  This action must be faithful because $r$ and $s$ are nonadjacent. Furthermore, since
$r$ is large with respect to $S$, we have $(r,6u(u^2-1))=1$. Applying Lemma~\ref{l:ZalCor}, we
obtain $r^{\gamma+1}\in\omega(G)=\omega(L)$. On the other hand, Proposition~\ref{p:structureK}
yields that $t(r,L)=2$, and so $r$ either equals to $p$ or divides $q^2-1$ by Lemmas~\ref{l:tpL}
and~\ref{l:trLeq2&3&4}.

Lemma~\ref{l:estimphirsbyn} provides the inequality $r>n/2$. Suppose $r=p$. Since $n^2/4>2n-1$ whenever $n\geqslant8$, we have $p^2>2n-1$. It follows from Lemma~\ref{l:bigk}(iv) that the $p$-exponent of $L$ does not exceed $p^2$. Then $p^2\geqslant p^{\gamma+1}$, whence $k_j(u)_{\{p\}}\leqslant p$. Suppose that $r$ divides $(q^2-1)$ and put $(q^2-1)_{\{r\}}=r^\delta$. Using the inequality $r>n/2$ and Lemma~\ref{l:rpart}, we obtain that the $r$-exponent of $L$ does not exceed $r^{\delta+1}$. Therefore, $r^{\delta+1}\geqslant r^{\gamma+1}$. Thus, for any $j\in J$  and every $r\in R_j(u)\cap\pi(K)$ the number $(k_j(u))_{\{r\}}$ divides $p(q^2-1)$.

Suppose now that $r$ does not divide~$|K|$. Then $|\overline{G}/S|_{\{r\}}=r^\kappa>1$. Therefore,
$\overline{G}$ includes a subgroup isomorphic to an extension of $S$ by an automorphism $\tau$ of
order~$r^\kappa$, where $\kappa\geqslant1$. Since $r$ is odd and coprime to
$|\operatorname{Inndiag}(S)/S|$, we may assume that $\tau$ is a field automorphism. If $u=v^\beta$
and $\beta=r^\nu\cdot l$, where $(r,l)=1$, then $\nu\geqslant\kappa\geqslant1$. If $r$ does not
divide $v^{lj}-1$, then $r$ does not divide $v^{r^\nu\cdot lj}-1=u^j-1$, which is false. So
$r^\gamma=(k_j(u))_{\{r\}}=(u^j-1)_{\{r\}}=(v^{r^\nu\cdot
lj}-1)_{\{r\}}=r^\nu(v^{lj}-1)_{\{r\}}>r^\kappa$ by Lemma~\ref{l:rpart}. Further, $r^\gamma$ is the
greatest power of~$r$ lying in~$\omega(S)$. Lemma~\ref{l:AutoOrderPreserve} yields that $r^\gamma$
is the $r$-exponent of~$\overline{G}$, so it is the $r$-exponent of~$G$. If $r\neq p$ and
$k=e(r,q)\geqslant3$, then the inequality $r>n/2$ implies that $(q^k-1)_{\{r\}}$ is the
$r$-exponent of~$L$. By Lemma~\ref{l:rsinLnotinS}, there is a large with respect to $L$ prime $s$
with $rs\in\omega(L)\setminus\omega(S)$. It follows from Lemma~\ref{l:criterion} that
$r^\gamma{s}\in\omega(L)$ as well. Since $(rs,|K|)=1$, the group $\overline{G}$ contains an element
$x$ of order $r^\gamma{s}$. Then the element $y=x^{r^{\kappa}}$ is of order $r^{\gamma-\kappa}{s}$
and belongs to~$S$, which is impossible because $\gamma>\kappa$. Thus, $r$ divides $p(q^2-1)$.

If $r$ divides $q^2-1$, then again the inequality $r>n/2$ and Lemma~\ref{l:rpart} implies
$r^\gamma\leqslant r(q^2-1)_{\{r\}}$. Suppose $r=p$. Then the $p$-exponent of $L$ does not
exceed~$p^2$ because $p>n/2$. Since the product of distinct primes from $R_j(u)$, dividing
$|\overline{G}/S|$, divides the number $\log_vu$, we obtain that $k_j(u)$ divides
$p(q^2-1)\log_vu$. Finally, since for all $j\in J$ numbers $k_j(u)$ are pairwise coprime, the
product $\prod_{j\in J}k_j(u)$ divides $p(q^2-1)\log_vu$.
\end{prf}

\begin{lemma}\label{l:qbyq} Let $d$ be a size of the set $J$ from Lemma~\emph{\ref{l:setJ}}. Then
\begin{equation}\label{eq:qbyq} p(q^2-1)> q^{\frac{10t(L)d}{3t(S)}}. \end{equation}
\end{lemma}

\begin{prf}
Fix $j\in J$ and let $r$ be a prime from~$R_j(u)$. If $S=L_m^\varepsilon(u)$, then
$k_j(u)=k_{\varphi(r,S)}(\varepsilon u)$ and, by Lemma~\ref{l:estimphirsbyn}, the inequality
$\varphi(r,S)\geqslant t(L)\geqslant 23$ holds. If $S$ is symplectic or orthogonal,
then $\varphi(r,S)=\eta(j)$ and $\eta(j)\geqslant
(2t(L)-4)/3\geqslant 13$ by Lemma~\ref{l:estimphirsbyn}. Anyway, Lemma~\ref{l:kjubyu7} yields that
$$k_j(u)> u^6\log_vu.$$ This equation and Lemma~\ref{l:estimkju} imply
\begin{equation}\label{eq:u6d} p(q^2-1)>u^{6d}.\end{equation}
By Lemma~\ref{l:bigk}(i)--(ii), the spectrum of $L$ contains a number $b$ such that $b\geqslant q^{10t(L)/9}$ and all its prime divisors are large with respect to~$L$. It follows from Proposition~\ref{p:largeonlyinS} that $b\in\omega(S)$, so $b\leqslant u^{2t(S)}$ by Lemma~\ref{l:bigk}(iii). Thus,
\begin{equation}\label{eq:ubyq} u^{2t(S)}\geqslant q^{10t(L)/9}. \end{equation}
Finally, (\ref{eq:u6d}) and (\ref{eq:ubyq}) yield the inequality~(\ref{eq:qbyq}). The lemma is proved.
\end{prf}

Now we show that the right side of~(\ref{eq:qbyq}) is greater than $q^3$, which leads us to a contradiction with the conclusion of Lemma~\ref{l:qbyq}. Recall that $t(L)=l$ and $t(S)=t$.

Assume that $l/t>9/10$. Then $$\frac{10ld}{3t}\geqslant \frac{10l}{3t}>\frac{10}{3}\cdot \frac{9}{10}=3.$$

Let now $l/t\leqslant 9/10$. Then $$\frac{10ld}{3t}=\frac{10l(t-l)}{3t} = \frac{10l}{3}\cdot \left(1-\frac{l}{t}\right)\geqslant \frac{10l}{3}\cdot \frac{1}{10}=\frac{l}{3}>3.$$

Thus, we got a contradiction. It follows that $k_a(u)$ has a prime divisor large with respect to~$L$, so $t(S)\leqslant t(L)$. An application of Proposition~\ref{p:largeonlyinS} completes the proof of Proposition~\ref{p:tSlesstLp2}.
\end{prf}

\begin{prop}\label{p:structureGS2}
Suppose that $t(L)\geqslant23$. If a prime $r$ distinct from $p$ divides the order of
$\overline{G}/S$, then $\varphi(r,L)\leqslant n/3$.
\end{prop}

\begin{prf}
It follows from Proposition~\ref{p:structureGS} that it remains to prove the assertion: if
$L=L^\varepsilon_n(q)$, $n\in[2^{\gamma+3},9\cdot2^\gamma)$ and
$e(r,\varepsilon{q})=\varphi(r,L)=3\cdot2^\gamma$ for some $\gamma\geqslant3$, then $r$ does not divide~$|\overline{G}/S|$.

Assume to the contrary, that $r$ divides $|\overline{G}/S|$. Let $b=2^{\gamma+3}$ and $l=e(r,q)=3\cdot2^\gamma$. Since $b$ is a multiple of~$4$, the equality $k_b(q)=k_b(\varepsilon{q})$ holds. If $s\in R_b(q)$, then $s$ is large with respect to~$L$ and, by Proposition~\ref{p:tSlesstLp2}, it is also large with respect to~$S$. Therefore, $\varphi(s,S)\geqslant m/2-1$ by Lemma~\ref{l:estimlarge}(ii).
Denote the $s$-part of $k_b(q)$ by $f$, and a Sylow $s$-subgroup of $S$ by~$F$.
Then $|F|\in\{f,f^2\}$. By Frattini argument, the group $N=N_{\overline{G}}(F)$ contains an element
$x$ of order~$r$. Since $b+l>n$ and $l$ does not divide $b$, the primes $r$ and $s$ are nonadjacent in~$GK(L)$. Therefore,
$F\rtimes\langle x\rangle$ is a Frobenius group. It follows that $r$ divides $f^2-1$. This is true for every prime divisor $s$ of $k_b(q)$, so the number $k_b(q)^2-1$ is a multiple of~$r$.
Further, $k_b(q)^2-1=\Phi_b(q)^2-1=(\Phi_b(q)-1)(\Phi_b(q)+1)$. Since $p\neq2$, it follows that $r$ divides the number
$$\frac{q^{2^{\gamma+2}}-1}{2}\cdot\frac{q^{2^{\gamma+2}}+3}{2}.$$
Obviously, $r$ does not divide the first product term because $e(r,q)=3\cdot2^{\gamma}$. Assume
that $r$  divides the second term and put $a=q^{2^\gamma}$. Then $r$ divides
$(a^4+3,a^3-1)=(a^3-1,a+3)=(a+3,28)$. However, $e(r,q)\geqslant24$, so $r$ cannot divide $28$, a
contradiction. The proposition is proved.
\end{prf}

\begin{prop}\label{p:tpSeqtpL}
Let $t(L)\geqslant23$. Suppose that the characteristic $p$ of the group $L$ satisfies the conditions:

\emph{(i)} $p$ does not divide the order of $K;$

\emph{(ii)} if $S=L_{m}^\varepsilon(u)$, then $p$ does not divide the number~$\varepsilon{u}-1.$

Then $p$ divides the order of $S$ and $t(p,S)=t(p,L)$.
\end{prop}

\begin{prf} Recall that $p\neq2$ by Proposition~\ref{p:solutionfor2}. Assume to the contrary that either $p$ does not divide $|S|$ or
$t(p,S)\neq t(p,L)$. Observe that if $p$ divides $|S|$ then $t(p,S)\geqslant t(p,L)$. Indeed, it follows from Lemma~\ref{l:tpL} that every prime $r$ nonadjacent to $p$ in $GK(L)$ is large with respect to~$L$, so $r$ does not divide~$|\overline{G}/S|\cdot|K|$.

\begin{lemma}\label{l:pdividesOut} The number $|\overline{G}/S|$ is divided by~$p$. Every $p$-element of
$\overline{G}/S$ is conjugated to a field automorphism of~$S$.
\end{lemma}

\begin{prf}
Suppose that $p$ does not divide $|\overline{G}/S|$. Then $p$ divides $|S|$ and, by our assumption,
$t(p,S)>t(p,L)$. Let $\rho$ be a $\{p\}$-coclique of size $t(p,S)$ in $GK(S)$. Put
$\rho'=\rho\setminus\{p\}$ and $I=E(p,S)=\{e(r,u)\mid r\in\rho'\}$. Then for any $i\in I$ every
prime divisor of $k_i(u)$ is large with respect to~$S$ by Lemma~\ref{l:rcoclique}.
Proposition~\ref{p:tSlesstLp2} yields that for every $i\in I$ there is a prime divisor $r_i$ of
$k_i(u)$ being large with respect to~$L$. The set $\{p\}\cup\{r_i\mid i\in I\}$ is a coclique
in~$GK(S)$. Furthermore, since none of primes from the coclique divides $|K|\cdot|\overline{G}/S|$,
then it is also a $\{p\}$-coclique in~$GK(L)$ of size $t(p,S)$, a contradiction.

Since $(p,2|\operatorname{Inndiag}(S)/S|)=1$, we may assume that every $p$-element of $\overline{G}/S$ is conjugated to a field automorphism of~$S$.
\end{prf}

\begin{lemma}\label{l:specialkju} There exists a positive integer $j$ satisfying the conditions:

\emph{(i)} for every prime $r$ from $R_j(u)$ the inequality $\varphi(r,S)>m/2$ holds and $rp\not\in\omega(S);$

\emph{(ii)} for every prime $s$ from $R_j(u)$, being large with respect to~$L$, the number $sp$ belongs to~$\omega(L)$.
\end{lemma}

\begin{prf}
First, we suppose that either $p$ does not divide $|S|$ or $p$ is large with respect to~$S$. Let
$t=t(L)=t(S)$, and $\rho$ be a coclique in $GK(S)$ of size~$t$ with $p\in\rho$ if $p$
divides~$|S|$. Put $I=E(\rho,S)$. It is clear that there exists a subset $Y$ of $I$ consisting of
at least $t-2$ elements and such that for any $j\in Y$ and every $r\in R_j(u)$ we have
$pr\not\in\omega(S)$ and $\varphi(r,S)>m/2$. Since $t(p,L)\leqslant4$, there is a subset $J$ of $Y$
consisting of at least $t-5$ elements and satisfying the condition: for any $j\in J$ every prime
$s$ from $R_j(u)$, being large with respect to $L$ (for any $j\in J$ such primes from $R_j(u)$
exist by Proposition~\ref{p:tSlesstLp2}), is adjacent to~$p$ in~$GK(L)$. The assumption
$t\geqslant23$ implies that $J$ is not empty.

Suppose now that $t(p,S)<t$,  $e(p,u)=a$, and $\rho$ is a $\{p\}$-coclique in $GK(S)$ of size
$t(p,S)$. Put $I=\{e(r,u)\mid r\in\rho\}\setminus\{a\}$. Then Lemma~\ref{l:rcoclique} yields
$\varphi(r,S)>m/2$ for any $i\in I$ and every $r\in R_i(u)$. By our assumption, $t(p,S)>t(p,L)$, so
there is $j\in I$ such that every $s$ from $k_j(u)$, being large with respect to~$L$, is adjacent
to~$p$ in~$GK(L)$.
\end{prf}

\begin{lemma}\label{l:kjpu0lessqsquare}
Suppose that $p$ divides $\beta$, where $u=v^{\beta}$, and $u_0=u^{1/p}$. Then for the number $j$, defined in Lemma~\emph{\ref{l:specialkju}}, the inequality
$k_{jp}(u_0)\leqslant(q^2-1)\log_vu$ holds.
\end{lemma}

\begin{prf}
Since $\overline{G}/S$ contains a field automorphism of order~$p$, the number $u_0$ is an integer. It follows from Lemma~\ref{l:kjpu0} that $k_{jp}(u_0)$ divides $k_j(u)$. Fix some prime divisor $r$ of~$k_{jp}(u_0)$. Suppose that $r$ does not divide $|\overline{G}/S|\cdot|K|$. Then $r$ is large with respect to~$L$, hence $rp\in\omega(L)$. Since $rp\not\in\omega(S)$, the prime $r$ must divide the order of the centralizer $H=C_{S}(f)$ in $S$ of a field automorphism $f$ of order $p$. The group $H$ is a classical group of the same type and same dimension as $S$, but defined over the field of order~$u_0$. Therefore, using  $\varphi(r,S)>m/2$, it is easy to check that $k_{jp}(u_0)$ is coprime to the order of $H$. Indeed, if $S=L_m^\varepsilon(u)$, then $\varphi(r,H)=\nu_{\varepsilon}(jp)=p\nu_{\varepsilon}(j)=p\varphi(r,S)>pm/2>m$, and if $S$ is symplectic or orthogonal, then $\varphi(r,H)=\eta(jp)=p\eta(j)>pm/2>m$.

Thus, every prime divisor $r$ of $k_{jp}(u_0)$ divides $|\overline{G}/S|\cdot|K|$. Furthermore, $p$ does not divide $k_j(u)$. Applying Lemma~\ref{l:estimkju}, we obtain that $(k_{jp}(u_0))_{\{r\}}$ divides $((q^2-1)\log_vu)_{\{r\}}$ for every prime divisor $r$ of $k_{jp}(u_0)$. The lemma is proved.
\end{prf}

Let us complete the proof of Proposition~\ref{p:tpSeqtpL}. Applying the inequality $t(S)\geqslant23$, Lemma~\ref{l:estimphirsbyn}, and arguing as in the beginning of the proof of Lemma~\ref{l:qbyq}, we obtain that $\eta(j)\geqslant11$ for $j$ defined in Lemma~\ref{l:specialkju}. It follows from Lemma~\ref{l:kjubyu7} that $k_{jp}(u_0)>u^4\log_vu$. Lemma~\ref{l:kjpu0lessqsquare} yields $q^2-1\geqslant k_{jp}(u_0)/\log_vu>u^4$. On the other hand, if we apply the equality $t(L)=t(S)$ to the inequality (\ref{eq:ubyq}), based on Lemma~\ref{l:bigk} and derived in the proof of Lemma~\ref{l:qbyq}, we obtain the inequality $u\geqslant q^{5/9}$. Thus, the impossible chain of inequalities $q^2-1>u^4\geqslant(q^{5/9})^4>q^2$ arises and leads us to a contradiction. The proposition is proved.
\end{prf}

\section{Proof: pigeons and holes}\label{s:transfer}

In this section we complete the proof of Theorem~\ref{t:main}. We write $t=t(L)$ and recall that $t(L)=t(S)$
by Proposition \ref{p:tSlesstLp2}.

\begin{lemma}\label{l:transfer} Let $t(L)\geqslant23$.
If $s\in\pi(L)$ is chosen so that $\varphi(s,L)>n/3$, then $t(s,L)=t(s,S)$.
\end{lemma}

\begin{prf} If $t(s,L)=t(L)$, then the desired result follows from Proposition~\ref{p:largeonlyinS}. We
may assume, therefore, that $s$ is small with respect to $L$. Let $\rho$ be an $\{s\}$-coclique of
size $t(s,L)$ in $GK(L)$. Since $s$ is small with respect to $L$, all other numbers in $\rho$ are
large with respect to $L$ by Lemma~\ref{l:rcoclique}, and hence they do not divide
$|\overline{G}/S|\cdot|K|$. As $\varphi(s,L)>n/3$, Proposition~\ref{p:structureGS2} implies that
$s$ does not divide $|\overline{G}/S|\cdot|K|$ either. Thus $\rho$ is a coclique in $GK(S)$, whence
$t(s,S)\geqslant t(s,L)$.

Let $\rho=\{s,s_1,\ldots,s_{k-1}\}$ be a coclique of size $k=t(s,S)$ in $GK(S)$. By
Lemma~\ref{l:rcoclique}, the numbers $s_1,\ldots,s_{k-1}$ are large with respect to $S$. By
Proposition~\ref{p:tSlesstLp2}, for every $i\in\{1,\ldots,k-1\}$, there is a prime divisor $w_i$
of $k_{e(s_i,u)}$ that is large with respect to $L$. Write
$\rho'=\{s,w_1,\ldots,w_{k-1}\}$. Then $\rho'$ is a coclique in $GK(S)$. Since no number
in $\rho'$ divides $|\overline{G}/S|\cdot|K|$, it follows that $\rho'$ is a coclique in $GK(L)$ of
size $k$, whence $t(s,S)\leqslant t(s,L)$. Thus $t(s,L)=t(s,S)$, and the proof is complete.
\end{prf}

Let $r\in\pi(L)\setminus\{p\}$, and suppose that $r$ is small with respect to $L$. Recall that, by
Lemma \ref{l:rcoclique}, the set $E(\rho',L)=\{e(r,q)\mid r\in\rho'\}$, where $\rho$ is an
$\{r\}$-coclique of greatest size and $\rho'=\rho\setminus\{r\}$, is independent of the choice of
$\rho$, and it is denoted by $J(r,L)$.

\begin{lemma}\label{l:correct}
Let $L$ be a classical simple group over a field of order $q$ and characteristic $p$, and suppose that
$\operatorname{prk}(L)=n\geqslant19$. If a prime divisor $r$ of $|L|$ other than $p$
is chosen so that $\varphi(r,L)>n/3$ or $t(r,L)>2t(L)/3$, then $e(r,q)>2$ and $e(r,q)\neq 6$. In
particular, if $r$ is small then $J(r,L)$ depends only on the number $e(r,q)$ and type of $L$.
\end{lemma}

\begin{prf} Let $t=t(L)$.
If $\varphi(r,L)>n/3$ then $\varphi(r,L)\geqslant 7$, and therefore $e(r,q)\geqslant 7$ or
$e(r,q)=5$. If $t(r,L)>2t/3$ then $t(r,L)\geqslant 8$, and hence $e(r,q)>2$ by Lemma
\ref{l:trLeq2&3&4}. Also it is easy to check that for $n\geqslant19$ and $q\neq2$ we have
$t(r_6(q),L)<2t/3$ (recall that $k_6(2)=1$). Now the result follows from Lemmas \ref{l:rcoclique} and
\ref{l:zsigmondy}.
\end{prf}

Define $$M(L)=\{i\mid\varphi(r_i(q),L)>n/3, t(r_i(q),L)<t\}$$ and
$$N(L)=\{i\mid 2t/3<t(r_i(q),L)<t\}.$$
Also define a function $\zeta_L: M(L)\cup N(L)\mapsto \mathbb{N}$ by setting
$$\zeta_L(i)=t(r_i(q),L).$$

By Lemma \ref{l:correct}, the sets  $M(L)$ and $N(L)$ and function $\zeta_L$ are well-defined.
We write $T(L)$ to denote $\zeta_L(M(L)\cap N(L))$. In four further lemmas, we describe
the sets $M(L)$, $N(L)$, and $T(L)$, as well as the function $\zeta_L$, for all types of classical groups $L$.

\begin{lemma}\label{l:tPLL}
Let $L=L_n^\varepsilon(q)$, where $n\geqslant 45$. If $n\not\equiv5 \pmod 6$, then $N(L)=M(L)$.
If $n\equiv5 \pmod 6$, then $N(L)=M(L)\setminus \{\nu_\varepsilon((n+1)/3)\}$. If $i\in M(L)$, then
$\zeta_L(i)=\nu_\varepsilon(i)$. In particular, $\zeta_L$ is injective and $t-j\in T(L)$ for every
$1\leqslant j\leqslant 7$.
\end{lemma}

\begin{prf} Note that the condition $t(r,L)<t$ is equivalent to the inequality
$\varphi(r,L)<n/2$. Suppose that $i\in M(L)\cup N(L)$, and consider the set $J(i)=J(r_i(q),L)$.
The adjacency criterion yields $J(i)=J_1(i)\setminus J_2(i)$, where $J_1(i)=\{j\mid
n-\nu_\varepsilon(i)<\nu_\varepsilon(j)\leqslant n \}$ and $J_2(i)=\{j\mid \nu_\varepsilon(j) \text{
is a multiple of } \nu_\varepsilon(i)\}$. Since $\nu_\varepsilon$ is bijective, we have $|J_1(i)|=
\nu_\varepsilon(i)$. Thus
$$\zeta_L(i)=1+|J(i)|=1+|J_1(i)|-|J_1(i)\cap
J_2(i)|=1+\nu_\varepsilon(i)-|J_1(i)\cap J_2(i)|.$$

Let $i\in M(L)$. Then $n/3<\nu_\varepsilon(i)<n/2$. Therefore $J_1(i)\cap J_2(i)=\{j\}$, where
$\nu_\varepsilon(j)=2\nu_\varepsilon(i)$, and hence $\zeta_L(i)=\nu_\varepsilon(i)$.
If $n$ is even, then $2t/3=n/3$, and so $\zeta_L(i)>2t/3$. If $n$ is odd and $i\neq
(n+1)/3$, then $\zeta_L(i)\geqslant (n+2)/3>2t/3$. If $n$ is odd and $i=(n+1)/3$, then
$\zeta_L(i)=2t/3$ and $i\not\in N(L)$. Thus $M(L)\subseteq N(L)$ if $n\not\equiv5 \pmod
6$, and $M(L)\setminus N(L)=\{(n+1)/3\}$ otherwise.

Let $i\in N(L)$. Suppose that $\nu_\varepsilon(i)\leqslant n/3$. If
$\nu_\varepsilon(i)\leqslant (n-3)/3$, then $\zeta_L(i)\leqslant 1+\nu_\varepsilon(i)\leqslant
n/3\leqslant 2t/3$, which is a contradiction. If $\nu_\varepsilon(i)\geqslant (n-2)/3$, then
$3\nu_\varepsilon(i)\in J_1(i)\cap J_2(i)$, therefore, $|J(i)|\leqslant \nu_\varepsilon(i)\leqslant
n/3\leqslant 2t/3$, and again we have a contradiction. Thus $i\in M(L)$, and hence $N(L)\subseteq M(L)$.

We conclude that $T(L)=\{x\in\mathbb{Z}\mid n/3<x<n/2$\} if $n$ is even, and $T(L)$ contains
$\{x\in\mathbb{Z}\mid (n+1)/3<x<n/2\}$ if $n$ is odd. It follows easily that $t-1<n/2$. If
$n$ is even, the inequality $t-7>n/3$ is equivalent to $n>42$. If $n$ is odd, the inequality
$t-7>(n+1)/3$ is equivalent to $n>43$. Thus $t-7\in T(L)$ when $n\geqslant 44$.
\end{prf}

\begin{lemma}\label{l:tPLS}
Let $L\in \{S_{2n}(q), O_{2n+1}(q)\}$,  where $n\geqslant 29$. Then $M(L)=N(L)$. If $i\in M(L)$,
then
$$\zeta_L(i)=\begin{cases}
[\frac{3\eta(i)+2}{2}]\text{ if }n\text{ is even,}\\
[\frac{3\eta(i)+3}{2}]\text{ if }n\text{ is odd.}
\end{cases}$$ Moreover,
$$T(L)\cap\{x\mid t-6\leqslant x\leqslant t-1\}=\begin{cases}
\{t-2, t-3, t-5, t-6\}\text{ if }n\equiv 0,3\pmod{4},\\
\{t-1, t-3, t-4, t-6\}\text{ if  }n\equiv 2\pmod{4},\\
\{t-1, t-2, t-4, t-5\}\text{ if }n\equiv 1\pmod{4}.\\
\end{cases}$$
\end{lemma}

\begin{prf} Using Table~\ref{tab:tL}, it is easy to see that
$[(n+1)/2]\leqslant 2t/3\leqslant [(n+1)/2]+2/3.$

Suppose that $i\in M(L)\cup N(L)$, and consider the set $J(i)=J(r_i(q),L)$. The adjacency criterion
yields $J(i)= J_1(i)\setminus J_2(i)$, where $J_1(i)=\{j\mid n-\eta(i)<\eta(j)\leqslant n\}$ and
$J_2(i)=\{j\mid j=li \text{ for an odd $l$}\}$. Therefore $\zeta_L(i)=1+|J_1(i)|-|J_1(i)\cap
J_2(i)|$. By Lemma \ref{l:eta}, the size of $J_1(i)$ is equal to $[3\eta(i)/2]$ if $n$ is even, and
$[(3\eta(i)+1)/2]$ if $n$ is odd.

Let $i\in M(L)$. Then $\eta(i)\geqslant (n+1)/3$. If $j=li$ with $l\geqslant 3$ odd, then
$\eta(j)=l\eta(i)>ln/3\geqslant n$. Therefore $J_1(i)\cap J_2(i)=\varnothing$, and so
$\zeta_L(i)=1+[3\eta(i)/2]=[(3\eta(i)+2)/2]$ when $n$ is even, and
$\zeta_L(i)=1+[(3\eta(i)+1)/2]=[(3\eta(i)+3)/2]$ when $n$ is odd.  Since
$\eta(i)\geqslant (n+1)/3$, we have $3\eta(i)/2\geqslant (n+1)/2$, and so $\zeta_L(i)\geqslant
1+[(n+1)/2]>2t/3$. Thus $M(L)\subseteq N(L)$.

Let $i\in N(L)$. By the above formula, $\zeta_L(i)\leqslant 1+|J_1(i)|\leqslant [(3\eta(i)+3)/2]$. On the
other hand, $\zeta_L(i)>[(n+1)/2]$. Therefore $3\eta(i)+3>n+1$, whence $\eta(i)>(n-2)/3$.
Suppose that $i\not\in M(L)$, that is $\eta(i)\leqslant n/3$. Then $\eta(3i)=3\eta(i)\leqslant
n$ and $\eta(3i)+\eta(i)=4\eta(i)> 4(n-2)/3\geqslant n$. It follows that $3i\in J_1(i)\cap J_2(i)$ and
$\zeta_L(i)\leqslant [(3\eta(i)+1)/2]\leqslant [(n+1)/2]\leqslant 2t/3$, a contradiction. Thus
$N(L)\subseteq M(L)$.

Let $a$ be the largest element of $\eta(M(L))$. Then $a\geqslant (n-2)/2$ when $n\not\equiv
3\pmod 4$, and $a=(n-3)/2$ when $n\equiv 3\pmod 4$. The condition $n\geqslant 29$ yields
$(n-3)/2-3>n/3$, and hence $a,a-1,a-2,a-3$ are in $\eta(M(L))$. Let $\eta(i)=a$ and
$\eta(j)=a-1$. Write $b$ and $c$ for $\zeta_L(i)$ and $\zeta_L(j)$ respectively. As
$[3(x-2)/2]=[3x/2]-3$, we obtain that the four largest elements of $T(L)$ are
$b$, $c$, $b-3$, and $c-3$, and $T(L)\cap\{x\mid t-6\leqslant x\leqslant t-1\}=\{b,c,b-3,c-3\}\cap\{x\mid
t-6\leqslant x\leqslant t-1\}$.

Suppose that $n\equiv 0\pmod 4$. Then $t=(3n+4)/4$ and $a=n/2-1$. If  $\eta(i)=n/2-1$, then
$\zeta_L(i)=[(3n-2)/4]=(3n-4)/4=t-2$. If $\eta(i)=n/2-2$, then $\zeta_L(i)=[(3n-8)/4]=t-3$.

Suppose that $n\equiv 2\pmod 4$. Then $t=(3n+2)/4$ and $a=n/2-1$. If $\eta(i)=n/2-1$, then
$\zeta_L(i)=[(3n-2)/4]=t-1$. If $\eta(i)=n/2-2$, then $\zeta_L(i)=[(3n-8)/4]=t-3$.

Suppose that $n\equiv 1\pmod 4$. Then $t=(3n+5)/4$ and $a=(n-1)/2$. If $\eta(i)=(n-1)/2$, then
$\zeta_L(i)=[(3n+3)/4]=t-1$. If $\eta(i)=(n-3)/2$, then $\zeta_L(i)=[(3n-3)/4]=t-2$.

Suppose that $n\equiv 3\pmod 4$. Then $t=(3n+3)/4$ and $a=(n-3)/2$. If $\eta(i)=(n-3)/2$, then
$\zeta_L(i)=[(3n-3)/4]=t-2$. If $\eta(i)=(n-5)/2$, then $\zeta_L(i)=[(3n-9)/4]=t-3$.
\end{prf}

\begin{lemma}\label{l:tPLOp}
Let $L=O_{2n}^+(q)$, where $n\geqslant 30$. Then
$$N(L)=\begin{cases}
M(L) \text{ if }n\not\equiv 6,9\pmod{12}, \\
M(L)\cup\{2n/3\}  \text{ if }n\equiv 6\pmod{12},\\
M(L)\cup\{2n/3,n/3\} \text{ if }n\equiv 9\pmod{12}.\\
\end{cases}$$

If $i\in M(L)$, then
$$\zeta_L(i)=\begin{cases}
[\frac{3\eta(i)+1}{2}]\text{ if }n\text{ is even,}\\
[\frac{3\eta(i)+2}{2}]\text{ if }n\text{ is odd and $i$ is even},\\
[\frac{3\eta(i)+3}{2}]\text{ if }n\text{ is odd and $i$ is odd}.
\end{cases}$$

If $n\equiv6,9\pmod{12}$ and $\eta(i)=n/3$, then $\zeta_L(i)=[(n+1)/2]=(2t+1)/3$.

If $n$ is odd, then $\zeta_L$ is injective on $M(L)\cap N(L)$, $t-j\in T(L)$ for every
$1\leqslant j\leqslant 6$, and $\zeta_L^{-1}(t-j)\equiv 4\pmod8$ for some $1\leqslant j\leqslant 6$.

If $n$ is even, then
$$T(L)\cap\{x\mid t-6\leqslant x\leqslant t-1\}=\begin{cases}
\{t-1, t-3, t-4, t-6\}\text{ if }n\equiv 0\pmod{4},\\
\{t-1, t-2, t-4, t-5\}\text{ if }n\equiv 2\pmod{4}.\\
\end{cases}$$
\end{lemma}

\begin{prf}
It is easy to verify that $[(n+1)/2]-1/3\leqslant 2t/3\leqslant [(n+1)/2]$.

Suppose that $i\in M(L)\cup N(L)$, and consider the set $J(i)=J(r_i(q),L)$. The adjacency
criterion yields $$J(i)= J_3(i)\cup J_1(i)\setminus J_2(i),$$ where $J_1(i)=\{j\mid
n-\eta(i)<\eta(j)\leqslant n, j\neq 2n\}$, $J_2(i)=\{j\mid j=li \text{ for an odd $l$}\}$, and
$J_3(i)=\{2(n-\eta(i))\}$ when $i$ is odd, $J_3(i)=\{n-\eta(i)\}$ when $i$ is even and
$n-\eta(i)$ is odd, and $J_3(i)=\varnothing$ when both $i$ and $n-\eta(i)$ are even. Note that $J_1(i)\cap
J_3(i)=\varnothing$. Therefore $\zeta_L(i)=1+|J_3(i)|+|J_1(i)|-|J_1(i)\cap J_2(i)|$. By Lemma
\ref{l:eta}, the size of $J_1(i)$ is equal to $[3\eta(i)/2]-1$ or $[(3\eta(i)+1)/2]-1$, depending on whether
$n$ is even or odd.

Thus if $n$ is even and $\eta(i)$ is odd, then
$$\zeta_L(i)=1+[3\eta(i)/2]-|J_1(i)\cap
J_2(i)|=(3\eta(i)+1)/2-|J_1(i)\cap J_2(i)|.$$ If both $n$ and $\eta(i)$ are even, then
$$\zeta_L(i)=[3\eta(i)/2]-|J_1(i)\cap
J_2(i)|=3\eta(i)/2-|J_1(i)\cap J_2(i)|.$$ If both $n$ and $i$ are odd, then
$$\zeta_L(i)=1+(3\eta(i)+1)/2-|J_1(i)\cap
J_2(i)|=(3\eta(i)+3)/2-|J_1(i)\cap J_2(i)|.$$ If $n$ is odd and $\eta(i)$ is even, then
$$\zeta_L(i)=1+3\eta(i)/2-|J_1(i)\cap
J_2(i)|=(3\eta(i)+2)/2-|J_1(i)\cap J_2(i)|.$$ Finally, if $n$ is odd, $i$ is even and $\eta(i)$
is odd, then $$\zeta_L(i)=(3\eta(i)+1)/2-|J_1(i)\cap
J_2(i)|.$$

Let $i\in M(L)$. Then $J_1(i)\cap J_2(i)=\varnothing$, and hence the formulas for $\zeta_L$ are proved.
The condition $\eta(i)>n/3$ implies that $(3\eta(i)+1)/2>(n+1)/2$, therefore,
$$\zeta_L(i)\geqslant [(3\eta(i)+1)/2]\geqslant [(n+1)/2]\geqslant
2t/3.$$ Suppose that all non-strict inequalities in the above chain are equalities. Then
$(n+1)/2$ is integer and $3\eta(i)=n+1$, and so $n$ is odd and $\eta(i)$ is even. But then
the first inequality is strict, a contradiction. Thus $M(L)\subseteq N(L)$.

Let $i\in N(L)$. Then $[(n+1)/2]-1/3\leqslant 2t/3<\zeta_L(i)\leqslant [(3\eta(i)+3)/2]$, and hence
$n\leqslant 3\eta(i)+3$. Suppose that $n-3\leqslant 3\eta(i)\leqslant n-1$. Then
$3i\in J_1(i)\cap J_2(i)$. If $3\eta(i)=n-3$ or $3\eta(i)=n-1$, then $n$ and $\eta(i)$ have different
parity, whence $\zeta_L(i)\leqslant [3\eta(i)/2]\leqslant [(n-1)/2]<2t/3$, a contradiction. If
$3\eta(i)=n-2$, then  $n$ and $\eta(i)$ have the same parity, and so $\zeta_L(i)\leqslant
[(3\eta(i)+1)/2]=[(n-1)/2]<2t/3$, which is not the case. Finally, assume that $\eta(i)=n/3$. If $n$ is odd,
then $J_1(i)\cap J_2(i)=\{3i\}$. If $n$ is even, then $J_1(i)\cap J_2(i)=\varnothing$. Furthermore,
$n-\eta(i)$ is even. It follows that $\zeta_L(i)=[(3\eta(i)+1)/2]=[(n+1)/2]$. If $n\equiv 0,3 \pmod 4$, then
$[(n+1)/2]=2t/3$, a contradiction.  If $n\equiv 1,2 \pmod 4$, then $[(n+1)/2]>2t/3$. Thus
if $n\equiv 6,9 \pmod {12}$ and $\eta(i)=n/3$, then $i\in N(L)\setminus M(L)$.

Let $a$ be the largest element of $\eta(M(L))$. Then $a\geqslant [n/2]-1$. Since
$n\geqslant 30$, we have $a-3>n/3$.

Let $n$ be even, and $b=\zeta_L(i)$, with $\eta(i)=a$, and $c=\zeta_L(j)$, with
$\eta(j)=a-1$. If $n\equiv 0\pmod 4$, then $t=3n/4$, $b=[(3n-4)/4]=t-1$, and $c=[(3n-10)/4]=t-3$. If
$n\equiv 2\pmod 4$, then $t=(3n-2)/4$, $b=[(3n-4)/4]=t-1$, and $c=[(3n-10)/4]=t-2$. Reasoning as
in the proof of Lemma \ref{l:tPLS}, we obtain  that $T(L)\cap\{x\mid t-6\leqslant
x\leqslant t-1\}=\{b,c,b-3,c-3\}$.

Let $n$ be odd. We will show that $\zeta_L$ is injective. Suppose that $i$ and $j$ have the same parity and
$[3\eta(i)/2]=[3\eta(j)/2]$. Then $3\eta(j)-1\leqslant 3\eta(i)\leqslant 3\eta(j)+1$, whence
$\eta(i)=\eta(j)$. By parity condition, it follows that $i=j$. Now suppose that $i$ is even and $j$ is odd,
and $[3\eta(i)/2]=[3\eta(j)/2]+1$. Then $[3i/4]=(3j+1)/2$, and so either $3i/2=3j+1$ or $3i/2=3j+2$.
Both equalities are clearly impossible.

If $n\equiv 1\pmod 4$, then $t=(3n+1)/4$, $a=(n-3)/2$ is odd, and $a,2a\in M(L)$. Therefore
$T(L)$ contains the following numbers: $\zeta_L(a)=(3a+3)/2=3(n-1)/4=t-1$, $\zeta_L(2a)=t-2$,
$\zeta_L(2(a-1))=(3a+1)/2=t-3$, $\zeta_L(a-2)=t-4$, $\zeta_L(2(a-2))=t-5$, and $\zeta_L(2(a-3))=t-6$.
If $n\equiv 3\pmod 4$, then $t=(3n+3)/4$ and $a=(n-1)/2$ is odd with $a\not\in M(L)$. It follows that
$T(L)$ contains the following numbers: $\zeta_L(2a)=(3a+1)/2=(3n-1)/4=t-1$, $\zeta_L(2(a-1))=(3(a-1)+2)/2=t-2$,
$\zeta_L(a-2)=t-3$, $\zeta_L(2(a-2))=t-4$, and $\zeta_L(2(a-3))=t-5$. Furthermore, $a-4>n/3$, and hence
$\zeta_L(a-4)=t-6$ and $\zeta_L(2(a-4))=t-7$ also lie in  $T(L)$. Thus the set $X=\{x\mid t-6\leqslant
x\leqslant t-1\}$ is a subset of $T(L)$. Moreover, one of the numbers $2(a-1)$ and $2(a-3)$ is congruent to 4 modulo
8 and $\zeta_L(2(a-1)), \zeta_L(2(a-3))\in X$.
\end{prf}

\begin{lemma}\label{l:tPLOm}
Let $L=O_{2n}^-(q)$, where $n\geqslant 30$. Then
$$N(L)=\begin{cases}
M(L) \text{ if }n\not\equiv 0,6,9\pmod{12}, \\
M(L)\cup\{2n/3\}  \text{ if }n\equiv 0\pmod{6},\\
M(L)\cup\{2n/3,n/3\} \text{ if }n\equiv 9\pmod{12}.\\
\end{cases}$$

If $i\in M(L)$, then
$$\zeta_L(i)=\begin{cases}
[\frac{3\eta(i)+4}{2}]\text{ if }n\text{ is odd,}\\
[\frac{3\eta(i)+2}{2}]\text{ if both }n\text{ and $i$ are odd},\\
[\frac{3\eta(i)+3}{2}]\text{ if }n\text{ is odd and $i$ is even}.
\end{cases}$$

If $\eta(i)=n/3$, then $$\zeta_L(i)=\begin{cases}
(n+2)/2=(2t+1)/3\text{ if }n\equiv0\pmod{12}\\
(n+2)/2=(2t+2)/3\text{ if }n\equiv6\pmod{12}\\
(n+1)/2=(2t+1)/3\text{ if }n\equiv9\pmod{12}\\
\end{cases}.$$

If $n$ is odd, then $\zeta_L$ is injective on $M(L)\cap N(L)$, $t-j\in T(L)$ for every
$1\leqslant j\leqslant 6$, and $\zeta_L^{-1}(t-j)\equiv 4\pmod8$ for some $1\leqslant j\leqslant 6$.

If $n\equiv 0\pmod{4}$, then
$$T(L)\cap\{x\mid t-7\leqslant x\leqslant t-1\}=
\{t-1, t-2, t-4, t-5, t-7\}.$$

If $n\equiv 2\pmod{4}$ and $n\geqslant 38$, then
$$T(L)\cap\{x\mid t-8\leqslant x\leqslant t-1\}=
\{t-2, t-3, t-5, t-6, t-8\}.$$

If $n=34$, then $T(L)=\{t-2, t-3, t-5, t-6\}$, and if $n=30$, then $T(L)=\{t-2, t-3, t-5\}$.
\end{lemma}

\begin{prf}
It is easy to verify that $n/2+1/6\leqslant 2t/3\leqslant n/2+2/3$.

Suppose that $i\in M(L)\cup N(L)$, and consider the set $J(i)=J(r_i(q),L)$. The adjacency criterion
yields
$$J(i)= J_3(i)\cup J_1(i)\setminus J_2(i),$$ where $J_1(i)=\{j\mid n-\eta(i)<\eta(j)\leqslant n, j\neq n\}$,
$J_2(i)=\{j\mid j=li \text{ for an odd $l$}\}$, and $J_3(i)=\{2(n-\eta(i))\}$ if $i$ is even,
$J_3(i)=\{n-\eta(i)\}$ if both $i$ and $n-\eta(i)$ are odd, and $J_3(i)=\varnothing$ if $i$
is odd, $n-\eta(i)$ is even. Note that $J_1(i)\cap J_3(i)=\varnothing$. Therefore
$\zeta_L(i)=1+|J_3(i)|+|J_1(i)|-|J_1(i)\cap J_2(i)|$. If $n$ is even, then $\eta(n)=n/2\leqslant
n-\eta(i)$, and so the size of $J_1(i)$ is equal to $[3\eta(i)/2]$. If $n$ is odd, then the size of
$J_1(i)$ is equal to $[(3\eta(i)+1)/2]-1$.

Thus if $n$ is even and $\eta(i)$ is odd, then
$$\zeta_L(i)=1+1+(3\eta(i)-1)/2-|J_1(i)\cap J_2(i)|=(3\eta(i)+3)/2-|J_1(i)\cap J_2(i)|.$$
If both $n$ and $\eta(i)$ are even, then
$$\zeta_L(i)=1+1+3\eta(i)/2-|J_1(i)\cap J_2(i)|=(3\eta(i)+4)/2-|J_1(i)\cap J_2(i)|.$$
If  $n$ is odd and $\eta(i)$ is even, then
$$\zeta_L(i)=1+3\eta(i)/2-|J_1(i)\cap J_2(i)|=(3\eta(i)+2)/2-|J_1(i)\cap J_2(i)|.$$
If both $n$ and $i$ are odd, then
$$\zeta_L(i)=(3\eta(i)+1)/2-|J_1(i)\cap J_2(i)|.$$
Finally, if $n$ is odd, $i$ is even, and $\eta(i)$ is odd, then
$$\zeta_L(i)=1+(3\eta(i)+1)/2-|J_1(i)\cap J_2(i)|=(3\eta(i)+3)/2-|J_1(i)\cap J_2(i)|.$$

Let $i\in M(L)$. Then $J_1(i)\cap J_2(i)=\varnothing$, and so the formulas for $\zeta_L$ are proved.
Since $\eta(i)>n/3$, it follows that $(3\eta(i)+3)/2>(n+3)/2$, whence
$$\zeta_L(i)\geqslant [(3\eta(i)+3)/2]\geqslant [(n+3)/2]> 2t/3.$$
Thus $M(L)\subseteq N(L)$.

Let $i\in N(L)$. Then $n/2+1/6\leqslant 2t/3<\zeta_L(i)\leqslant (3\eta(i)+4)/2$, and hence
$n\leqslant 3\eta(i)+3$. Suppose that $n-3\leqslant 3\eta(i)\leqslant n-1$. Then $3i\in
J_1(i)\cap J_2(i)$. If $n$ and $\eta(i)$ have the same parity, then $\zeta_L(i)\leqslant
(3\eta(i)+2)/2\leqslant n/2<2t/3$, a contradiction. If $n$ and $\eta(i)$ have different parity, then
$\zeta_L(i)\leqslant (3\eta(i)+1)/2\leqslant n/2<2t/3$, which is not the case. Suppose that
$\eta(i)=n/3$. If $n$ is even, then $J_1(i)\cap J_2(i)=\{2n\}$, and so $\zeta_L(i)=
(3\eta(i)+2)/2=(n+2)/2>2t/3$. It follows that $2n/3\in N(L)\setminus M(L)$ when $n\equiv 0\pmod
6$. Assume that  $n$ is odd. Then $J_1(i)\cap J_2(i)=\{2n\}$ when $i=2n/3$ and $J_1(i)\cap
J_2(i)=\varnothing$ when $i=n/3$. Therefore $\zeta_L(i)=(3\eta(i)+1)/2=(n+1)/2$. If $n\equiv 3\pmod
4$, then $(n+1)/2=2t/3$, a contradiction. If $n\equiv 1\pmod 4$, then $(n+1)/2>2t/3$. Thus
if $n\equiv 9 \pmod {12}$, then $n/3, 2n/3\in N(L)\setminus M(L)$.

Let $a$ be the largest element of $\eta(M(L))$. Then $a\geqslant [n/2]-1$ when $n\not\equiv 2\pmod
4$ and $a=n/2-2$ when $n\equiv 2\pmod 4$. Therefore if $n>30$, then $a-3>n/3$, and if, in addition, $n$ is even
and $n\neq 34$, then $a-4>n/3$.

If $n\equiv 0\pmod 4$, then $t=(3n+4)/4$ and $a=n/2-1$, whence $\zeta_L(a)=(3a+3)/2=3n/4=t-1$ and
$\zeta_L(2(a-1))=(3a+1)/2=t-2$. If $n\equiv 2\pmod 4$, then $t=(3n+2)/4$ and $a=n/2-2$, and hence
$\zeta_L(a)=(3a+3)/2=(3n-6)/4=t-2$ and $\zeta_L(2(a-1))=(3a+1)/2=t-3$. The completion of the proof
is similar to that of Lemma \ref{l:tPLOp}.
\end{prf}

\begin{lemma}\label{l:TLinTS} If either $S=O_{2m}^+(u)$ with $m\equiv6,9\pmod{12}$ or $S=O_{2m}^-(u)$ with
$m\equiv0,9\pmod{12}$, then $T(L)\subseteq T(S)\cup\{(2t+1)/3\}$. If $S=O_{2m}^-(u)$ with
$m\equiv6\pmod{12}$, then $T(L)\subseteq T(S)\cup\{(2t+2)/3\}$. In other cases, $T(L)\subseteq
T(S)$. In particular, if $x\in T(L)$ and $x>(2t+2)/3$, then $x\in T(S)$.
\end{lemma}

\begin{prf} Choose $x\in T(L)\setminus T(S)$. Let $s\in\pi(L)$ and
$t(s,L)=x$. As $e(s,q)\in M(L)\cap N(L)$, by Lemma~\ref{l:transfer}, we have $t(s,S)=x$, and hence
$e(s,u)\in N(S)$. If $e(s,u)\in M(S)$, then $x\in T(S)$. If
$e(s,u)\not\in M(S)$, then applying  Lemmas~\ref{l:tPLOp} and~\ref{l:tPLOm}, we see that
the assertion is true in this case too.
\end{prf}

Let $\mathcal{L}$ be the class of all classical groups $L$ with $t(L)\geqslant23$. In what follows,
$L$ is always a group lying in this class. Note that the equality $t(L)=t(S)$ implies that
the nonabelian composition factor $S$ of $G$ also lies in~$\mathcal{L}$. We divide
$\mathcal{L}$ on several subclasses, according to the behavior of the function $\zeta_L$. We write $\mathcal{X}=\mathcal{X}_1\cup\mathcal{X}_2$, where
$$\mathcal{X}_1=\{L\in\mathcal{X}\mid L=L^\varepsilon_n(q)\};$$
$$\mathcal{X}_2=\{L\in\mathcal{X}\mid L=O_{2n}^\varepsilon(q), n\equiv1\pmod2\}.$$
Note that $\mathcal{X}$ joins all groups $L$ from $\mathcal{L}$ such that $\zeta_L$ is injective.

Furthermore, put $\mathcal{Y}=\mathcal{L}\setminus\mathcal{X}=\mathcal{Y}_1\cup\mathcal{Y}_2
\cup\mathcal{Y}_3,$ where
\begin{multline*}\mathcal{Y}_1=\{L\in\mathcal{Y}\mid L=S_{2n}(q),O_{2n+1}(q)\text{ and }n\equiv0,3\pmod4,\text{ or } \\ L=O_{2n}^-(q)\text{ and }n\equiv2\pmod4\};\end{multline*}
\begin{multline*}\mathcal{Y}_2=\{L\in\mathcal{Y}\mid L=S_{2n}(q),O_{2n+1}(q)\text{ and }n\equiv2\pmod4,\text{ or } \\ L=O_{2n}^+(q)\text{ and }n\equiv0\pmod4\};\end{multline*}
\begin{multline*}\mathcal{Y}_3=\{L\in\mathcal{Y}\mid L=S_{2n}(q),O_{2n+1}(q)\text{ and }n\equiv1\pmod4,\text{ or } \\ L=O_{2n}^+(q)\text{ and }n\equiv2\pmod{4},
\text{ or }L=O_{2n}^-(q)\text{ and }n\equiv0\pmod4\}.\end{multline*}
Observe that $L$ belongs to $\mathcal{Y}_i$ if and only if $t(L)-i\not\in T(L)$.

It is easy to see that we have the partitions:
$$\mathcal{L}=\mathcal{X}\sqcup\mathcal{Y}=\mathcal{X}_1\sqcup\mathcal{X}_2\sqcup\mathcal{Y}_1\sqcup\mathcal{Y}_2
\sqcup\mathcal{Y}_3.$$

\begin{lemma}\label{l:tmaless2t3} If $t\geqslant23$, then for every natural $a\leqslant6$, we have
$t-a>(2t+2)/3$.
\end{lemma}

\begin{prf} The inequality $t-a<(2t+2)/3$ is equivalent to $a<(t-2)/3$.
\end{prf}

\begin{lemma}\label{l:SinYi} If $i=1,2,3$ and $S\in\mathcal{Y}_i$, then $L\in\mathcal{Y}_i$. In particular,
if $L\in\mathcal{X}$, then $S\in\mathcal{X}$.
\end{lemma}

\begin{prf}
Assume that the conclusion is false. Choose $i\in\{1,2,3\}$ and write $x=t-i$. By assumption,
$L\not\in\mathcal{Y}_i$, therefore, $x\in T(L)$ by Lemmas~\ref{l:tPLL}--\ref{l:tPLOm}. On the other hand,
$x=t-i>(2t+2)/3$ by Lemma~\ref{l:tmaless2t3}, and so Lemmas~\ref{l:tPLOp} and~\ref{l:tPLOm}
yield $x\not\in T(S)\cup\{(2t+1)/3,(2t+2)/3\}$. This contradicts Lemma~\ref{l:TLinTS}.
\end{prf}

Our next goal is to show that $p$ does not divide $|K|$ and, in the case $S=L_m^\varepsilon(u)$,
that it does not divide $\varepsilon{u}-1$ either. For this purpose, we
need the following result.

\begin{lemma}\label{l:rsnotprs} Let $j\in N(L)\cap M(L)$. There is a natural number $i$
such that $r_i(q)$ is large with respect to $L$,  $\varphi(r_i(q),L)<2n/3$, $r_i(q)r_j(q)$ lies in
$\omega(L)$, and $pr_i(q)r_j(q)$ does not. If, in addition, $L=L_n^\varepsilon(q)$ and
$\nu_\varepsilon(j)<(n-1)/2$, or $L=S_{2n}(q), O_{2n+1}(q)$ and $n-\eta(j)$ is odd, or
$L=O_{2n}^\varepsilon(q)$ and $\eta(j)<(n-1)/2$, then there are two distinct $i$ and $i'$
satisfying these conditions, and such that  $r_i(q)r_{i'}(q)\not\in\omega(L)$.
\end{lemma}

\begin{prf}

Our proof repeatedly uses Lemma~\ref{l:criterion}, and we apply this lemma without further mention.

If $L=L_n^\varepsilon(q)$, then $i=\nu_\varepsilon(n-\nu_\varepsilon(j))$ is the desired number. If,
in addition, $\nu_\varepsilon(j)<(n-1)/2$, then $i'=\nu_\varepsilon(n-1-\nu_\varepsilon(j))$ also has
the desired properties.

Let $L=S_{2n}(q), O_{2n+1}(q)$ and write $a=n-\eta(j)$. Then $i=2a$ is the desired number.
If $a$ is odd, then $i'=a$ also has the desired properties, since $\eta(i')=\eta(i)=a$.

Let $L=O_{2n}^\varepsilon(q)$ and $\eta(j)<n/2-1$. Writing $a=n-1-\eta(j)$, we see that $i=2a$ has
the desired properties, and if $a$ is odd, then so does $i'=a$. If $a$ is even, then
$a+1=n-\eta(j)$ is odd. If $j=\eta(j)$ and $L=O_{2n}^+$, or $j=2\eta(j)$ and $L=O_{2n}^-$, then
$i'=a+1$ is the desired number. And if $j=2\eta(j)$ and $L=O_{2n}^+$, or $j=\eta(j)$ and
$L=O_{2n}^-$, then so is $i'=2(a+1)$. Now suppose that $\eta(j)=n/2-1$. If $n\equiv0\pmod4$, then
$a=n/2$ is even, and we may choose the required numbers $i$ and $i'$ in a similar manner. If
$n\equiv2\pmod4$, then $L=O_{2n}^+$ (otherwise $r_j(q)$ is large with respect to $L$), and hence
$i=n/2$ and $i'=n+1$ are the desired numbers. Finally, let $\eta(j)=(n-1)/2$. Since $r_j(q)$ is
small, the only possibility is that $n\equiv3\pmod4$ and, moreover, $L=O_{2n}^+(q)$,
$j=2\eta(j)=n-1$ or $L=O_{2n}^-(q)$, $j=\eta(j)=(n-1)/2$. In both cases $i=n+1=2(n-\eta(j))$ has
the desired properties.
\end{prf}

\begin{lemma}\label{l:finishK} The order of $K$ is not divisible by $p$.
\end{lemma}

\begin{prf} We derive a contradiction by assuming that $p$
divides the order of the soluble radical $K$, which is nilpotent by
Proposition~\ref{p:structureK}. Let $P$ be a Sylow $p$-subgroup of $K$, and let $V$
be the factor group $P/\Phi(P)$. The group $S$ acts on $V$ via
conjugation. If this action is not faithful, then $p$ is adjacent to all primes that are large
with respect to $L$, which is not the case. Thus $S$ acts faithfully on $V$.

By Proposition~\ref{p:structureK} and Lemma~\ref{l:tpL}, we may assume that
$L\in\{S_{2n}(q),$  $O_{2n+1}(q)\mid n\equiv0\pmod2\}$. Suppose first that $n\equiv0\pmod4$, and choose
$s\in\pi(L)$ so that $\varphi(s,L)=\eta(e(s,q))=(n-2)/2$. By Lemma~\ref{l:tPLS}, it follows that
$t(s,L)=[(3\eta(e(s,q))+2)/2]=(3n-4)/4=t-2$. If $n\equiv2\pmod4$, then we choose $s$ so that
$\varphi(s,L)=\eta(e(s,q))=(n-4)/2$, and hence $t(s,L)=[(3e(s,q)+2)/2]=(3n-10)/4=t-3$. In both
cases, we have $e(s,q)\in M(L)\cap N(L)$. By Lemma~\ref{l:rsnotprs}, there is $a$ such that
$r_a(q)$ is large with respect to $L$, $\varphi(r_a(q),L)<2n/3$, $r_a(q)s$ lies in $\omega(L)$, and
$pr_a(q)s$ does not, and we write $r=r_a(q)$ and $i=e(r,u)$.

Lemma~\ref{l:transfer} yields $x=t(s,S)=t(s,L)\in\{t-2,t-3\}$, and in particular,  $x\in T(S)$.
Therefore $j=e(s,u)\in M(S)\cap N(S)$ and $m/3<\varphi(s,S)<m/2$.  Since $r$ is large with
respect to $L$, it is also large with respect to $S$.  As $rs\in\omega(L)$ and
$rs$ does not divide $|K|\cdot|\overline{G}/S|$, it follows that $rs\in\omega(S)$.

Suppose first that $\varphi(s,S)+\varphi(r,S)>m$. If $S$ is symplectic or orthogonal, then
adjacency of $r$ and $s$ in $GK(S)$ implies that $i/j$ is an odd integer. Since $r$ is large with
respect to $S$ and $s$ is small, we have $i\neq j$. Therefore $i\geqslant3j$, whence
$\eta(i)\geqslant3\eta(j)>3m/3=m$, which is a contradiction. Let $S=L_m^\varepsilon(u)$. Then
adjacency of $r$ and $s$ in $GK(S)$ together with the inequality $\varphi(s,S)+\varphi(r,S)>m$
yields $\nu_\varepsilon(i)=2\nu_\varepsilon(j)$. In particular, $\nu_\varepsilon(i)$ is even. As
$s$ is small with respect to $S$,  we have $\nu_\varepsilon(j)\neq m/2$. If
$\nu_\varepsilon(j)=(m-1)/2$, then $t(s,S)=t-1\not\in\{t-2,t-3\}$, which contradicts the choice
of~$s$. Thus $\nu_\varepsilon(i)=2\nu_\varepsilon(j)<m-1$. By \cite[Lemma~5]{VasGr} and
\cite[Lemma 5]{GrCovCl}, there is a subgroup of $S$ that is a Frobenius group with kernel being a
$v$-group and complement being a cyclic group of order $|(\varepsilon{u})^{\nu_\varepsilon(i)}-1|/d$, where $\pi(d)\subseteq\pi(u^2-1)$, and in particular, of order divisible by $rs$. Applying
Lemma~\ref{l:action} to this subgroup of $S$ acting on $V=P/\Phi(P)$, we conclude that
$prs\in\omega(G)$, which contradicts the choice of $r$ and $s$. Thus $\varphi(r,S)\leqslant
m-\varphi(s,S)<2m/3$.

Now observe that we chose $s$ so that $n-\eta(e(s,q))$ is odd. By Lemma~\ref{l:rsnotprs},
there are at least two large with respect to $L$ and nonadjacent numbers
$r$ and $w$ such that $sr$ and $sw$ lie in $\omega(L)$, and $psr$ and $psw$ do not.
For definiteness, we assume that $\varphi(r,S)\leqslant\varphi(w,S)$. It follows by the result
of the previous paragraph that $\varphi(r,S)\leqslant\varphi(w,S)\leqslant m-\varphi(s,S)$.

Writing $k=\varphi(s,S)$, we have  $m/3< k<m/2$. By \cite[Propositions~4.1.3, 4.1.4, 4.1.6]{KL},
$S$ contains a central product of subgroups $\overline{A}$ and
$\overline{B}$ having nonabelian composition factors $A$ and $B$, respectively, such that
$A$ and $B$ are also simple classical groups over the same field of order~$u$. The groups
$\overline{A}$ and $\overline{B}$ can be chosen as follows. If $S=L_m^\varepsilon(u)$, then
$A\simeq L_k^\varepsilon(u)$ and $B\simeq L_{m-k}^\varepsilon(u)$. If $S=S_{2m}(u)$, then $A\simeq
S_{2k}(u)$ and $B\simeq S_{2(m-k)}(u)$. If $S=O_{2m+1}(u)$ and $j=k$, then $A\simeq O_{2k}^+(u)$ and
$B\simeq O_{2(m-k)+1}(u)$. If $S=O_{2m+1}(u)$ and $j=2k$, then $A\simeq O_{2k}^-(u)$ and $B\simeq
O_{2(m-k)+1}(u)$. If  $S=O_{2m}^\varepsilon(u)$ and $j=k$, then $A\simeq O_{2k}^+(u)$ and $B\simeq
O_{2(m-k)}^\varepsilon(u)$. If $S=O_{2m}^\varepsilon(u)$ and $j=2k$, then $A\simeq O_{2k}^-(u)$ and $B\simeq
O_{2(m-k)}^{-\varepsilon}(u)$.

By our choice, it follows that $s\in\pi(A)$. We claim that $r,w\in\pi(B)$. If $S\neq
O_{2m}^\varepsilon(u)$, then since $\varphi(r,S)\leqslant\varphi(w,S)\leqslant m-k$, we have
$r,w\in\pi(B)$. Let $S=O_{2m}^+(u)$ and $j=k$. If $\varphi(w,S)<m-k$, then $r,w\in\pi(B)$.
Suppose that $\varphi(w,S)=m-k$. As $sw\in\omega(S)$, the adjacency criterion yields
$e(w,u)=\varphi(w,S)$, and therefore, $w\in\pi(B)$. If $\varphi(r,S)=\varphi(w,S)=m-k$, then by
similar reasoning, $e(r,u)=\varphi(r,S)=e(w,u)$, which is not the case because $rw\not\in\omega(L)$.
Thus $\varphi(r,S)<\varphi(w,S)$ and $r\in\pi(B)$. The other cases can be handled in a similar
manner.

The numbers $s$, $r$, and $w$ are coprime to $u^2-1$, and hence to the orders of the centers of
$\overline{A}$ and $\overline{B}$ either. Therefore, $s$ divides $|A|$, while $r$ and $w$ divide $|B|$
and are large with respect to $B$. As $r$ and $w$ are not adjacent in $GK(L)$, they are not
adjacent in $GK(B)$ either, and by Lemma~\ref{l:adjanisotrop}, at least one of them must divide the
order of some proper parabolic subgroup $F$ of $B$. For definiteness, we assume that this number is
$r$ (there is no loss in making this assumption, since we will not use the condition
$\varphi(r,S)\leqslant\varphi(w,S)$ or the number~$w$ itself). Let $y$ be an element of order $r$
lying in the preimage of $F$ in $\overline{B}$, and let $x$ be an element of order $s$ lying
in~$\overline{A}$. We again consider the action of $S$ on $V$ via conjugation. By the choice of
$x$, it clearly lies in some proper parabolic subgroup of $S$. By Lemmas \ref{l:hallhigman}
and~\ref{l:good}, the degree of the minimal polynomial of $x$ on $V$ is equal to $s$, and
therefore, its centralizer $U$ in $V$ is nontrivial. Since $\overline{B}$ acting on $V$ normalizes
$U$, we may consider its action on $U$. If the kernel $J$ of this action does not lie in
$Z(\overline{B})$, then its image in the factor group $\overline{B}/Z(\overline{B})$ includes $B$,
so, in particular, $y$ centralizes $U$. If $J\leqslant Z(\overline{B})$, then $B$ acts faithfully
on $U$, and again by Lemmas~\ref{l:hallhigman} and~\ref{l:good}, the centralizer $C_U(y)$ is
not trivial. In both cases, there is  $z\in U$ such that $z^y=z$. Therefore, the element $g=zxy$ of
the natural semidirect product $V\rtimes S$ has order $psr$. By Lemma~\ref{l:semidirect},
the group $\widetilde{G}=G/O_{p'}(K)\Phi(P)$ has an element of the same order, and therefore so
does $G$. On the other hand, by our choice of $s$ and $r$, there is no element of order $psr$ in
$L$, a contradiction.
\end{prf}

\begin{lemma}\label{l:pdoesnotdivideum1}
If $S=L_{m}^\varepsilon(u)$, then $p$ does not divide $\varepsilon{u}-1$.
\end{lemma}

\begin{prf}

By Lemmas~\ref{l:tPLL}--\ref{l:tPLOm}, we can choose a prime $s\in\pi(L)$ with $e(s,q)\in M(L)\cap N(L)$ so that
$(2t+2)/3<t(s,L)<t-1$, and if $L\in\{S_{2n}(q),O_{2n+1}(q)\}$, then also $n-\varphi(s,L)$
is odd. By Lemma~\ref{l:rsnotprs}, there are at least two large with respect to $L$ and nonadjacent
numbers $r$ and $w$ such that $sr,sw\in\omega(L)$, but $psr,psw\not\in\omega(L)$.

As in the proof of the previous lemma, if $\varphi(r,S)+\varphi(s,S)>m$, then $\varphi(r,S)=2\varphi(s,S)$. As
$t(s,S)=t(s,L)<t-1$, it follows that $i=\varphi(r,S)<m-1$. By \cite[Corollary~3]{ButLU},  the set $\omega(S)$
contains a number $[(\varepsilon{u})^i-1,\varepsilon{u}-1]$, which is divisible by $prs$,
and this contradicts the choice of $r$ and $s$. Therefore we may assume that
$\varphi(r,S)\leqslant\varphi(w,S)\leqslant m-\varphi(s,S)$. Since $rw\not\in\omega(S)$ and
$S$ is linear or unitary, we have $\varphi(r,S)<\varphi(w,S)$, and hence
$\varphi(r,S)+\varphi(s,S)<m$. Writing $i=\varphi(r,S)$ and $j=\varphi(s,S)$, by
\cite[Corollary~3]{ButLU} we conclude that $S$ has an element of order
$[(\varepsilon{u})^i-1,(\varepsilon{u})^j-1,\varepsilon{u}-1]$, which is divisible by $prs$. This contradiction
completes the proof.
\end{prf}

\begin{prop}\label{p:tpLeqtpS} The number $p$ divides $|S|$, and $l=t(p,S)=t(p,L)\in\{2,3,4\}$. Furthermore,
$k=e(p,u)$ lies in the set $K(l,S)$ defined in Table~\emph{\ref{tab:KlL}}.
\end{prop}

\begin{prf} Lemmas~\ref{l:finishK}, \ref{l:pdoesnotdivideum1} and Proposition~\ref{p:tpSeqtpL} imply that $p$
divides $|S|$ and $l=t(p,S)=t(p,L)$. Lemma~\ref{l:tpL} yields $l\in\{2,3,4\}$. The final assertion follows
by Lemma~\ref{l:trLeq2&3&4} and Table~\ref{tab:trL}.
\end{prf}

  \begin{table}[!th]
  \caption{The set $K(l,S)$}\label{tab:KlL}
  \begin{center}{\small
  \begin{tabular}{|l|l|l|l|}
  \hline
  $S$ & $l=2$ & $l=3$ & $l=4$\\
  \hline
  $L_m(u)$& $\{2\}$ & $\{3\}$ & $\{4\}$\\
  \hline
  $U_m(u)$& $\{1\}$ & $\{6\}$ & $\{4\}$\\
  \hline
  $Sp_{2m}(u)$& $\{1,2\}$&$\{4\}$, if  $m\equiv 2,3\pmod 4$& $\{4\}$, if $m\equiv 0,1,5,8,9\pmod{12}$\\
  $O_{2m+1}(u)$&&&$\{3,6\}$, if $m\equiv 10\pmod{12}$\\
  &&&$\{3,4,6\}$, if $m\equiv 4\pmod{12}$\\
  \hline
  $O_{2m}^+(u)$&$\{1,2\}$&$\{4\}$, if $m\not\equiv 1\pmod 4$&$\{3,6\}$, if $m\equiv 4\pmod 6$\\
  &&& $\{4\}$, if $m\equiv
  1,9\pmod{12}$\\
  &&&$\{6\}$, if $m\equiv 11\pmod{12}$\\
  &&&$\{4,6\}$, if $m\equiv 5\pmod{12}$\\
  \hline
  $O_{2m}^-(u)$&$\{1,2\}$&$\{4\}$, if $m\equiv 3\pmod 4$& $\{4\}$, if $m\not\equiv 3,5,7,11\pmod{12}$\\
  &&&$\{3\}$, if $m\equiv 11\pmod{12}$\\
  &&&$\{3,4\}$, if $m\equiv 5\pmod{12}$\\
  \hline
  \end{tabular}}
  \end{center}
  \end{table}

Let $j$ and $k$ be natural numbers and $j\geqslant k$. Observe that $j/k$ is an odd integer if
and only if  $j\equiv k\pmod{2k}$.

\begin{lemma}\label{l:kdividesj}
Let $k=e(p,u)$. Suppose that $s,r\in\pi(L)\setminus\{p\}$ are different primes such that
$(sr,u(u^2-1)\cdot|K|\cdot|\overline{G}/S|)=1$, $sr\in\omega(L)$, and $psr\not\in\omega(L)$, and suppose that
$j=e(s,u)$ and $i=e(r,u)$. If $S=L_n^\varepsilon(u)$, then
neither $\nu_\varepsilon(j)$ nor $\nu_\varepsilon(i)$ is divisible by $\nu_\varepsilon(k)$. If $S$ is symplectic or
orthogonal, then neither $i$ nor $j$ is congruent to $k$ modulo $2k$. In particular,
if $a\in M(L)\cap N(L)$ and $j=e(r_a(q),u)$, then $\nu_\varepsilon(j)$ is not divisible by
$\nu_\varepsilon(k)$ when $S=L_n^\varepsilon(u)$, and $j$ is not congruent to $k$ modulo $2k$
when $S$ is symplectic or orthogonal.
\end{lemma}

\begin{prf}
First, suppose that $S=L_m^\varepsilon(u)$. To avoid unwieldy notation, we assume that $S$
is a linear group, that is $\nu_\varepsilon$ is the identity function. The proof for unitary groups
is analogous, we have only to use the function $\nu$ instead. Recall that
by Lemma~\ref{l:pdoesnotdivideum1}, the number $p$ does not divide $u-1$, that is we have $k\neq1$. It follows
from the equality $(sr,u(u^2-1)\cdot|K|\cdot|\overline{G}/S|)=1$ that $S$ has a semisimple element of
order $sr$. Now Lemma~\ref{l:criterion} implies that $S$ has an element of order $[u^j-1, u^i-1]/d$, where
$d$ is some divisor of $u-1$. If $k$ divides either of $j$ and $i$, then $prs$ divides
$[u^j-1, u^i-1]/d$, which contradicts the fact that $prs\not\in\omega(L)$.

Now let  $S$ be a symplectic or orthogonal group. Recall that $p\neq2$. The equality
$(sr,u(u^2-1)\cdot|K|\cdot|\overline{G}/S|)=1$ implies that $S$ has a semisimple  element of
order $sr$. By Lemma~\ref{l:criterion},  $S$ has an element of order $[u^{\eta(j)}+(-1)^j,
u^{\eta(i)}+(-1)^i]/d$, where $d$ divides $4$. Assume that one of the numbers $j$ and $i$,
say $j$, is congruent to $k$ modulo $2k$. Then $u^{\eta(k)}+(-1)^k$ divides
$u^{\eta(j)}+(-1)^j$, and hence $prs$ divides $[u^{\eta(j)}+(-1)^j, u^{\eta(i)}+(-1)^i]/d$,
contrary to the hypothesis.

Let $a\in M(L)\cap N(L)$ and let $j=e(r_a(q),u)$. By Lemma~\ref{l:rsnotprs}, there is a number $b$
such that $r_b(q)$ is large with respect to  $L$,  $r_a(q)r_b(q)$ lies in $\omega(L)$, and
$pr_a(q)r_b(q)$ does not. Writing $s=r_a(q)$, $r=r_b(q)$ and observing that
$(sr,u(u^2-1)\cdot|K|\cdot|\overline{G}/S|)=1$, we see that the assertion of the lemma holds for
$j=e(r_a(q),u)=e(s,u)$.
\end{prf}

\begin{lemma}\label{l:LnotinX}
If $L\in\mathcal{X}$, then the conclusion of Theorem~{\em\ref{t:main}} holds.
\end{lemma}

\begin{prf}
Let $L\in\mathcal{X}$. Then Lemma~\ref{l:SinYi} implies that $S$ also lies in $\mathcal{X}$.
Suppose that $S=L^\varepsilon_m(u)$. Since $L\in\mathcal{X}$, we have $t(p,L)=3$, and hence $k=e(p,u)=\nu_\varepsilon(3)$
(see Table~\ref{tab:KlL}). One of the numbers $t-1$, $t-2$, $t-3$ is a multiple of $3$, and we denote this
number by $c$. By Lemma~\ref{l:tmaless2t3}, we have $(2t+2)/3<c<t$, and so $c\in
T(L)$. Therefore, there is $a\in M(L)\cap N(L)$ such that $\zeta_L(a)=c$. Then
$c=t(r_a(q),L)=t(r_a(q),S)$ due to Lemma~\ref{l:transfer}. By Lemma~\ref{l:tPLL}, the function $\zeta_S$ is injective and $\zeta_S(\nu_\varepsilon(c))=c$, and
hence $j=e(r_a(q),u)=\nu_\varepsilon(c)$. Since $\nu_\varepsilon=\nu_\varepsilon^{-1}$,
we see that $\nu_\varepsilon(j)=c$ is divisible by $\nu_\varepsilon(k)=3$, which is a
contradiction by Lemma~\ref{l:kdividesj}.

Let $S$ be an orthogonal group lying in $\mathcal{X}$. Then  $k=e(p,u)=4$. By Lemmas~\ref{l:tPLOp} and~\ref{l:tPLOm}, the set
$C=\{t-x\mid x=1,2,\ldots,6\}$ contains $c$ such that $t(r_j(u),S)=c$ and $j\equiv4\pmod8$.
Since $C\subseteq T(L)$, it follows that $c\in T(L)$. Thus there is $a\in M(L)\cap N(L)$ such that
$\zeta_L(a)=c$. Then $t(r_a(q),S)=c$, and hence $e(r_a(q),u)=j$, which is not the case by Lemma~\ref{l:kdividesj}.
\end{prf}

\begin{lemma}\label{l:SisnotX1} If $S\in\mathcal{X}_1$, then the conclusion of Theorem~{\em\ref{t:main}} holds.
\end{lemma}

\begin{prf}
Let $k=e(p,u)$. Suppose that there is $c\in T(L)$ such that $\nu_\varepsilon(k)$ divides $c$.
As in the previous lemma, there is $a\in M(L)\cap N(L)$ such that $t(r_a(q),S)=c$, and by Lemma~\ref{l:tPLL},
it follows that $j=e(r_{a}(q),u)=\nu_\varepsilon(c)$. Hence $\nu_\varepsilon(j)=c$ is divisible by $\nu_\varepsilon(k)$, which is
a contradiction due to Lemma~\ref{l:kdividesj}. Thus it suffices to show that there is $j\in T(L)$ that
is a multiple of $\nu_\varepsilon(k)$.

By the previous result, we may assume that $L\in\mathcal{Y}$. In this case, an additional
difficulty is that the members of $T(L)$ are not consecutive integers. Nevertheless, for each group
$L$ we will prove that $T(L)$ contains a number with the desired property.

Let $L\in\{S_{2n}(q),O_{2n+1}(q)\}$ and $n$ be even. Then $t(p,L)=2$. Since
$S=L_{m}^\varepsilon(u)$, and therefore, the characteristic $p$
of $L$ does not divide $\varepsilon{u}-1$ by Lemma~\ref{l:pdoesnotdivideum1},
it follows that $k=e(p,u)=\nu_\varepsilon(2)$, and hence $\nu_\varepsilon(k)=2$.
If $n\equiv0\pmod4$, then $t-2,t-3\in T(L)$, and if
$n\equiv2\pmod4$, then $t-3,t-4\in T(L)$. In each case,  $T(L)$ contains two consequent numbers,
and one of them is a multiple of $2$, as required.

Let  $L=O_{2n}^-(q)$ and $n$ be even. Then $t(p,L)=4$, whence $k=e(p,u)=\nu_\varepsilon(4)=4$.
Now the proof is somewhat similar to the proof in the previous paragraph if we take
into account the following observation: $T(L)$ contains the numbers $t-2,t-4,t-5,t-7$
when $n\equiv0\pmod{4}$ and $t-2,t-4,t-5,t-7$ when $n\equiv2\pmod{4}$ and $n>34$,
and the set of the residues of these numbers module 4 is all of $\{0,1,2,3\}$;
furthermore, if $n=30$, then $t-3=23-3=20$ is a multiple of $4$ and lies in
$t(L)$, and if  $n=34$, then so does $t-2=26-2=24$.

The remaining cases are $L=S_{2n}(q)$, $O_{2n+1}(q)$ with $n$ odd and
$L=O_{2n}^+(q)$ with $n$ even, where we have $t(p,L)=3$, whence $k=e(p,u)=\nu_\varepsilon(3)$ and
$\nu_\varepsilon(k)=3$. If $L=S_{2n}(q)$, $O_{2n+1}(q)$ and $n\equiv1\pmod4$, then $t\equiv2\pmod3$.
Thus $t-2$ lies in $T(L)$ and is a multiple of $3$. If either $L=S_{2n}(q)$, $O_{2n+1}(q)$ with $n\equiv3\pmod4$
or $L=O_{2n}^+(q)$ with $n\equiv0\pmod4$, then $t\equiv0\pmod3$, and $t-3$ lies in $T(L)$ and is a multiple of $3$.
If $L=O_{2n}^+(q)$ and $n\equiv2\pmod4$, then $t\equiv1\pmod3$, and  $t-1$ is the desired number.
\end{prf}

Now we may assume that $L\in\mathcal{Y}$ and $S\in\mathcal{Y}\cup\mathcal{X}_2$, and in particular,
both $L$ and $S$ are symplectic or orthogonal groups.  Moreover, by
Lemma~\ref{l:SinYi}, each $i\in\{1,2,3\}$ satisfies the following condition: if
$L\in\mathcal{Y}_i$, then $S\in\mathcal{Y}_i\cup\mathcal{X}_2$. Also observe that $[(3x+2)/2]$ is an injective
function, and hence $\eta(i)$ is uniquely determined by $\zeta_L(i)$. Furthermore, if $i\equiv 0\pmod 4$, then
$i$ is uniquely determined by $\eta(i)$, and so by $\zeta_L(i)$ as well.

\begin{lemma}\label{l:tpLneq3}
If $t(p,L)=3$, then the conclusion of Theorem~{\em\ref{t:main}} holds.
\end{lemma}

\begin{prf}
Let $t(p,L)=3$. Then either $L=S_{2n}(q)$, $O_{2n+1}(q)$ with $n$ odd or $L=O_{2n}^+(q)$  with
$n$ even. Furthermore, since $S$ is symplectic or orthogonal, we have $k=e(p,u)=4$.
Also note that if $S=S_{2m}(u)$, $O_{2m+1}(u)$, then $m\equiv2,3\pmod4$; if $S=O_{2m}^+(u)$, then
$m\not\equiv1\pmod4$; and if $S=O_{2m}^-(u)$, then $m\equiv3\pmod4$ (see Table~\ref{tab:KlL}).

Let $L=S_{2n}(q)$, $O_{2n+1}(q)$ and $n\equiv1\pmod4$. Suppose that $S=S_{2m}(u)$ or
$S=O_{2m+1}(u)$. Since $t(L)=t(S)$, it follows that $m\in\{n,n+1\}$.
Now the congruence $m\equiv2,3\pmod4$ yields $m=n+1$. But then $S\in\mathcal{Y}_2$,
which is not the case because $L\in\mathcal{Y}_3$. If $S=O_{2m}^+(u)$,
then for every $m$, we have that $t(S)$ is not congruent to $2$ modulo $3$. On the
other hand, $t(L)\equiv2\pmod3$,  and hence $S\neq O_{2m}^+(u)$.
Similarly, $S\neq O_{2m}^-(u)$, since otherwise $m\equiv3\pmod4$ and $t(S)\equiv0\pmod3$.

Let $L=S_{2n}(q)$, $O_{2n+1}(q)$ and $n\equiv3\pmod4$. Suppose that $S=S_{2m}(u)$ or
$S=O_{2m+1}(u)$. The equality $t(L)=t(S)$ yields $m=n$. Both $n-3$ and $n-7$ are divisible by $4$,
therefore, one of them is congruent to $4$ modulo $8$. Denote this number by $a$. Since
$t(r_{n-3}(q),L)=t-2$ and $t(r_{n-7}(q),L)=t-5$, we have $x=t(r_a(q),L)>2t/3$, whence $x\in
T(L)$. Thus $t(r_a(q),S)=x$. As $\zeta_S$ is invertible on the set of multiples of $4$ and $S$
has the same type as $L$, it follows that $e(r_a(q),u)$ is equal to $a$, and so it is
congruent to $4$ modulo $8$. Applying Lemma~\ref{l:kdividesj}, we derive a contradiction. Thus
$S$ is an orthogonal group of even dimension. If $m=n$, then $S=O_{2n}^\varepsilon(u)$.
By Lemmas~\ref{l:tPLS}--\ref{l:tPLOm}, both $t(r_{n-3}(q),L)$ and $t(r_{n-3}(u),S)$ are equal to $t-2$, so $t-2=t(r_{n-3}(u),S)=t(r_{n-3}(q),L)=t(r_{n-3}(q),S)$. Similarly, $t-5=t(r_{n-7}(u),S)=t(r_{n-7}(q),L)=t(r_{n-7}(q),S)$. Repeating the previous argument
we derive a contradiction by Lemma~\ref{l:kdividesj}. If $m\neq n$, then the equality $t(L)=t(S)$ implies that the
only remaining possibility is $S=O_{2n+2}^+(u)$.  Then $S\in\mathcal{Y}_2$, which
is not the case because $L\in\mathcal{Y}_1$.

Let $L=O_{2n}^+(q)$ and $n\equiv2\pmod4$. If $S\in\{S_{2m}(u), O_{2m+1}(u)\}$, then
$t(L)=t(S)$ yields $m=n-2$, whence $m\equiv0\pmod4$, which is a contradiction.
If $S=O_{2m}^-(u)$, then $m\in\{n-2,n-1\}$, and this is also a contradiction since  $m$
must be congruent to $3$ modulo~$4$. If $S=O_{2m}^+(u)$, then $m\in\{n-1,n\}$. Since $m$ cannot
be congruent to $1$ modulo~$4$, it follows that $m\neq n-1$. Thus $S$ has the same
type and the same Lie rank $n$ as $L$. Both $n-2$ and $n-6$ are multiples of $4$, so
one of them is congruent to 4 modulo 8, and we denote this number by $a$. As $t(L)\geqslant23$, we have $n\geqslant31$ and
$\eta(n-2)>\eta(n-6)=n/2-3>n/3$. Therefore $a\in M(L)\cap N(L)$, and applying Lemma~\ref{l:kdividesj},
we derive a contradiction.

Let  $L=O_{2n}^+(q)$ and $n\equiv0\pmod4$. The equality $t(L)=t(S)$ leaves us with two possibilities:
either $m=n-1\equiv3\pmod4$ or $m=n$ and $S=O_{2n}^+(u)$. If $S$ is a symplectic group or an orthogonal
group of odd dimension, then $S\in\mathcal{Y}_1$, while $L\in\mathcal{Y}_2$.
Suppose that $S=O_{2n}^+(u)$, that is $S$ has the same type and the same Lie rank $n$ as $L$.
Since $n-4$, $n-8$ are multiples of $4$, we can choose a number $a$ of them with $a\equiv4\pmod8$.
Furthermore, $\eta(n-4)>\eta(n-8)=n/2-4>n/3$ because $t(L)\geqslant23$ and so $n\geqslant31$.
Thus $a\in M(L)\cap N(L)$,  and applying Lemma~\ref{l:kdividesj},
we derive a contradiction.

The remaining case for $L=O_{2n}^+(q)$ with $n\equiv0\pmod4$, is $S=O_{2(n-1)}^\varepsilon(u)$,
and this case requires the most effort. Let $s=r_{n-4}(q)$ and $a=n-1-\eta(n-4)=(n+2)/2$.
Similarly to the proof of Lemma~\ref{l:rsnotprs} we conclude that
$r_1=r_{a}(q)$ and $r_2=r_{2a}(q)$ have the following properties: for $i=1,2$, the number $r_i$
is large with respect to $L$, $sr_i\in\omega(L)$, and $psr_i\not\in\omega(L)$. Since $n\equiv0\pmod4$,
we conclude that $r_3=r_{n+4}(q)=r_{2(a+1)}(q)$ has the same properties. Now $r_1,r_2,r_3$ constitute
a coclique in $GK(L)$ because for any $i,j\in\{1,2,3\}$ we have $\varphi(r_i,L)>n/2$ and $e(r_i,q)\neq
e(r_j,q)$ whenever $i\neq j$. Therefore these numbers are large with respect to $S$ and constitute a
coclique in $GK(S)$. Since $\zeta_S$ is injective and $t(s,S)=t(s,L)=t-3$,
the value of $e(s,u)$ depends only on $e(s,q)$. If
$S=O_{2(n-1)}^+(u)$, then $j=e(s,u)=(n-6)/2$, and if $S=O_{2(n-1)}^-(u)$, then $j=e(s,u)=n-6$. In both
cases, $\eta(e(s,u))=\eta(j)$ is odd. Let $r\in\{r_1,r_2,r_3\}$ and write $i=e(r,u)$. As
$rs\in\omega(S)$, we have $\eta(i)\leqslant m-\eta(j)$. If $\eta(i)=m-\eta(j)$, then since
$m$ is odd, it follows that $i$ is even. The condition $rs\in\omega(S)$ implies that
$S$ has an element of order $[u^{\eta(j)}+(-1)^j, u^{\eta(i)}+(-1)^i]/2$, which is
not the case because $i$ and $j$ have opposite parity when $S=O_{2m}^+(u)$, and the same parity when
$S=O_{2m}^-(u)$. Thus $\eta(i)<m-\eta(j)$. If $\eta(i)<m-\eta(j)-2$, then $S$
has an element of order $[u^{\eta(j)}+(-1)^j, u^{\eta(i)}+(-1)^i,u^2+1,u+(-1)^{i+1}]$, which
is divisible by $prs$, contrary to the fact that $prs\not\in\omega(L)$. Let $\eta(i)=m-\eta(j)-2$. Then $i$
is even, and hence $S$ has an element of order $[u^{\eta(j)}+(-1)^j,
u^{\eta(i)}+(-1)^i,u^2+1]$, which is divisible by $prs$. Finally, let
$\eta(i)=m-\eta(j)-1$. But then $e(r,u)$ has only two different possible values, and
this is a contradiction since the numbers $e(r_i,u)$ with $i=1,2,3$ must be distinct.
\end{prf}

\begin{lemma}\label{l:tplneq2}
If $t(p,L)=2$, then the conclusion of Theorem~{\em\ref{t:main}} holds.
\end{lemma}

\begin{prf}
If $t(p,L)=2$, then $L=S_{2n}(q), O_{2n+1}(q)$ and $n$ is even. Since $S$ is symplectic or orthogonal,
we have $k=e(p,u)\in\{1,2\}$.

Suppose $n\equiv2\pmod4$. Then $L\in\mathcal{Y}_2$.  The equality $t(L)=t(S)$ implies that
$S=S_{2n}(u), O_{2n+1}(u)$. Let $j=n/2-2$ and $s=r_{j}(q)$. Observe that
$\eta(j)=j$ is odd. As in the proof of Lemma~\ref{l:rsnotprs}, we write
$a=n-\eta(j)=n/2+2$. Then $r_1=r_a(q)$ and $r_2=r_{2a}(q)$ have the following
properties: for $i=1,2$, the number $r_i$ is large with respect to $L$, $sr_i\in\omega(L)$, and
$psr_i\not\in\omega(L)$. Furthermore, $r_1$ and $r_2$ are not adjacent in $GK(L)$. Since
$t(s,S)=t(s,L)$ and $S$ has the same type and the same Lie rank as $L$, it follows that
$\eta(e(s,u))=\eta(e(s,q))=j$. Let $r\in\{r_1,r_2\}$ and write $i=e(r,u)$. As
$rs\in\omega(S)$, we have $\eta(i)\leqslant n-\eta(j)$. If $\eta(i)<n-\eta(j)$, then $S$ contains
elements of orders $[u^{\eta(j)}+(-1)^j, u^{\eta(i)}+(-1)^i,u-1]$ and $[u^{\eta(j)}+(-1)^j,
u^{\eta(i)}+(-1)^i,u+1]$, one of which must be divisible by $prs$ because $p$ divides either $u-1$
or $u+1$. But then $prs\in\omega(S)\setminus\omega(L)$, a contradiction. Thus
$\eta(i)=n-\eta(j)=a$. Since $a$ is odd, there are two possible values for $e(r,u)$, and these are $a$
and $2a$. As $e(r_1,u)\neq e(r_2,u)$, we have $\{e(r_1,u),e(r_2,u)\}=\{a,2a\}$. If $k=1$, then
$a\equiv k\pmod{2k}$, and if $k=2$, then $2a\equiv k\pmod{2k}$. We denote by $r$ the number
in the set  $\{r_1, r_2\}$ for which $i=e(r,u)\equiv k\pmod{2k}$. Now $s$ and $r$ satisfy the
hypothesis of Lemma~\ref{l:kdividesj}, which is a contradiction because $i\equiv k\pmod{2k}$.

Now suppose that $n\equiv0\pmod4$. The conditions $t(L)=t(S)$ and  $L\in\mathcal{Y}_1$ leave us with two
further possibilities: either $S=S_{2n}(u), O_{2n+1}(u)$ or $S=O_{2m}^\varepsilon(u)$ with $m=n+1$. In the
first case, we use the argument of the previous paragraph with $j=n/2-1$ and $a=n-\eta(j)=n/2+1$.  Now assume that $S=O_{2(n+1)}^\varepsilon(u)$ and write $b=n/2-1$, $s=r_{b}(q)$,
$a=n-b=n/2+1$, $r_1=r_a(q)$, and $r_2=r_{2a}(q)$. Since $\eta(b)=b$ is odd, it follows that $r_1$ and $r_2$ satisfy the conclusion of Lemma~\ref{l:rsnotprs}, that is  $r_1$ and $r_2$ are not adjacent in $GK(L)$, they are large with respect to $L$, and for $i=1,2$ the number $r_is$ belongs to $\omega(L)$, but $pr_is$ does not. Let $r\in\{r_1,r_2\}$, $i=e(r,u)$, and $j=e(s,u)$. Then the equality $t(s,S)=t(s,L)=t-2$ yields
$j=m-3$ when $S=O_{2m}^+(u)$, and $j=(m-3)/2$ when $S=O_{2m}^-(u)$. As
$rs\in\omega(S)$, it follows that $\eta(i)\leqslant m-\eta(j)$. If $\eta(i)< m-\eta(j)$, then $\omega(S)$ contains
$[u^{\eta(j)}+(-1)^j, u^{\eta(i)}+(-1)^i, u+(-1)^{i-1}]$ which is divisible by the odd part of $u^2-1$,
and so by $prs$, which is a contradiction. Thus since $m\equiv1\pmod4$,
the number $\eta(i)=m-\eta(j)=m-(m-3)/2=(m+3)/2$ is even, and therefore there is only one
possible value for $e(r,u)$, namely, $m+3$. But then $e(r_1,u)=e(r_2,u)$, a contradiction.
\end{prf}

\begin{lemma}\label{l:tplneq4} If $t(p,L)=4$, then the conclusion of Theorem~{\em\ref{t:main}} holds.
\end{lemma}

\begin{prf} Note that the equality $t(p,L)=4$ holds if and only if $L=O_{2n}^-(q)$ and $n$ is even.

Suppose that $n\equiv0\pmod4$. If $S=S_{2m}(u)$ or $S=O_{2m+1}(u)$, then $t(L)=t(S)$ yields
$m=n$. But then $S\in\mathcal{Y}_1$ and $L\in\mathcal{Y}_3$, a contradiction. If
$S=O_{2m}^-(u)$, then $m\in\{n,n+1\}$. Suppose that $m=n$. Then $k=e(p,u)=4$. Since $a=n-4$ and $b=n-8$
lie in $M(L)\cap N(L)$ and are multiples of $4$, it follows that one of the numbers $e(r_a(q),u)=e(r_a(q),q)=a$
and $e(r_b(q),u)=e(r_b(q),q)=b$ is congruent to $4$ modulo $8$. But this contradicts Lemma~\ref{l:kdividesj}.
If $S=O_{2m}^+(u)$, then $m\in\{n+1,n+2\}$. Let $m=n+2$. Then
$m\equiv2\pmod4$, and, on the other hand, it follows from Table~\ref{tab:KlL} that $m\equiv4\pmod6$. Thus
$m\equiv10\pmod{12}$, and hence $m=12c+10$, $n=12c+8$ for some integer  $c$. Writing $X=\{a\mid
n/2\leqslant \eta(a)< 2n/3\}$ and noting that $n/2=6c+4$ is even,  we conclude that
$|X|=\left[3((8c+6)-(6c+4))/2\right]=3c+3$ by Lemma \ref{l:eta}. There is a coclique $\rho$ of size $|X|$
in $GK(L)$ such that $e(r,q)\in X$ for every $r\in\rho$. Let $r\in\rho$. The prime $r$ is large with respect to $L$,
and so $r$ is large with respect to $L$ too, whence $\varphi(r,S)\geqslant m/2$ (see Table~\ref{tab:tL}). Furthermore,
there is $b\in M(L)=N(L)$ such that $rr_b(q)\in\omega(L)$. Since $e(r_b(q),u)\in
N(S)$ and $N(S)=M(S)$, we have $\varphi(r,S)\leqslant m-\varphi(r_b(q),S)<2m/3$. Thus writing
$Y=\{a\mid m/2\leqslant\eta(a)<2m/3\}$, we see that $e(r,u)\in Y$ for every $r\in\rho$. By Lemma
\ref{l:eta}, the size of $Y$ is equal to $\left[(3((8c+7)-(6c+5))+1)/2\right]=3c+3=|X|$, and therefore
$\{e(r,u)\mid r\in\rho\}=Y$. It follows that $\rho$ contains $r$ and $r'$ such that $e(r,u)=m/2$ and
$e(r',u)=m$. But then $rr'\in\omega(S)$, a contradiction.

It remains to consider the situation where $S=O_{2m}^\varepsilon(u)$, $m=n+1$, and either
$k=e(p,u)=4$ or $m\equiv5\pmod{12}$ and $(\varepsilon, k)\in\{(+,6),(-,3)\}$. First, assume that
$m\equiv5\pmod{12}$ and $k\neq4$. Observe that the congruence and the condition $t\geqslant23$
imply that $m=n+1\geqslant45$. In particular,  we have $(n+12)/2>n/3$, and so
$n-4, n-8, n-12\in M(L)\cap N(L)$ with $t(r_{n-4}(q),L)=t-2$, $t(r_{n-8}(q),L)=t-5$,
and $t(r_{n-12}(q),L)=t-8$. If $S=O_{2m}^+(u)$, then it follows that the numbers $e(r_{n-4}(q),u)=m-3$,
$e(r_{n-8}(q),u)=m-7$, and $e(r_{n-4}(q),u)=m-11$ are congruent to $2$ modulo $4$. Moreover, they
have different residues modulo $3$, and hence one of them is congruent to $6$ modulo~$12$.
We denote this number by $j$ and take $a$ to be such that $e(r_a(q),u)=j$. Applying
Lemma~\ref{l:kdividesj} to $a\in M(L)\cap N(L)$, we conclude that $j\not\equiv6\pmod{12}$, a contradiction. If $S=O_{2m}^-(u)$, then the numbers $e(r_{n-4}(q),u)=(m-3)/2$, $e(r_{n-8}(q),u)=(m-7)/2$,
and $e(r_{n-4}(q),u)=(m-11)/2$ are odd and one of them is a multiple of $3$. Since $k=3$, we again derive
a contradiction by Lemma~\ref{l:kdividesj}. Thus we may assume that $k=e(p,u)=4$.
Let $s_1=r_{(n-2)/2}(q)$, $s_2=r_{n-2}(q)$, $r_1=r_{n+2}(q)$, $r_2=r_{(n+2)/2}$, and
$w=r_{n}(q)$. It follows that $(n-2)/2,n-2\in M(L)\cap N(L)$, the primes $r_1,r_2,w$ are large with respect to $L$
and for $l=1,2$, we have $s_lr_l,s_lw\in\omega(L)$, but $ps_lr_l,ps_lw\not\in\omega(L)$. Furthermore,
$t(s_1,L)=t(s_2,L)=t-1=t(s_1,S)=t(s_2,S)$. Therefore, $j=e(s_1,u)=e(s_2,u)$ is equal to
$(m-3)/2$ when $S=O_{2m}^+(u)$ and to $m-3$ when $S=O_{2m}^-(u)$. Let $r\in\{r_1,r_2,w\}$ and
$i=e(r,u)$. As for $l=1,2$, we have $s_lr_l,s_lw\in\omega(S)$, it follows that $\eta(i)\leqslant
m-\eta(m-3)=m-\eta((m-3)/2)=(m+3)/2$. If $\eta(i)=(m+3)/2$, then $\eta(i)$ is even, and so
$i=m+3$. But then $i+2\eta(j)=2\eta(i)+2\eta(j)=2m$ and either $j=\eta(j)$ and $S=O_{2m}^+(u)$ or
$j=2\eta(j)$ and $S=O_{2m}^-(u)$, and by adjacency criterion, $r$ and $s_l$ are not adjacent in
$S$ for $l=1,2$, which is not the case. Thus $\eta(i)<(m+3)/2$. If $\eta(i)=(m-1)/2$, then
$O_{2m}^+(u)$ has an element of order $a=[u^{\eta(i)}+1,u^j-1,u^2+1]$ and
$O_{2m}^-(u)$ has an element of order $b=[u^{\eta(i)}+1,u^j+1,u^2+1]$. This is a contradiction because $prs_l$ divides both $a$ and $b$
for $l=1,2$. If $\eta(i)<(m-1)/2$, then it is easy to construct an element of required order. So $\eta(i)=(m+1)/2$ and there are only two possible values for $e(r,u)$.
This is impossible because $\{r_1,r_2,w\}$ is a coclique in~$GK(L)$.

Now suppose that $n\equiv2\pmod4$. The group $S$ cannot be isomorphic to $O_{2m}^+(u)$ since otherwise
$t(L)$ is congruent to $2$ modulo $3$ and $t(S)$ does not. Let $S\in\{S_{2m}(u),O_{2m+1}(u)\}$.
Then $t(L)=t(S)$ yields $m\in\{n,n-1\}$. If $m=n$, then $S\in\mathcal{Y}_2$, and if
$m=n-1$, then $S\in\mathcal{Y}_3$, which contradicts to the fact that $L\in\mathcal{Y}_1$. Thus
$S=O_{2m}^-(u)$, and it follows from $t(L)=t(S)$ that $m=n$. Then $m\equiv2\pmod4$, and so
$k=e(p,u)=4$. Assume that $t>23$, that is $m=n>30$. Then $\eta(n-6)>\eta(n-10)=n/2-5>n/3$, and hence
$n-6,n-10\in M(L)\cap N(L)$ with one of them being a multiple of $4$ but not of $8$. Now Lemma~\ref{l:kdividesj}
gives the desired contradiction.

Finally, assume that $L=O_{60}^-(q)$ and $S=O_{60}^-(u)$. Let $s_1=r_{11}(q)$, $s_2=r_{22}(q)$,
$r_1=r_{38}(q)$, $r_2=r_{19}(q)$, and $w=r_{18}(q)$. Then $11,22\in M(L)\cap N(L)$,
$\{r_1,r_2,w\}$ is a coclique in $GK(L)$ and it consists of numbers large with respect to $L$, and
for $l=1,2$, we have $s_lr_l,s_lw\in\omega(L)$, but $ps_lr_l,ps_lw\not\in\omega(L)$. Writing
$j_1=e(s_1,u)$, $j_2=e(s_2,u)$ and using Lemma~\ref{l:transfer}, we obtain
$\eta(j_1)=\eta(j_2)=11$. Choose $r\in\{r_1,r_2,w\}$ and take $i=e(r,u)$. Since for $l=1,2$,
we have $s_lr_l,s_lw\in\omega(S)$, it follows that $\eta(i)\leqslant19$. It is easy to see that
if $\eta(i)\leqslant16$, then $\omega(S)$ contains a number divisible by $prs_l$ for each
$l=1,2$, which is impossible. So we assume that $\eta(i)\geqslant17$. If $\eta(i)=18$, then $i=36\equiv4\pmod8$,
and applying Lemma~\ref{l:kdividesj} either to the pair $(s_l,r_l)$ or to the pair $(s_l,w)$ for some
$l\in\{1,2\}$, we have the desired contradiction. Let $\eta(i)\in\{17,19\}$ and
$j\in\{j_1,j_2\}$. If $j=11$, then since $S$ has no elements of order $[u^{11}-1, u^{19}-1]/d$ with
$d\in\{1,2,4\}$ and does have element of order $[u^2+1,u^{11}-1,u^{17}-1]$, it follows that $i\in\{34,38\}$.
Similarly, if $j=22$ then $i\in\{17,19\}$. If $j_1=j_2$, then $i=e(r,u)$ can take at most two
different values, but there must be at least three of them as $\{r_1,r_2,w\}$ is a coclique in $GK(L)$.
And if $j_1\neq j_2$, then $e(w,u)\in\{17,19\}\cap\{34,38\}=\varnothing$, and this contradiction
completes the proof.
\end{prf}

If $L=L_n^\varepsilon(q)$ with $n\geqslant45$,  $L\in\{S_{2n}(q),O_{2n+1}(q)\}$ with $n\geqslant29$,
$L=O_{2n}^-(q)$ with $n\geqslant30$, or $L=O_{2n}^+(q)$ with $n\geqslant31$, then
$t(L)\geqslant23$. We eliminated all possibilities for the group $S$ under the assumption that $t(L)\geqslant23$,
and so the proof of Theorem~\ref{t:main} is complete.

\begin{rem} Although the final part of the proof is technically complicated, its idea is transparent and
based on an application of the well-known pigeonhole principle. Here pigeons are from $M(L)\cap
N(L)$ and holes are from $N(S)$. The number of pigeons and number of holes are almost the same and,
what is sufficient, the difference between them is a small constant which does not depend on
$n=\prk(L)$. However, there are prohibited holes, and the number of those holes increases when so
does~$n$. It provides a contradiction that becomes evident with the growth of~$n$.
\end{rem}

\section{Proof of Theorems~1 and~2}

As written in Introduction, Theorems~\ref{t:L&U} and~\ref{t:S&O} follow from Theorem~\ref{t:main}
and a series of previously obtained results. The next assertion summarizes the main results
of~\cite{VasGrMaz1,VasGrSt}.

\begin{lemma}\label{l:othercases}
Let $L$~be a finite simple classical group over a field of characteristic $p$, and
$L\not\in\{L_2(9), L_3(3), U_3(3), U_3(5), U_5(2), S_4(3)\}$. Suppose that $G$ is a finite group
with $\omega(G)=\omega(L)$, and $S$ is a unique nonabelian composition factor of~$G$. Then one of
the following holds:

\emph{(1)} $S\simeq L;$

\emph{(2)} $L=S_4(q)$, where $q>3$, and $S\simeq L_2(q^2);$

\emph{(3)} $L\in\{S_6(q), O_7(q), O_8^+(q)\}$ and $S\in\{L_2(q^3), G_2(q), S_6(q), O_7(q)\};$

\emph{(4)} $n\geqslant4$, $L\in\{S_{2n}(q), O_{2n+1}(q)\}$ and $S\in\{O_{2n+1}(q), O_{2n}^-(q)\};$

\emph{(5)} $n\geqslant6$ is even, $L=O_{2n}^+(q)$ and $S\in\{S_{2n-2}(q), O_{2n-1}(q)\};$

\emph{(6)} $S$ is a group of Lie type over a field of characteristic distinct to~$p$.
\end{lemma}

This assertion and Theorem~\ref{t:main} yield Theorem~\ref{t:S&O} immediately.

\smallskip

A \emph{cover} of a finite group $G$ is a finite group having $G$ as a factor group. A cover is
proper if the corresponding factor group is proper. If $\omega(G)\neq\omega(H)$ for every proper
cover $H$ of $G$, then $G$ is said to be \emph{recognizable by spectrum among covers}. It follows
from~\cite[Lemma~9]{ZavLU} that a finite group $G$ is recognizable by spectrum among covers if and
only if $\omega(G)\neq\omega(H)$ for every splitting extension $H=V\rtimes G$, where $V$ is
an absolutely irreducible finite-dimensional $G$-module over a finite field of characteristic~$r$.
For every simple linear and unitary group $L=L_n^\varepsilon(q)$ with $n\neq 4$ and every
$L$-module $V$, if the characteristic $r$ of $V$ coincides with the characteristic $p$ of $L$, then
it is proved that $\omega(L)\neq\omega(V\rtimes L)$ \cite{ZavLU}. If $V$ is defined over
the field of characteristic distinct to $p$ and $L=L_n(q)$, then any extension $V\rtimes L$
contains an element whose order does not belong to~$\omega(L)$ \cite[Lemma~11]{ZavLU}. For other
classical groups the following general result was recently obtained.

\begin{lemma}[\cite{GrCovCl}]\label{l:GrCov}
Let $L$ be one of the simple groups $U_n(q)$, where $n\geqslant 4$, $S_{2n}(q)$, where $n\geqslant
3$, $O_{2n+1}(q)$, where $n\geqslant 3$, and $O^\pm_{2n}(q)$, where $n\geqslant 4$. Suppose that
$V$ is a finite-dimensional $L$-module over a field of characteristic $r$ prime to $q$. Then either
$\omega(V\rtimes L)\neq \omega(L)$ or $L=U_5(2)$ and $r=3$. If $L=U_5(2)$, then there is a
$10$-dimensional $L$-module $V$ over a field of characteristic $3$ such that
$\omega(V\rtimes L)= \omega(L)$.
\end{lemma}

If $L=L_n^\varepsilon(q)$, $n\geqslant45$, and $G$ is a finite group isospectral to~$L$, then
Theorem~\ref{t:main} and Lemma~\ref{l:othercases} provide
$L\leqslant G/K\leqslant\operatorname{Aut}L$, where $K$ is the soluble radical of~$G$.
An application of the aforementioned results on covers of linear and unitary groups completes the
proof of Theorem~\ref{t:L&U}.

\section*{Acknowledgments}

The author is grateful to M.A. Grechkoseeva for her reading earlier drafts of this article. Her valuable comments and helpful suggestions have enabled
us to improve the text significantly.

%\bibliographystyle{elsarticle-num}
%\bibliography{fsg}
\end{document}